%% file: octut.tex
\documentclass[12pt,a4paper]{article}\usepackage[top=20mm,left=12mm]{geometry}\textwidth18.2cm\textheight25.1cm
\usepackage{amsmath,amssymb,amsthm}
\usepackage{float}\floatstyle{plaintop}\restylefloat{table}
\usepackage[tableposition=top]{caption}
\usepackage{graphicx,pstricks,fancyvrb,multirow,lineno,url}
\usepackage{listings}\usepackage{paralist, longtable}
\usepackage{setspace}
\newtheoremstyle{myplain}
  {-\baselineskip\topsep}   
  {\topsep}   
  {\itshape\setstretch{1.05}}  
  {0pt}       
  {\bfseries} 
  {.}         
  {5pt plus 1pt minus 1pt} 
  {}       
\theoremstyle{myplain}
\newtheorem{theorem}{Theorem}[section]
\newtheorem{remark}[theorem]{Remark}
\input{./hudef.tex}\input{./p2pdef.tex}

\newlength{\tew}
\newcommand{\argmax}{\operatornamewithlimits{argmax}}
\def\ig{\includegraphics}

\renewcommand{\arraystretch}{1}\renewcommand{\baselinestretch}{1.05}
\def\medskip{}\def\bigskip{}
\usepackage{hyperref}

\newfloat{Algorithm}{ht}{lop}
\def\alv{\vec\al}
\def\ug{u_{\text{guess}}}\def\lamh{\hat\lam}\def\uend{\uh_0}
\def\epsi{\eps_\infty}

\input{lst}

\def\Dx{\dd x}\def\evalat#1{#1}
\def\medskip{}\def\bigskip{}\def\k{q}
\def\taskip{\renewcommand{\arraystretch}{1}\renewcommand{\baselinestretch}{1}}
\def\teskip{\renewcommand{\arraystretch}{1.1}
\renewcommand{\baselinestretch}{1.2}}
\def\hulst#1#2{\taskip\lstinputlisting[#1]{#2}\teskip}
\def\hutab#1{\taskip\begin{\table}#1\end{table}\teskip}

\usepackage{hyperref}
\allowdisplaybreaks

\begin{document}
\mbox{}\vspace{-0mm}\begin{center}\Large
Infinite time horizon 
spatially distributed optimal control problems with \pdep\ -- 
algorithms and tutorial examples\\[2mm]
\normalsize Hannes Uecker$^1$, Hannes de Witt$^2$\\[2mm]
\footnotesize
$^1$ Institut f\"ur Mathematik, Universit\"at Oldenburg, D26111 Oldenburg, 
hannes.uecker@uni-oldenburg.de\\
$^2$ Institut f\"ur Mathematik, Universit\"at Oldenburg, D26111 Oldenburg, 
hannes.de.witt1@uni-oldenburg.de\\[2mm]
\normalsize
\today
\end{center}
\noi
\begin{abstract} 
We use the continuation and bifurcation package \pdep\  
to numerically analyze infinite time horizon optimal control  
problems for parabolic systems of PDEs. The basic idea is  a
two step approach to the canonical systems, derived from \PMAXP. 
First we find branches of steady or time--periodic states of the 
canonical systems, i.e., canonical steady states (CSS) 
respectively canonical periodic states (CPS), and then 
use these results 
to compute 
time-dependent canonical paths 
connecting to a CSS or a CPS with the so called saddle point property. 
This is a (high dimensional) boundary value problem in time, 
which we solve by a continuation algorithm in the initial states. 
We first explain the algorithms and then the implementation 
via some example problems and associated 
\pdep\ demo directories. The first two examples deal 
with the optimal management of a distributed shallow lake, and 
of a vegetation system, both with (spatially, and temporally) 
distributed controls. These examples show interesting 
bifurcations of so called patterned CSS, including 
patterned {\em optimal} steady states. As a third example we discuss 
optimal boundary control of a fishing problem with boundary catch. 
For the case of CPS--targets we first focus on an ODE toy model to 
explain and validate the method, and then discuss an 
optimal pollution mitigation PDE model. 
\end{abstract}
\noindent

\tableofcontents 

\section{Introduction}\label{i-sec}
We consider optimal control (OC) problems for partial differential 
equations (PDEs) of the form 
\begin{subequations} \label{oc1}
  \hual{
&\pa_t v=-G_1(v,\k):=D\Delta v+g_1(v,\k), \text{ in $\Om\subset\R^d$ (a bounded domain)}, 
}
with initial condition $v|_{t=0}=v_0$, and suitable boundary conditions. 
Here $v:\Om\times[0,\infty)\ra\R^N$ denotes a (vector of) state variables, $\k:\Om\times[0,\infty)\ra\R$ is a (distributed) control, 
$D\in\R^{N\times N}$ is a diffusion matrix, and 
$\Delta=\pa_{x_1}^2+\ldots+\pa_{x_d}^2$ is the Laplacian. 
The goal is to find 
\hual{
V(v_0)&:=\max_{\k(\cdot,\cdot)}J(v_0(\cdot),\k(\cdot,\cdot))
} 
\end{subequations}
for the discounted time integral
\hual{\label{jdef}
& J(v_0(\cdot),\k(\cdot,\cdot)):=\int_0^\infty\er^{-\rho t}
J_{ca}(v(t),\k(t)) \dd t,
}
where $J_{ca}(v(\cdot,t),\k(\cdot,t))=\frac 1{|\Om|}
\int_\Om J_c(v(x,t),\k(x,t))\dd x$ 
is the spatially averaged current value objective function, 
with the local current value $J_c:\R^{N+1}\ra\R$ a given function, and 
$\rho>0$ is the discount rate which corresponds to a long-term 
investment rate.  The discounted time integral $J$ 
is typical for economic
problems, where ``profits now'' weight more than mid or far 
future profits. 
The $\max$ (formally $\sup$) in (\ref{oc1}b) runs over all {\em admissible} controls $\k$; this will be specified in more detail in the examples below. 
Additionally, we  also give one example of a boundary control, where 
$\k:\Om_\k\subset\pa\Om\to \R$, and where \reff{oc1} and \reff{jdef} 
are modified accordingly. 
In applications, $J_c$ and $G_1$ of course often also 
depend on a number of parameters, 
which however for simplicity we do not display 
here.\footnote{$G_1$ in (\ref{oc1}a) can in fact be of a much more general form, 
but for simplicity here we stick to  (\ref{oc1}a).}

In this tutorial we explain by means of four examples  how to numerically 
study problems of type \reff{oc1} with \pdep\footnote{see \cite{p2pure} for 
background, and \cite{p2phome}  
for download of the package, demo files, and various documentation and tutorials, including a quick start guide also giving installation 
instructions}. The examples are from \cite{U16, GU17, GUU19, hotheo}, and 
we mostly refer to these works and the references therein 
for modeling background and (bioeconomic) interpretation of the results, and for general references and comments on OC for PDE problems, here keeping these aspects to a minimum. 
%
%
In the first example, $\k:\Om\times[0,\infty)\ra\R$ is the phosphate load in 
a model describing the phosphorus contamination of a shallow lake 
by a scalar PDE (\ref{oc1}a) with homogeneous Neumann BCs 
$\pa_\nu v=0$, $\nu$ the outer normal.
Similarly, in the second example,  $\k:\Om\times[0,\infty)\ra\R$ is 
the harvesting effort
in a vegetation-water system, such that (\ref{oc1}a) is a 
two component reaction diffusion, again with homogeneous Neumann BCs, while 
in the third example $\k:\pa\Om_\k\times[0,\infty)\ra\R^2$ 
with $\Om_\k\subset\pa\Om$ is a boundary control, namely the fishing effort on (part of) the shore of a lake. In these examples, so--called 
canonical steady states (CSSs), i.e., steady states of the 
so--called canonical system (see below) play an important role. 
The fourth example considers optimal pollution (mitigation), where the 
states are the emissions of some firms and 
the pollution stock, and the control is the pollution abatement investment. 
Here, canonical periodic states (CPSs) play an important role. 
To keep this tutorial simple, for all examples we restrict to
the 1D case of $\Om\subset\R$ an interval. 
Generalizations to domains in $\R^2$ are
straightforward, and also straightforward to implement in \pdep, and
for the shallow lake and the vegetation systems 
have also been studied in \cite{U16,GU17} respectively.

In the remainder of this introduction, we first briefly review the 
derivation of the canonical system as a necessary first order 
optimality condition from \reff{oc1} via 
Pontryagin's maximum principle. Then we explain the basic algorithms to 
treat the canonical system, i.e., to 
first find CSSs and CPSs and then canonical paths (CPs), i.e., 
solutions of the canonical system which connect some given 
initial states to some CSS or CPS with the so--called 
saddle point property. These CPs yield candidates for 
solutions of \reff{oc1}.  
  In \S\ref{ex-sec} we 
present the example problems and the \pdep\ implementation details 
for controlling to a  CSS, and in \S\ref{cps-sec} 
for controlling to a CPS. In \S\ref{dsec} we close with 
a brief summary and outlook. 
We assume that the reader has a basic knowledge of 
\mlab, has installed \pdep\ (see also \cite{qsrc}), 
and also has some previous experience with 
the software. If this is not the case, then we recommend to at least 
briefly look at one of the simpler problems discussed in, 
e.g., \cite{actut}. 

\brem\label{legrem}{\rm The \pdep\ library described here is in 
{\tt libs/oclib}. There is also an older OC lib, namely {\tt libs/oc}, 
which we mainly keep for downward compatibility, and the associated 
old demos are in {\tt ocdemos/legacy}. The upgrade of {\tt oc} to 
{\tt oclib} includes the computation of CPs to CPSs, and the option 
to free the truncation time $T$, and we strongly recommend to switch to 
this setting.}\eex\erem 

\subsection{Pontryagin's Maximum Principle and the canonical system}
We first consider the case of spatially distributed 
controls $\k:\Om\times[0,\infty)\ra\R$, and assume homogeneous Neumann 
BCs for $v$, i.e., $\pa_\nu v=0$ on $\pa\Om$, $\nu$ the outer normal, 
and introduce the costates $\lam:\Om\times(0,\infty)\ra \R^{N}$ and the 
(local current value) Hamiltonian
\hual{
\CH&=\CH(v,\lam,\k)=J_c(v,\k)+\lam^T(D\Delta v+g_1(v,\k)).} 
By Pontryagin's Maximum Principle (see Remarks \ref{lrem} and \ref{lagrem}) 
for the intertemporal Hamiltonian $\tilde{\CH}=
\int_0^\infty \er^{-\rho t} \ov{\CH}(t)\dd t$ with the spatial integral 
\huga{\label{fullH} 
\ov{\CH}(t)=\int_\Om \CH(v(x,t),\lam(x,t),\k(x,t))
\dd x, 
}
an optimal solution $(v,\lam)$ (or equivalently 
$(v,\k):\Om\times[0,\infty)\ra \R^{N+1})$) has to solve the canonical system (CS) 
\begin{subequations} \label{cs}
\hual{
\pa_t v&=\pa_\lam\CH=D\Delta v+g_1(v,\k), \quad v|_{t=0}=v_0, \\
\pa_t \lam&=\rho\lam-\pa_v\CH=\rho\lam+g_2(v,\k)-D\Delta\lam, 
\intertext{
where  $\k=\argmax_{\tilde{\k}}\CH(v,\lam,\tilde{\k})$, which generally we assume to 
be obtained from solving}
\pa_\k\CH(v,\lam,\k)&=0.
\intertext{Under suitable concavity conditions on $J_c$ this holds due to the absence 
of control constraints.  
The costates $\lam$ also fulfill zero flux BCs, and in the derivation 
of \reff{cs} we imposed the so called inter--temporal transversality condition}
\lim_{t\to\infty}&\er^{-\rho t}\int_\Om \spr{v(t,x),\lam(t,x)}\dd x=0. \label{tcond}
}
\end{subequations}

In principle we want so solve \reff{cs} for $t\in[0,\infty)$, but 
in (\ref{cs}a) we have initial data for only half the variables, 
and in (\ref{cs}b) we have anti--diffusion, such that \reff{cs} 
is ill--posed as an initial value problem. 
For convenience we set%
\footnote{the notation $u=(v,\lam)$ for the vector of state and 
costate variables is used here 
as $u$ generally denotes the vector of unknowns in \pdep; in optimal 
control $u$ is often used as the notation for the control, which here we 
denote by $q$;} 
\huga{
u(t,\cdot):=\bpm v(t,\cdot)\\ \lam(t,\cdot)\epm: \Om\ra\R^{2N}, 
} 
and write \reff{cs} as 
 \hual{\label{cs2}
\pa_t u&=-G(u,\eta):=\CD\Delta u+f(u,\eta), \quad 
\CD=\bpm D&0\\0&-D\epm, \quad f(u,\eta)=\bpm g_1(u,\eta)\\ g_2(u,\eta)\epm, 
}
with BCs $\pa_\nu u=0$, where $\eta\in\R^p$ stands for parameters present,
 which for instance include the discount rate $\rho$. 
A solution $u$ of the canonical system \reff{cs2} is called a 
\emph{canonical path} (CP), 
a fixed point of \reff{cs2} (which automatically fulfills \reff{tcond}) 
is called a 
\emph{canonical steady state (CSS)} and a time-periodic solution of \reff{cs2} is called \emph{canonical periodic states (CPS)}. With a slight abuse of notation 
we also call $(v,\k)$ with $\k$ given by (\ref{cs}c) a canonical 
path.  

\brem\label{lrem}{\rm 
For general background on OC in an ODE setting with a focus on the 
infinite time horizon see \cite{grassetal2008} or \cite{T15}. 
For the PDE see, \cite{Tr10} and the 
references therein, or specifically \cite{RZ99, RZ99b} and 
\cite[Chapter5]{AAC11} for \PMAXP{} 
for OC problems for semi-linear diffusive models. However, these 
works are in a 
finite time horizon setting, and often the objective function 
is linear in the control and there are control constraints, e.g., 
$\k(x,t)\in Q$ with some bounded set $Q$. Therefore $\k$ is not obtained from 
the analogue of (\ref{cs}c), but rather takes the values from $\pa Q$, which 
is often called bang--bang control.  
Here we do not (yet) consider (active) 
control or state constraints, and no terminal time, but the infinite 
time horizon. 
Our distributed OC models are motivated by \cite{BX08, BX10}, which also 
discuss \PMAXP{} in this setting. 
}\eex \erem 

\brem\label{lagrem}{\rm
The use of the Hamiltonian $\tilde\CH$ is the standard way of dealing 
with intertemporal OC problems in economics. Equivalently, the canonical system \reff{cs} 
is formally obtained as the first variation of the Lagrangian 
\huga{\label{L1}
\CL=\frac 1{|\Om|}\int_0^\infty \er^{-\rho t}\left(\int_\Om J_c(v,q)
-\spr{\lam, \pa_t v+G_1(v,q)}\dd x\right) \dd t, 
}
where $G_1(v,q)=-D\Delta v-g_1$, and where $\lam=(\lam_1,\ldots,\lam_N)$ can 
be identified as Lagrange multipliers to the constraint (\ref{oc1}a), i.e., 
$\pa_t v+G_1(v,q)=0$. Using integration by parts in $x$ with the Neumann BCs 
$\pa_n v=0$ and $\pa_n \lam=0$ we have 
$
\int_\Om \spr{\lam, D\Delta v}\dd x=\int_\Om \spr{D\Delta\lam, v}\dd x, 
$
and using integration by parts in $t$ with transversality condition 
\reff{tcond} yields 
$
-\int_0^\infty \er^{-\rho t}\int_\Om \spr{\lam,\pa_t v}\dd x\dd t
=\int_\Om \spr{\lam(x,0),v(x,0)}\dd x+\int_0^\infty \er^{-\rho t}\spr{\pa_t \lam-\rho\lam, v}\dd x\dd t. 
$
Thus, $\CL$ can also be written as 
\hual{
\CL=&\frac 1{|\Om|}\biggl[\int_\Om \spr{\lam(x,0),v(x,0)}\dd x\label{L2}\\
&+\int_0^\infty \er^{-\rho t}\left(\int_\Om J_c(v,q)
{+}\spr{\pa_t\lam{+}\rho\lam{+}D\Delta\lam,v}{+}\spr{\lam,g_1(v,q)}
\dd x\right) \dd t\biggr],\notag 
}
and \reff{cs} are the first variations of $\CL$ with respect to $\lam$ 
(using \reff{L1}) and  $v$ (using \reff{L2}) with $v(0,x){=}v_0(x)$. 
Both computations (with $\tilde{\CH}$ and $\CL$) are somewhat 
formal, 
and in particular the necessity of the transversality condition \reff{tcond},  
is subject to active research.  See also \cite{GUU19} and the 
references therein for a discussion of rigorous results for 
infinite time horizon OC problems with PDE constraints. 
}
\eex 
\erem 

\subsection{The general setup and the algorithms for 
canonical paths}
To study \reff{cs2} 
we proceed in two steps, which can be seen as a variant of the ``connecting 
orbit method'', see, e.g., \cite{BPS01}, \cite[Chapter 7]{grass2014} 
and Remark \ref{corem}: 
first we compute (branches of) CSSs and CPSs, and 
second we compute CPs connecting some initial states to some CSSs or CPSs.  
Thus we take a somewhat broader perspective than aiming at computing just one optimal control, 
given an initial condition $v_0$. 
Instead, we aim to give a somewhat global 
picture by identifying the optimal CSS/CPS and their respective domains 
of attraction. 

\paragraph{(a) Branches of CSSs and CPSs.} 
We compute (approximate) CSSs of \reff{cs2}, i.e., solutions $\uh$ of 
\huga{\label{cs21}
G(u,\eta)=0, 
}
together with the spatial BCs, by discretizing \reff{cs21} via the 
finite element 
method (FEM) and then treating the discretized system as a
 continuation/bifurcation problem.%
\footnote{The~ $\hat{}$\,\, notation is often used in OC for CSS, and is 
not related to Fourier transform in any way; we use the notation $\uh$ 
for CPS $t\mapsto\uh(t)$ in an analogous sense.} 
This gives branches $s\mapsto (\uh(\eta),\eta(s))$ of solutions, parameterized 
by a (pseudo-) arclength, which is in particular 
useful to possibly find 
several solutions $\uh^{(l)}(\eta)$, $j=l,\ldots, m$ at fixed $\eta$. 
CPS are usually not computed directly, but via Hopf bifurcation 
from branches of CSS. Thus, after finding such Hopf bifurcations, 
we do a branch switching with an appropriate initial guess 
for the rescaled problem
\begin{subequations}
	\hual{\label{cps21}
		\pa_t u&=-T_pG(u,\eta), \\
		u(0)&=u(1),}
\end{subequations}
where the period $T_p$ becomes an additional unknown, see \cite{hotheo}.

By computing the associated $J_{ca}(\hat v,\hat\k)$ 
we can identify which of the CSSs and CPSs is optimal amongst the CSSs 
and CPSs. 
Given a CSS $\uh$, for simplicity we also write 
$J_{ca}(\uh):=J_{ca}(\hat v^{(l)},\k^{(l)})$, and moreover have, 
by explicit evaluation of the time integral, 
\huga{J(\uh)=J_{ca}(\uh)/\rho.}
For a CPS $\uh$ with period length $T_p$ we 
evaluate the time integral
\huga{
J(\uh)=\int_0^{\infty} e^{-\rho t} J_{ca}(\uh(t)) dt
=\frac{1}{1-e^{-\rho T_p}}  
\int_0^{T_p} e^{-\rho T_p} J_{ca}(\uh(t)) dt.
}
Due to the discounting this integral may highly depend on the 
phase of the CPS, and thus we do not have a single objective value 
for a CPS, but a continuum of objective values. However, 
when computing CPs to a CPS it generally turns out that the 
values of the CPs are independent of the chosen phase of the CPS, 
see Remark \ref{prem1} below.

\paragraph{(b) Canonical paths to canonical steady states.}
In a second step, using the results from (a), we compute CPs 
connecting chosen initial states to a CSS $\uh$ (or a CPS $\uh$, 
see below), and the objective 
values of the canonical paths. For paths to a CSS we choose a truncation time 
$T$ and modify \reff{tcond} to the condition 
that $u(T)\in W_s(\uh)$ and near $\uh$, where $W_s(\uh)$ denotes the stable 
manifold of $\uh$. In practice, we approximate $W_s(\uh)$ by the stable 
eigenspace $E_s(\uh)$, and thus consider the time-rescaled BVP 
\begin{subequations}\label{bvp0}
\hual{
\pa_t u&=-TG(u),\\
v|_{t=0}&=v_0,\\
u(1)&\in E_s(\uh),  
} 
\end{subequations}
and $\|u(1)-\uh\|$ small in a suitable sense, further discussed below. 
If the mesh in the FEM discretization from (a) consists of $n$ nodes, then 
$u(t)\in\R^{2Nn}$, and (\ref{bvp0}a) yields a system of $2Nn$ ODEs 
in the form (with a slight abuse of notation) 
\begin{subequations}\label{bvp1}
\hual{
M\ddt u&=-TG(u),
}
while the initial and transversality conditions become%
\footnote{recall that we put the 
truncation time $T$ into $\pa_t u=-TG(u)$ 
such that the end point of a CP is at $t=1$}
\hual{
v|_{t=0}&=v_0,\\
\Psi(u(1)-\uh)&=0.  
} 
\end{subequations}
Here, $M\in\R^{2Nn\times 2Nn}$ is the mass matrix of the FEM mesh, 
(\ref{bvp1}b) consists of $Nn$ initial conditions for 
the states, while the costates $\lam$ (and hence the control $\k$) 
are free, and 
$\Psi\in \R^{Nn\times 2Nn}$ defines the projection onto the 
unstable eigenspace $E_u(\uh)$, 
where due to the convention that $\pa_t u=-TG(u)$, the stable 
eigenspace is spanned by the (generalized) eigenvectors of $\pa_u G(u)$ 
to eigenvalues $\mu$ with {\em positive} real parts. 
Thus, to have $2Nn$ BCs altogether we need dim$E_s(\uh)=Nn$. 
On the other hand, we always have dim$E_s(\uh)\le Nn$, see 
\cite[Appendix A]{GU17}. We define the defect 
\huga{\label{ddef}
d(\uh):={\rm dim}E_s(\uh)-Nn}
and call a CSS $\uh$ with $d(\uh)=0$ a CSS with the saddle--point 
property (SPP). At first sight it may appear that $d(\uh)$ depends 
on the spatial discretization, i.e., on the number of $n$ of nodes. 
However, $d(\uh)$ remains constant for finer and finer meshes, 
see  \cite[Appendix A]{GU17} for further comments. 

For $\uh=(\vh,\hat\lam)$ with the SPP, 
and $\|v_0-\vh\|$ sufficiently small, we may expect the existence 
of a solution $u$ of \reff{bvp1}, which moreover can be found from 
a Newton loop for \reff{bvp1} with initial guess $u(t)\equiv \uh$. 
Here, as a first guess for the truncation time $T$ we may use 
the longest decay length of the stable directions, i.e., 
\huga{\label{T0}
T=\re(\mu_2)^{-1}, 
}
where $\mu_2$ is the stable eigenvalue with the smallest real part.  

For larger $\|v_0-\vh\|$ a solution of \reff{bvp1} 
may not exist, or a good initial guess may be hard 
to find, and therefore we use a continuation process for \reff{bvp1}. 
In the simplest setting, assume that for some $\al\in[0,1]$
we have a solution $u_\al$ of \reff{bvp1} with (\ref{bvp1}b) replaced by 
\huga{\label{albc}
	v(0)=\al v_{0}+(1-\al)\hat v,
}
 (e.g., for $\al=0$ we have $u\equiv \uh$). We then increase $\al$ by some stepsize 
$\del_\al$ and use $u_\al$ as initial guess for (\ref{bvp1}a), 
(\ref{bvp1}c) and \reff{albc}, ultimately aiming at $\al=1$. 
To ensure that $\|u(1)-\uh\|$ is small, the truncation time $T$ may be set 
free if 
\huga{\label{csss0} \|u(1)-\uh\|_\infty\le \epsi
} 
is violated, and the additional 
boundary condition 
\huga{\label{cssTbc} 
\|u(1)-\uh\|_2^2=\eps^2}
 with fixed $ \eps $ is added, where 
$\ds\|u\|_2=\left(\frac 1 {n_u}\sum_{i=1}^{n_u} u_i^2\right)^{1/2}$ is a weighted (discrete) 
$L^2$ norm, and we initialize 
\huga{\label{epsini}
\eps=\frac 1 {10}\|u(1)-\uh\|_\infty,
} 
which of course is only a rough estimate and highly problem dependent, 
and, like all numerical parameters, can be reset by the user. 
See Remark \ref{derem} for further comments. 

To discretize in time and then solve 
(\ref{bvp1}a), (\ref{bvp1}c) and \reff{albc} (including \reff{cssTbc} if 
$T$ is free) 
we use the BVP solver TOM \cite{mazS2002,MT04,MT2009}%
\footnote{see also \url{www.dm.uniba.it/~mazzia/mazzia/?page_id=433}}, 
in a version {\tt mtom} which accounts for the mass matrix $M$ on the 
lhs of (\ref{bvp1}a)%
\footnote{and which is also used for finding time-periodic orbits in \pdep, see \cite{hotheo,hotutb}}.%
\footnote{We also use a simple BVP solver {\tt bvphdw}, which we mainly 
set up for testing, but which is sometimes
  more robust than the sophisticated methods 
(error estimation and mesh refinement) of {\tt mtom}. In the 
following discussion, we write mtom for the BVP solver, 
but {\tt bvphdw} can similarly be used.} 
This  predictor ($u_\al$) -- corrector ({\tt mtom} for $\al+\del_\al$) 
continuation method corresponds to a ``natural'' parametrization 
of a canonical paths branch by $\al$. 
We also give the option to use a secant predictor 
 \huga{\label{undef0}
u^j(t)=u^{j-1}(t)+\del_\al \tau(t), \quad 
\tau(t)=\bigl(u^{(j-1)}(t)-u^{(j-2)}(t)\bigr)
/\|u^{(j-1)}(\cdot)-u^{(j-2)}(\cdot)\|_2, 
} 
where $u^{j-2}$ and $u^{j-1}$ are the two previous steps. However, 
the corrector still works at fixed $\al$, in contrast to 
the arclength predictor--corrector described next. 

It may happen that no solution of 
(\ref{bvp1}a), (\ref{bvp1}c) and \reff{albc} 
 is found for $\al>\al_0$
for some $\al_0<1$, i.e., that the continuation to the intended
initial states fails. In that case, often the CP branch shows a
fold in $\al$, and we use a modified continuation process, 
letting $\alpha$ be a free parameter and using a pseudo--arclength 
parametrization by $\sig$ in the BCs at $t=0$. 
We set  $ \alpha $ free and add BCs at continuation step $j$,  
 \huga{\label{iccarc}
	\left\langle s,(u(0)-u^{(j-1)}(0)\right\rangle
+s_\alpha(\alpha-\alpha^{(j-1)})=\sig,
}
      with $(u^{(j-1)}(\cdot),\al^{j-1})$ the solution from the
      previous step, and 
$(s,s_\alpha)\in\R^{2Nn}\times \R$ appropriately
      chosen with $\|(s,s_\alpha)\|_*=1$, where $\|\cdot\|_*$ is a 
suitable norm in $\R^{2Nn+1}$, which may contain different 
weights of $v$ and $v_\alpha$.  For $s=0$ and
      $s_\alpha=1$ we find natural continuation with stepsize $\del_\al=\sig$ again. 
To get around folds we may use the secant
$$s:=
\xi\bigl(u^{(j-1)}(0)-u^{(j-2)}(0)\bigr)/\|u^{(j-1)}(0)-u^{(j-2)}(0)\|_2
\text{ and } s_\al=1-\xi 
$$
with small $\xi$, and also a secant predictor 
\huga{\label{undef}
(u^j,\al^j)^{\text{pred}}
=(u^{j-1},\al^{j-1})+\sig \tau} 
for $t\mapsto u^j(t)$ with 
\huga{\label{spred} 
\tau=\xi\bigl(u^{(j-1)}(\cdot)-u^{(j-2)}(\cdot)\bigr)
/\|u^{(j-1)}(\cdot)-u^{(j-2)}(\cdot)\|_2
\text{ and } \tau_\al=1-\xi. 
}
This essentially follows \cite[\S7.2]{grassetal2008}.

Finally, given $\uh$, to calculate $\Psi$, at startup we 
solve the generalized adjoint eigenvalue problem 
\huga{\label{psievp}
\pa_u G(\uh)^T\Phi=\Lambda M \Phi
} for the eigenvalues $\Lambda$ 
and (adjoint) eigenvectors $\Phi$, which also gives the defect $d(\uh)$ 
by counting the negative eigenvalues in $\Lam$. If 
$d(\uh)=0$, then from $\Phi\in \C^{2Nn\times 2Nn}$ 
we generate a real base 
of $E_u(\hat u)$ which we sort into the matrix  $\Psi\in \R^{Nn\times 2Nn}$. 
Algorithm 1 summarizes our method to compute a CP to a CSS. 

\begin{Algorithm}
\fbox{\parbox{0.98\textwidth}{ 
\bci
\item[0. Preparation.] Compute $\Psi$ and the defect $d$ at the CSS. 
If $d\ne 0$, then return (not a saddle point). Otherwise, set $j=0$, 
and, if no initial guess for $T$ is given, compute $T$ from \reff{T0}. 
\eci
 Repeat until $\al=1$ or $j=m$ or until convergence failure. 
 \bci 
 \item[1. BVP solution.]  
 Solve \reff{bvp1} and check \reff{csss0}. If \reff{csss0} is 
violated (or from the start), then 
 augment \reff{bvp1} by \reff{cssTbc}, free $T$, and solve again. 
 \item[2. Next prediction (or stop).] 
 If $\al=1$ or $j=j_{\max}$ then stop and return solution $u$. \\
 If no solution found:
\bci
\item[] 
If arc=0, then stop and return the last solution. 
\item[] If arc=1 and $\del>\del_{\min}$ (else stop and 
return the last solution), then decrease 
$\del$ and go to 1 with new predictor from \reff{undef}. 
\eci 
If solution found: 
\bci
\item[] Let $j=j+1$. 
\item[] If arc=0, then let $\al=\alv_j$, let 
$v_0=\al v_0^*+(1-\al)\vh$, $\ug=u^{(j-1)}$ or set $\ug$ according 
to \reff{undef} (secant predictor), and go to 1.
\item[] If arc=1, then make a new $(\al,u)$ predictor via \reff{undef}, 
set $v_0=\al v_0^*+(1-\al)\vh$, and go to 1. 
\eci 
\eci 
}}
\caption{Continuation algorithm to compute CPs to a CSS $\uh$. Input 
$\uh=(\vh,\lamh)$, initial states $v_0^*$, 
vector $\alv=(\al_1,\al_2,\ldots,\al_m)$ 
of $\al$ values for the 'natural' continuation with $m$ steps, 
or $(\al_1,\al_2)$ and 
number $m$ of arclength steps. Optionally truncation time guess $T$. 
For consecutive calls, Step 0 is omitted, 
and the new predictor is generated from the already computed data. 
Furthermore, let arc=0,1 be the switch controlling whether natural 
or arclength continuation is used. 
\label{CSSalg}}
\end{Algorithm}

\brem\label{derem}{\rm Writing (the discretized version of) \reff{bvp1} 
(and, if switched on, \reff{cssTbc}) as 
\huga{\label{css23} 
H(U)=0,\quad U=(u,T,\al), 
} 
then $u$ is a (numerical) solution of \reff{css23} if $\|H(U)\|_*\le 
{\tt tol}$, where the $M\ddt u-G(u)$ component 
of $H(U)$ is essentially measured in the $\|\cdot\|_\infty$ 
norm (for {\tt mtom} we use the relative error), 
and thus we also use $\|\cdot\|_\infty$ in \reff{csss0}. 
On the other hand, for the active condition \reff{cssTbc} we choose 
the euclidean norm with derivative $\frac 2{n_u}(u(1)-\uh)$ 
(as a row vector) instead of the at first sight more natural 
condition $\|u(1)-\uh\|_2^2=0$, because we thus 
altogether obtain a well conditioned Jacobian 
$\pa_U H(U)$ for the extended system (see \reff{cssjac}).  
For $\eps$ on the order of {\tt tol}, 
\reff{cssTbc} and $\|u(1)-\uh\|_2^2=0$ are 
essentially equivalent, but the additional flexibility via $\epsi$  
(and $\eps$ derived from $\epsi$) 
with for instance $\epsi$ on the order of $10^{-3}\|\uh\|_\infty$ 
turns out to be useful to obtain fast and robust results. Finally, starting with such 
a possibly rather large $\eps$ then also allows to decrease $\eps$ 
a posteriori in a few steps, see the examples in \S\ref{ex-sec}. 
In detail, we solve
\huga{\label{css24}
H(U)=\bpm \Phi(U) \\ \Theta(U) \\
  \CG(u,T) \epm=\bpm 0 \\ 0 \\ 0\epm \in \R^{n_um+k_1+k_2}, 
} 
where $n_u=2Nn_p$ is the number of spatial degerees of freedom ($N=$number 
of states, $n_p=$number of spatial discretization points), 
where $m$ denotes the number of time steps, where the boundary conditions 
are written as 
\huga{
\Phi=\bpm v|_{t=0}-v_0 \\ \Phi_2(u)\epm, \text{ where } 
\Phi_2(u) = \left\{ \barr{cll}
    \Psi(u(1)-\hat{u}) & \in \R^{n_u/2}, & \text{CSS, fixed T,} \\[2mm]
    \bpm \Psi(u(1)-\hat{u}) \\ \|u(1)-\hat{u}\|^2-\varepsilon^2 \epm
    & \in \R^{n_u/2+1}, & \text{CSS, free T,} \\[4mm]
    P(u(1)-\hat{u}) & \in \R^{n_u/2+1}, & \text{CPS, see \reff{bv_cps},}
    \earr \right.
}  
where $\Theta$ contains the arclength boundary
condition \reff{iccarc}, if switched on, and where 
$\CG$ is the discretization of $\left(\ref{bvp1}a\right)$. 
 Thus, $k_1=0$ for CPs to CSSs
with fixed time, and $k_1=1$ for CPs to CSSs with free time or CPs to
CPSs, and $k_2=0$ (natural continuation) or $k_2=1$ (arclength). 
To solve \reff{css24}, given a guess $U_0$ 
from a previous step, we use Newton's Method, i.e., 
\huga{\label{cnew}
  U_{j+1} =U_j - \CA(U_j)^{-1}
  H(U_j), 
} 
\huga{\label{cssjac} \CA=\pa_UH=\bpm 
(\text{{\small $I_{\frac{n_u}2\times \frac{n_u}2}, 0_{\frac{n_u}2\times 
\frac{n_u}2})$}} 
& 0_{\frac{n_u}2\times n_u} &\dots&\dots&\dots&\dots&0_{\frac{n_u}2}&
0_{\frac{n_u}2}\\ 
0 & 0 &\dots&\dots&\dots& D_u\Phi_2&0&0\\ 
  s(1,\ldots,1) &0 &0 &0 &0 & 0 & 0 & s_\alpha \\
  M_1 & H_1 & 0 & \dots & \dots & 0 & (\partial_T \CG)_1 & 0 \\
  0 & M_2 & H_2 & 0 & \dots & 0 & \vdots & 0 \\
  \vdots & 0 & \ddots & \ddots & \dots & 0 & \vdots & 0 \\
  \vdots & \vdots & 0 & \ddots & \ddots & \vdots & \vdots & 0 \\
  0 & \dots & \dots & 0 & M_{m-1} & H_{m-1} & (\partial_T \CG)_{m-1} & 0 
\epm, 
} 
where the first line consist of $n_u/2$ rows, 
where the second line consists of 
$n_u/2+k_1$ rows, and the 
$[s(1,\ldots,1)\ \ 0\ \ \ldots\ 0\ \ s_\al]$---row and the last column of $\CA$ are only present in the 
arclength setting.%
\footnote{We only indicate in the first line the 
dimension of the $0$.} Moreover, 
$M_j=h_j^{-1}M+\frac{1}{2} TG(u_j)$,
$H_j=-h_j^{-1}M+\frac{1}{2} TG(u_{j+1})$,
$(\partial_T \CG)_j=\frac{1}{2}\left(G(u_j) + G(u_{j+1})\right)$, where 
$h_j$ is the $j$-th time step, and $u_j$ the field at time
$t_j$. 
Thus, for $T$ free, $D_u\Phi_2$ contains the row
$2(u(1)-\hat{u})$, which would be ill
conditioned if $u(1)=\hat{u}$, and \reff{cssTbc} is much more robust. 
Of course, $\CA^{-1}$ in \reff{cnew} stands for the 
linear system solver used.
}
\eex \erem

\paragraph{(c) Canonical paths to canonical periodic states.}
For CPs to CPSs the basic idea is a BVP of style (\ref{bvp0}) 
again. However (\ref{bvp0}c) has to be adapted to the CPS case.  The
theoretical truncated BC is 
\huga{
  u(1) \in W_s(\uh) \text{ (and $\|u(1)-\uh_0\|$ small)}, 
} 
where  $W_s(\uh)$ is the stable manifold of the CPS $\uh$, 
and $\uh_0$ is some point on $\uh$. In practice, 
to have a boundary condition analogous to (\ref{bvp1}c), we 
fix an end-point $\uh_0$ on the CPS. 
We then want a boundary condition of the form 
\huga{\label{bv_cps}
  P(u(1)-\uh_0)=0,
} 
with $P\in \R^{(n_u/2+1)\times n_u}$ to yield $n_u/2+1$ boundary conditions. 
For a CPS, the analog to the linearization at a CSS is  
 the monodromy matrix $M_p$,  
which describes the linear effect
of small deviations with respect to one period. It corresponds 
to the time $T_p$ (period of the CPS) map of the variational equation 
\huga{\label{mon123}
\pa_t v=-\pa_uG(\uh(t))v, \quad v(0)=v_0.
} 
The eigenvalues of $M_p$ are called Floquet multipliers, and are
independent of the choice of $\uh_0$, but the eigenvectors depend
on $\uh_0$. Since in \reff{cps21} we start with an autonomous system, 
we always have the (trivial) multiplier  $\ga_1=1$, which corresponds 
to a time shift of $\uh$, and this trivial multiplier can be used 
for assessing the numerics.%
\footnote{For instance, 
we give a warning if for the trivial multiplier we numerically 
have $|\ga_1-1|>{\tt tolfloq}$, with  default 
setting ${\tt tolfloq}=10^{-8}$.} 

The monodromy matrix $M_p$ can
be computed (approximated) as the product of the 
linearizations of (\ref{bvp0}a) at
every $t$--mesh point of $\uh$. However, for OC problems, in 
particular PDE OC problems, we often have both very large (due to 
the anti--diffusion in the co-states) and 
very small (due to diffusion in the states) multipliers%
\footnote{for instance $|\ga|\approx 10^{80}$ for the largest 
multiplier in a small scale PDE problem}, and thus we need a 
particularly stable method for their computation. Here we use a 
periodic Schur decomposition, see \cite{kressner01} and 
\cite{hotheo} and the references therein 
for details. This produces a set of matrices 
\huga{ 
M_p=E \CD E^T 
} with an orthogonal matrix $E$ and an upper triangular matrix
$\mathcal{D}$ with the multipliers on the diagonal. Moreover, 
the adjoint monodromy matrix can be computed without extra effort, 
and sorting of the multipliers in $\mathcal{D}$ is
possible, and as the first $k$ columns of $E$ are a basis of the span of
the first $k$ eigenvectors, we can compute projections on
eigenspaces this way. 

The projection $P$ in \reff{bv_cps} needs to provide 
$n_u/2+1$ boundary conditions  
by projecting onto the center--unstable eigenspace, 
associated with multipliers $\gamma$ with $\abs{\gamma}\ge 1$, 
where the translational eigenspace associated with the 
trivial multiplier $\gamma=1$ is included because 
 we want to fix the truncation point $\uh_0$ on the CPS. 
As for the CSS case, 
the dimension of the center--stable eigenspace 
is at most $Nn$ and thus this is the only case in which a
canonical path can be computed. This is called saddle point
property (SPP) of the CPS, 
see \cite{grassetal2008}. 
Given a CPS with the SPP we thus compute the matrix $P$ as the
projection on the center-unstable eigenspace, i.e.~on the eigenspace
spanned by Floquet-multipliers $\gamma$ with $\abs{\gamma}\geq1$. 

Thus, altogether we have $n_u/2$ BCs 
(\ref{bvp0}b) for the initial states, $n_u/2+1$ BCs \reff{bv_cps}, 
and $n_u$ ODEs (\ref{bvp0}a) for the $n_u$ unknowns $(u_i)$, 
and the free truncation time $T$ is the $(n_u+1)^{\text{ths}}$ 
unknown, i.e., 
\begin{subequations}\label{cps0}
\hual{
\pa_t u&=-TG(u),\\
v|_{t=0}&=v_0,\\
u(1)&\in E_s(\uh_0). 
} 
\end{subequations}
Moreover, as in \reff{csss0} we additionally require, 
for a given $\epsi>0$, 
\huga{\label{cps1}
\|u(1)-\uh_0\|<\epsi.
}

The continuation method in the initial states 
is the same as for CPs to CSS, i.e.~we have natural parametrization 
with and without secant predictor, and arclength parametrization. 
This includes the adaptation of the truncation time $T$ to ensure 
(\ref{cps0}c). 
Similar to \reff{T0}, 
an estimate for $T$ can be obtained from the largest (in modulus) stable 
multiplier $\ga$, which we denote by $\ga_2$, 
with the trivial multiplier denoted by $\ga_1=1$. 
In the linear regime (small deviation from $\uh$) 
we then have 
\huga{\label{Tp0}
\|u(1)-\uh_0\|\sim \er^{|\ga_2|T_p}\|u(0)-\uh_0\|.
}
However, for small $\al$ we may use $T=lT_p$ as default initialization, 
with rather small $l$ ($l=2,3$)  
and then for larger $\al$ add periods of the CPS at the end of the 
canonical path within the
continuation process if necessary, i.e., 
if the deviation from the CPS becomes to large. In detail, 
to ensure \reff{cps1} (which is {\em never} used to extend 
\reff{cps0}) we use an add hoc step additional to the 
discretization in time of \reff{cps0} and the solution of the 
obtained algebraic system by Newton's method: 
\huga{\label{addT}
\begin{split}&\text{After solving \reff{cps0} we check \reff{cps1} 
If \reff{cps1} is violated, then we add multiples of the}\\
&\text{period $T_p$ of the CPS to the truncation time $T$, 
extend the last computed CP}\\&\text{by extra copies of the CPS, and run the Newton 
loop again.}
\end{split}
} This appears to be a new idea, which 
allows to start with rather small initial $T$, 
and works very well in all our applications, see \S\ref{cps-sec} for 
details and examples.

\begin{Algorithm}
\fbox{\parbox{0.98\textwidth}{ 
\bci
\item[0. Preparation.] Choose a point $\uend=(\vh_0,\lamh_0)$ on the CPS and 
compute the projection 
$ P $ onto the center-unstable eigenspace of the monodromy matrix 
at $\uend$ via \reff{bv_cps}, i.e.~the eigenspace associated to 
Floquet-multiplier $\gamma$ with $\abs{\gamma}\geq1$, including 
the defect $d$. \\
If $d\ne 0$, then return (not a saddle point CPS).\\
Let $v_0=\al v_0^*+(1-\al)\vh_0$, and compute 
a guess for a canonical path to $\uend $, initially 
($l$ copies of, if $T=lT_p$) the CPS itself. Set $j=0$. 
\eci
 Repeat until $\al=1$ or $j=m$ or until convergence failure. 
\bci 
 \item[1. BVP solution.]  
 Solve \reff{cps0} for $u$. 
\item[2. Target check.] Check \reff{cps1}. If \reff{cps1} is violated, 
then proceed as in \reff{addT}, i.e., extend $T$ and go to 1. 
 \item[3. Next prediction (or stop).] 
 If $\al=1$ or $j=j_{\max}$ then stop. \\
 If no solution found:
\bci
\item[] 
If arc=0, then stop and return the last solution. 
\item[] If arc=1 and $\del>\del_{\min}$ (else stop and 
return the last solution), then decrease 
$\del$ and go to 1.~with new predictor from \reff{undef}. 
\eci 
If solution found:
\bci
\item[] Let $j=j+1$. 
\item[] If arc=0, then let $\al=\alv_j$, let 
$v_0=\al v_0^*+(1-\al)\vh_0$, $\ug=u^{(j-1)}$ or set $\ug$ according 
to \reff{undef} (secant predictor), and go to 1.
\item[] 
If arc=1, then make a new $(\al,u)$ predictor via\reff{undef}, 
set $v_0=\al v_0^*+(1-\al)\vh$, and go to 1. 
\eci 
\eci 
}}
\caption{Continuation algorithm to compute CPs to a CPS $\uh$. Input 
argument CPS $\uh=(\vh,\lamh)$, otherwise as for Algorithm \ref{CSSalg}.  
Again, for consecutive calls, Step 0 is omitted, 
and the new predictor is generated from the already computed data. 
arc=0,1 again determines whether natural 
or arclength continuation is used. 
\label{CPSalg}}
\end{Algorithm}

\brem\label{corem}{\rm There are further (and more sophisticated) 
methods for computing connecting orbits, including (homo-- and) 
heteroclinic orbits which also converge to some prescribed 
solutions as $t\to-\infty$. See, e.g., \cite{beyn90} for a detailed 
analysis of the 'standard' projection boundary condition, 
\cite{Pam01} and \cite{BPS01} for the so-called 
boundary corrector method, and \cite{DKK18a,DKK18b} for 
the special case of connecting orbits in $\R^3$. In particular, 
for connecting orbits to cycles (periodic orbits) 
these methods use a free truncation time $T$ together 
with certain phase conditions to ensure that $\|u(1)-\uh\|$ is 
small, where $\uh$ may vary on the cycle. 

Here, we fix $\uh$ and thus the (asymptotic) phase, 
and use \reff{addT} 
to ensure $\|u(1)-\uh\|<\eps$.  
From the application point of view, 
it is important to keep the defining systems for CPs as small 
as possible, and in particular to put the computation of the 
CPS and the projections at some target point $\uh$ on the CPS 
into a preparatory step. 
Algorithm 2 summarizes our method to compute a CP to a CPS.  
}\eex\erem 

\section{Examples and implementation details for 
CPs to CSSs}\label{ex-sec}
To explain how to use \pdep\ to calculate CSS and canonical paths 
we consider three example problems of type \reff{oc1}: The {\tt sloc} 
(shallow lake OC) problem from \cite{GU17}, the  
{\tt vegoc} (vegetation OC) problem from \cite{U16}, 
and the boundary fishing problem {\tt lvoc} (Lotka-Volterra OC) from \cite{GUU19}.  For all 
examples we first briefly sketch the models, and then summarize the 
contents of the respective demo folder and explain the most important files 
in some detail. For the first model we also explain the general 
setup how to initialize the spatial domain and discretization, the 
rhs, the computation of CSS, and the 
OC related routines of \pdep. For all models we give some plots, 
but for details and interpretation of the results we refer 
to the respective papers. 

\subsection{Optimal distributed control of the phosphorus in a shallow lake}\label{slsec}
\def\dname{sloc}\def\dhome{ocdemos}
Following \cite{BX08}, in \cite{GU17} we consider a model 
for phosphorus $v=v(x,t)$ in a shallow lake, and phosphate load $\k=\k(x,t)$ 
as a control. In 0D, i.e., in the ODE setting, this has been analyzed 
in detail for instance in \cite{KW10}. 
Here we explain how we set up the spatial problem in \pdep. 
The model reads 
\begin{subequations} 
\label{sldiff1}
  \begin{align}
    &V(v_0(\cdot))\stackrel{!}{=}\max_{\k(\cdot,\cdot)}J(v_0(\cdot),\k(\cdot,\cdot)), \qquad 
    J(v_0(\cdot),\k(\cdot,\cdot)):=\int_0^\infty\er^{-\rho t}
J_{ca}(v(t),\k(t)), 
\dd t\label{sldiffusion_model_obj2}\\
\intertext{where $J_c(v,\k)=\ln \k-\ga v^2$, 
$J_{ca}(v(t),\k(t))=\frac 1 {|\Om|}
\int_\Om J_c(v(x,t),\k(x,t))\!\Dx$ as in (\ref{oc1}b),  and 
$v$ fulfills the PDE 
}
 &\pa_t v(x,t)=D\Delta v(x,t)+\k(x,t)-bv(x,t)+\frac{v(x,t)^2}{1+v(x,t)^2},\label{sldiffusion_model_dyn}\\
    &\pa_\nu v(x,t)_{\pa\Omega}=0,\label{sldiffusion_model_bc}\qquad 
    v(x,t)_{t=0}=v_0(x),\quad
    x\in\Om\subset\R^d. 
  \end{align}
\end{subequations}
The parameter $b>0$ is the phosphorus degradation rate, and $\ga>0$ 
are ecological costs of the phosphorus contamination $v$. 
One wants a low $v$ for ecological reasons, but for economical reasons a high phosphate load $\k$, for instance from fertilizers used by farmers. 
Thus, the objective function consists of the concave
increasing function $\ln \k$, and the concave decreasing function
$-\ga v^2$.  In the demo directory {\tt sloc} we consider the parameters 
\hual{ D=0.5,\ \rho=0.03,\ \ga=0.5,\ b\in (0.5,0.8) \text{ (primary
    bif.~param.)}.\label{parsel}
}
With the co-state $\lam$, the canonical system for \reff{sldiff1} becomes, 
with $\ds \k(x,t)=-\frac{1}{\lam(x,t)}$, 
\begin{subequations}
  \label{slcan2}
  \begin{align}
    {\pa_t} v(x,t)&=\k(x,t)-bv(x,t)+\frac{v(x,t)^2}{1+v(x,t)^2}
+D\Delta v(x,t),\label{sldiffusion_model_cansys1}\\
    \pa_t \lam(x,t)&=2\ga v(x,t)+\lam(x,t)\left(\rho+b-\frac{2v(x,t)}{\left(1+v(x,t)^2\right)^2}\right)-D\Delta \lam(x,t),\label{sldiffusion_model_cansys2}\\
    \pa_\nu v&=\pa_\nu \lam=0\text{ on }\pa\Om,\label{sldiffusion_model_cansys3}\\
    \evalat{v(x,t)}_{t=0}&=v_0(x),\quad x\in\Omega. 
  \end{align}
\end{subequations}

\subsubsection{Canonical steady states}\label{csssec}
To compute CSS we use a standard \pdep\ setup for (the steady version) 
of \reff{slcan2}. 
As an overview, Table \ref{sloctab1} lists the scripts and functions in {\tt ocdemos/sloc}; these will be explained in more detail below, 
but for new users of \pdep\ we refer to \cite{actut} 
for an introduction into the basic \pdep\ data structures and setup 
of elliptic systems. 

\taskip\begin{table}
{\small\begin{tabular}{p{0.14\tew}|p{0.8\tew}}
bdcmds1D&script to compute bifurcation diagrams of CSS\\%
cpdemo1D&script to compute CPs\\%
skibademo&script to compute Skiba paths, see \S\ref{skibasec}\\
\hline
slinit&init routine; set the \pdep\ parameters 
to standard values, then {\em some} parameters to problem specific values; 
initialize the domain, set an initial guess $u$, 
and find a first steady state by a Newton loop \\%
oosetfemops&set FEM matrices (stiffness K and mass M) \\
slsG; slsGjac&$G(u)$ resp. the Jacobian $\pa_u G(u)$ for \reff{slcan2}\\
slcon; sljcf&extract control from states/costates; 
compute the current value $J_c$\\
cssvalf; psol3D&print CSS value and characteristics; 
mod of psol3D to plot several solutions in one fig. 
\end{tabular}}
\caption{{\small Scripts and functions in {\tt ocdemos/sloc} (for the 1D case;  
 some additional functions  for the 2D case are also in the folder). 
}\label{sloctab1}}
\end{table}\teskip

The \pdep\ FEM setup converts the PDE \reff{slcan2} into the ODE system 
(or algebraic system for steady states)  
\huga{\label{sla}
M\ddt u=-G(u),\text{ where } G(u)=-\CK u-Mf(u). 
}
In \reff{sla}, 
$M$ is the mass matrix of the FEM, $\CK=\bpm K&0\\0&-K\epm$ is the 
stiffness matrix, where $K$ is the 1--component stiffness matrix 
corresponding to the scalar (Neumann--)Laplacian, or, more precisely, 
$M^{-1}K\approx -\Delta$, and we put 'everything but diffusion' 
into the 'nonlinearity' $f$. 
As usual, our basic structure is a \mlab\ struct 
{\tt p} as in {\tt p}roblem, which has a number of fields 
(and subfields), e.g., {\tt p.fuha, p.pdeo, p.u, p.hopf, 
p.file, p.plot, p.sw, p.nc} 
 which contain for instance function handles to the right hand side 
({\tt p.fuha}), the FEM mesh ({\tt p.pdeo}), the current solution {\tt p.u}, 
data for Hopf orbits (CPS, {\tt p.hopf}), filenames/counters ({\tt p.file}), 
plotting controls ({\tt p.plot}), and switches ({\tt p.sw}) 
and numerical constants ({\tt p.nc}) used in the 
numerical solution, such as {\tt p.nc.tol}, where $u$ is taken 
as a solution of $G(u)=0$ if $\|G(u)\|_\infty<{\tt p.nc.tol}$. 
However, most of these can be set to standard values 
by calling {\tt p=stanparam(p)}. 
In standard problems the user only has 
to provide: 
\bcen
\item The geometry of the domain $\Om$, and in the \oop\ setting used here 
 a function {\tt oosetfemops} used to generate 
the needed FEM matrices. 
\item Function handles {\tt sG} implementing $G$, and, for speedup, {\tt sGjac}, implementing 
the Jacobian. 
\item An initial guess for a solution $u$ of $G(u)=0$, i.e., an 
initial guess for a CSS. 
\ecen 

Typically, the steps 1-3 are put into an init routine, 
here {\tt p=slinit(p,lx,ly,nx,sw,ndim)}, where {\tt lx,ly,nx} and {\tt ndim} 
are parameters to describe the domain size and discretization%
\footnote{ly (irrelevant) and ndim(=1) play no role in this tutorial, 
but we kept them in {\tt slinit} because the same init routine is also used 
in 2D}, 
and {\tt sw} is used to set up different initial guesses, 
see Listing \ref{l1}. For CSS computations the only 
additions/modifications to the standard \pdep\ setting 
are as follows: (the additional function handle) 
{\tt p.fuha.jcf} should be set to the local current 
value objective function, here {\tt p.fuha.jc=@sljcf} (see Listing \ref{l7}), 
and {\tt p.fuha.outfu} to {\tt ocbra}, 
i.e., {\tt p.fuha.outfu=@ocbra}. This automatically puts $J_{ca}(u)$ 
at position 4 of the calculated output--branch. Finally, it is useful 
(for instance for plotting) to set 
{\tt p.fuha.con=@slcon}, 
where {\tt \k=slcon(p,u)} (see Listing \ref{l6}) extracts the control $\k$ from the states $v$, 
costates $\lam$ and parameters $\eta$, all contained in the vector
{\tt u}.\footnote{We do not use {\tt slcon} in {\tt slsG}. However, putting this function into {\tt p} has the advantage that for instance plotting and extracting the value of the control can easily be done by calling some convenience functions of \poc.}

\hulst{
caption={{\small First 4 lines of {\tt \dname/slinit.m}, which collects some typical initialization commands. {\tt p=stanparam(p)} in line 2 sets the \pdep\ parameters, switches and numerical constants to standard 
values; these can always be overwritten afterwards, and some typically are. 
For instance, in line 3, besides setting the function handles to the rhs 
(necessarily problem dependent), here we overwrite the standard branch-output 
 p.fuha.outfu=@stanbra with the OC standard output {\tt ocbra}. 
The remainder of slinit follows the general rules of initialization 
in the \oop\ setting, explained in \cite{actut}. We only comment that 
while we here no further discuss the 2D case, 
the same init-file is used for the 1D and 2D cases, controlled via the 
input argument {\tt ndim}, and that via the switch {\tt sw} 
the user can choose an initial guess $u$ near one of the two spatially homogeneous branches: {\tt sw=1} leads to the so called 'flat state clean' (FSC) 
which turns into 'flat state intermediate' (FSI), and 
{\tt sw=2} leads to 'flat state muddy' (FSM); see \cite{GU17} for 
this nomenclature.}},
label=l1,language=matlab,stepnumber=5,linerange=1-4}{\dhome/sloc/slinit.m}

\hulst{label=l6,language=matlab,stepnumber=5, firstnumber=1}
{\dhome/sloc/oosetfemops.m}
\hulst{
label=l6,language=matlab,stepnumber=5, firstnumber=1}
{\dhome/sloc/slcon.m}

\hulst{caption={{\small {\tt oosetfemops.m, slcon.m} and {\tt sljcf.m} 
from {\tt \dname}. 
{\tt oosetfemops} assembles and stores the needed FEM matrices, 
{slcon} computes the control (here very simple), and {\tt sljcf} 
computes the 
current value $J$ from the states/costates. }},
label=l7,language=matlab,stepnumber=5, firstnumber=1}
{\dhome/sloc/sljcf.m}

\hulst{caption={{\small {\tt \dname/slsG.m}. The rhs of (\ref{slcan2}a,b): first we 
extract the parameters (line 3) and the fields $v,\lam$ (line 4) (with l for $\lam$)  from the full solution vector u which contains the pde variables $(v,\lam)$ and the parameters. Then we compute the nonlinearity $f$ (here and usually containing everything but diffusion), and then we compute the residual 
(line 9) using the pre-assembled stiffness and mass matrices in {\tt p.mat.K} 
and {\tt p.mat.M}, see {\tt oosetfemops.m}. The Jacobian in {\tt \dname/slsGjac.m} works accordingly.}},
label=l4,language=matlab,stepnumber=5, firstnumber=1}
{\dhome/sloc/slsG.m}

\begin{figure}[ht]
\begin{tabular}{lll}
(a) BD of CSS&(b) BD, current values $J_{ca}$& (c) example CSS\\
\ig[width=0.3\tew]{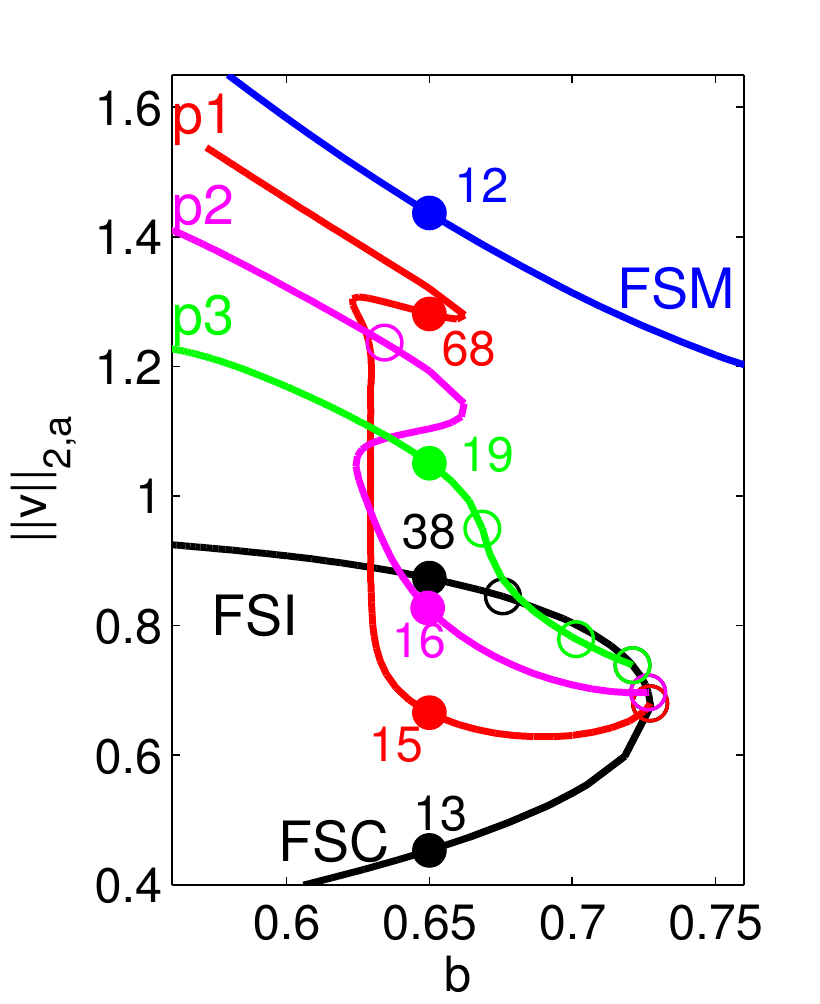}&
\ig[width=0.3\tew]{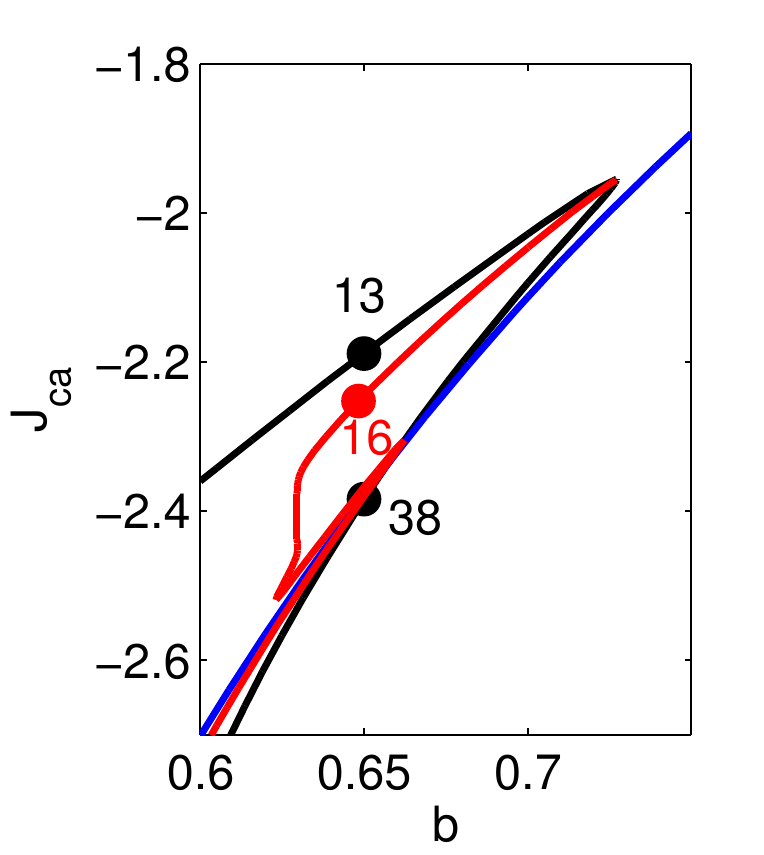}&
\raisebox{30mm}{\begin{tabular}{l}
\ig[width=0.15\tew]{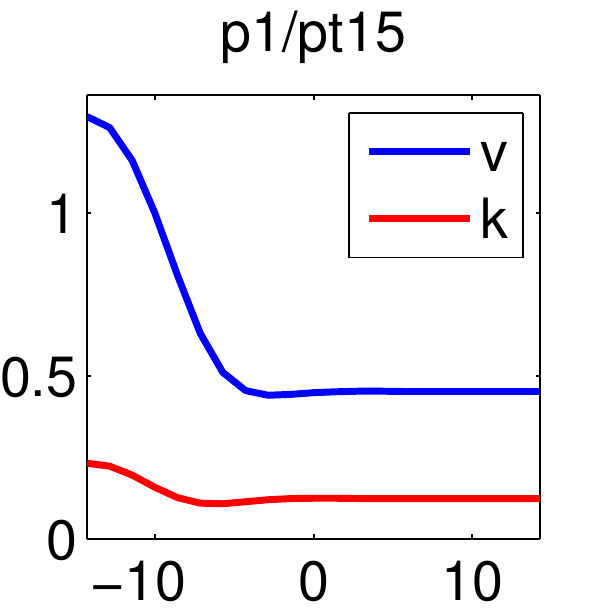}\ig[width=0.15\tew]{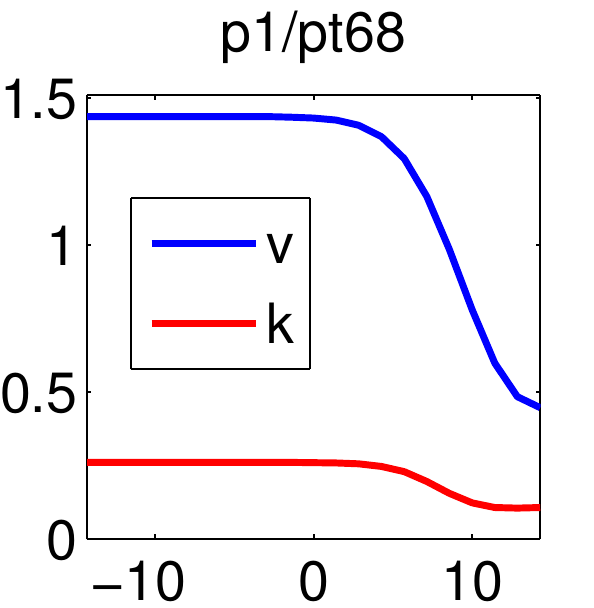}\\
\ig[width=0.15\tew]{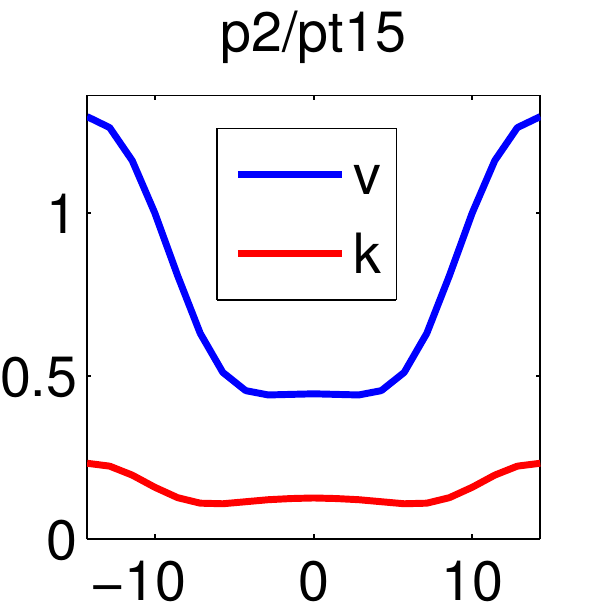}\ig[width=0.15\tew]{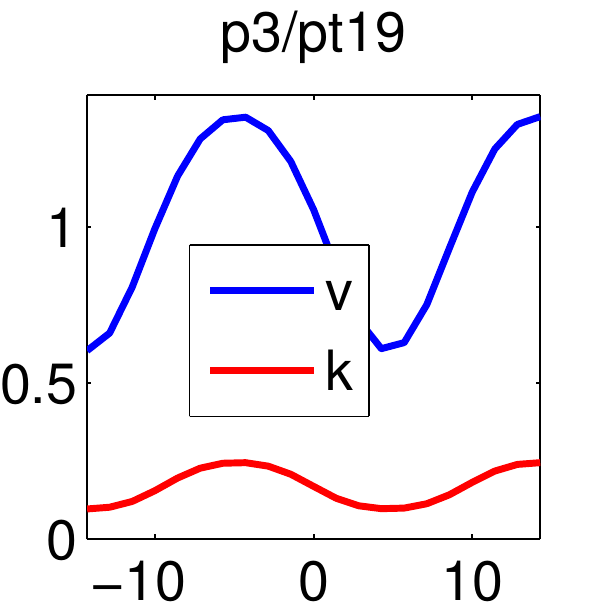}
\end{tabular}}
\end{tabular}
\caption{{\small Example bifurcation diagrams and solution plots from running 
{\tt bdcmds.m}. For $b<b_{\text{fold}}\approx 0.73$ there are 
three branches of 
FCSS, here called FSC (Flat State Clean, low $v$), FSI (Flat State Intermediate), 
and FSM (Flat State Muddy, high $v$). On FSC$\cup$FSI there are a number 
of bifurcations to patterned CSS branches. 
See \cite{GU17} for further discussion.} }
  \label{slfig1}
\end{figure}

At the end of {\tt slinit}, we call a Newton--loop to converge 
to a (numerical) CSS, which is called 'flat', i.e., spatially 
homogeneous. 
By calling {\tt p=cont(p)}  we can continue such a state in some parameter. If {\tt p.sw.bifcheck>0}, 
then \pdep\ detects, localizes and saves to disk 
bifurcation points on the branch. Afterwards, the bifurcating branches 
can be computed by calling {\tt swibra} and {\tt cont} again. 
These (and other) \pdep\ commands (continuation, branch switching, and plotting) 
are typically put into a script file, here {\tt bdcmds1D.m}, 
see Listing \ref{l2}, which we recommend to organize into cells. 
There are some modifications to the standard \pdep\ plotting commands, 
see, e.g., {\tt plot1D.m}, to plot $v$ and $\k$ simultaneously. These work as usual by overloading the respective 
\pdep\ functions by putting the adapted file in the current directory. 
See Fig.~\ref{slfig1} 
for example results of running {\tt bdcmds1D}. 

\hulst{caption=
{{\small {\tt \dname/bdcmds.m} (first 8 lines), 
following standard \pdep\ principles. 
The remainder of {\tt bdcmds.m} deals with plotting, see the source code.}},
label=l2,language=matlab,stepnumber=5, linerange=3-10}
{\dhome/sloc/bdcmds1D.m}

\subsubsection{Canonical paths}\label{cpsloc}
For OC problems, the computation of CSS is just the first step.  
The next goal is to
calculate CPs from some starting state $v(0)$ to a CSS
$\uh_1$ with the SPP. For this we use the continuation algorithm {\tt
  isc} which is essentially a wrapper for the BVP solvers {\tt mtom} and 
{\tt bvphdw}, and which for CSS as targets implements Algorithm \ref{CSSalg}. 
The data is again stored in a problem structure {\tt p} which has
a number of general options/parameters in {\tt p.oc}, options 
specific to the behavior of the BVP solvers stored in {\tt
  p.tomopt}, and solution data stored in {\tt p.cp}. 
In particular, the {\tt oclib} routines re--use the data and functions 
(FEM data, function handles) already set up for the computation 
of the CSS (or CPS), and no new functions 
need to be set up. 
The convenience function {\tt ocinit} sets most OC parameters to
standard values and, if provided with the corresponding data, 
the starting states and end point of the canonical path, 
i.e., the target CSS $\uh$, or the target $\uh_0$ on a CPS $\uh$. 
This function is the analog of {\tt stanparam} in the
CSS setting. For a first call the user only has to set the
parameters at the top of Table \ref{tab_user_ocinit}, the estimated
truncation time and the (initial) number of mesh points. 
However, as usual, the user can, and sometimes has to, change
some of the standard options.

\input{slfig2}

After setting up the data structures (via {\tt ocinit} or modifications 
and possibly further commands) in the struct {\tt p}, 
the computation of CPs is started by a call of {\tt p=isc(p,alvin,varargin)} 
with default input parameters {\tt p} and {\tt alvin}, 
where {\tt alvin} is a vector
of continuation steps. For arclength continuation, the third input
may contain the number of arclength continuation steps. 
For consecutive arclength calls, one can also set {\tt alvin=[]} to
directly start with arclength continuation. However a secant and two
values for $\alpha$ have to be given via {\tt p.oc.usec} and {\tt
  p.hist.alpha} (if natural continuation was done before, these fields
are filled). See Section \ref{oc-lib-sec} for examples and further 
description of {\tt isc}, other OC related
\pdep\ functions, and the parameters in {\tt p.oc} and {\tt
  p.tomopt}. Here we continue with a brief description of 
CP results for the SLOC problem.  
The canonical path related computations are done in the command files
{\tt cpdemo*}. We first compute paths from a patterned CSS to
a flat CSS and vice versa in {\tt cpdemo1D},  see Listing \ref{cpdemo1D} and
Figure \ref{slfig2}, while {\tt cpdemo2D} does
the same in $2D$. The file {\tt skibademo}, to be run after {\tt cpdemo1D}, 
computes 
some Skiba paths, see \S\ref{skibasec}. 

\input{octab}

\taskip 
\begin{table}[ht]  
{\small\begin{tabular}{p{0.1\tew}|p{0.85\tew}} 
 name& description \\  \hline
 cp.u& solution (CP) generated by {\tt isc}; used as initial guess 
for next continuation step, if already set by previous call to 
{\tt isc} or if set externally. \\
 cp.t &time mesh generated by {\tt isc} (or set externally)\\ 
cp.par&solution parameters, i.e., truncation time $T$ in par(1), 
current $\alpha$ in par(2) (for arclength)\\
hist.alpha& vector of the $\al$ values in the continuation \\ 
hist.vv & vector of the objective values of the canonical path for the $\al$ stored in {\tt hist.alpha} \\ 
hist.u &CPs at continuation steps stored in hist.alpha\\
hist.t &time-meshes (normalized to $[0,1]$) of continuation steps\\
hist.par & parameter values ($T$ and $\al$) of the continuation steps\\ 
\end{tabular}} 
\caption{{\small Additional fields in {\tt p}, typically 
set/maintained/updated by {\tt isc} 
({\tt p.hist} only maintained/updated if {\tt oc.retsw=1}). 
\label{tab_isc}}}
\end{table} 
\teskip

\taskip\hulst{caption= {{\small Cells 1 and 2 
from {\tt \dname/cpdemo1D.m} (to be  
run after {\tt  bdcmds1D.m}). 
We compute a CP  from the states of {\tt p3/pt19} to the CSS {\tt f1/pt13} by ``natural continuation'' in the initial states. In the remainder of {\tt cpdemo1D} we illustrate arclength continuation (also in preparation of 
{\tt skibademo}), and 
the extraction of data from the continuation history.}},
              label=cpdemo1D,language=matlab,stepnumber=5,
                 linerange=1-18}{\dhome/sloc/cpdemo1D.m} \teskip

\subsubsection{Main {\tt oclib} functions}\label{oc-lib-sec}
The functions {\tt ocinit} and {\tt isc} are the main user interface
functions for CP numerics. 
Essentially, after having set up {\tt p} as in \S\ref{csssec} for the
CSS, including {\tt p.fuha.jcf}, the user does not need to
set up any additional functions to calculate canonical paths and
their values. We give the   
signatures and some general remarks on the arguments and behavior 
of  {\tt ocinit} and {\tt isc}, 
with the Cells referring to Listing \ref{cpdemo1D}.

\bci
\item[{\tt p=ocinit(p,varargin)} (Cell 1).] Convenience function
  (similar to {\tt p=stanparam(p)}) to generate a standard problem
  structure {\tt p} where most parameters are set to
  standard values. These are parameters for {\tt mtom} (see {\tt
    tom/tomset.m}), and parameters for {\tt isc}, see Table
  \ref{tab_user_ocinit} for an overview. If {\tt
    varargin=\{sd0,sp0,sd1,sp1\}}, then {\tt ocinit} also sets the states of
  the solution stored in {\tt sd0/sp0} as the starting states and the
  solution stored in {\tt sd1/sp1} as the aimed CSS/CPS.  Typically
  some of the options should be overwritten in the further setup.
\item[{\tt p=isc(p,alvin,varargin)} (Cell 1).] {\tt p} is the problem
  structure containing the options/parameters described above 
in {\tt p.oc} and {\tt p.tomopt},
  see Tables \ref{tab_user_ocinit} and \ref{tab_isc}, and the solution
  in {\tt p.cp}. {\tt alvin} is the
  vector of desired $\al$ values for the continuation, for
  instance {\tt alvin=[0.25 0.5]}. If {\tt varargin=nsteps} is given
  as a third input, arclength continuation is started after the last
  value of {\tt alvin} with {\tt nsteps} steps, or until $\al=1$ is reached. 
If previous calls of
  {\tt isc} are present, then we can directly 
start arclength continuation  by calling
  {\tt isc} with empty {\tt alvin}. After the first call of {\tt isc}
  some additional fields are set in {\tt p.oc}, containing, 
e.g., the current secant, the last
  starting point in the continuation, and, if desired via 
{\tt p.oc.retsw=1}, the
  solutions at different continuation steps, see Table \ref{tab_isc}.
  \eci

  \brem\label{tomo}{\rm Concerning the original TOM options we
    typically run {\tt isc} with weak error bounds and what appears to
    be the fastest monitor and order options, i.e., {\tt
      tomopt.Monitor=3; tomopt.order=2}. Once continuation is successful (or
    also if it fails at some $\al$), we can always postprocess by
    calling {\tt mtom} again with a higher order, stronger error
    requirements, and different monitor options. See the original TOM
    documentation. The most convenient way to do so is to call {\tt
      isc} with {\tt alvin=1} again after resetting TOM-related
    options.  }\eex\erem

   There are a number of
  additional functions for internal use, and some convenience
  functions, which we briefly review as follows:

  \bci
\item[{\tt [Psi,mu,d,t]=getPsi(s1)}.] For CPs to CSSs only:
  compute $\Psi$, the eigenvalues {\tt mu}, the defect d, and a
  suggestion for $T$. This becomes expensive with large
  $2nN$ (number of spatial DoF).
\item[{\tt [Fu1,Fu2,d]=floqpsmatadj(p.opt.s1).}] For CPs to 
CPSs only: computes 
(by periodic Schur decomposition) the
  projection on the center-unstable eigenspace in {\tt Fu2}
  and the Floquet-multiplicators in {\tt d}. 
  Expensive, but has to be done only once
  for each CPS. 
\item[{\tt [sol,info]=mtom(ODE,BC,solinit,opt,varargin)}.]
  Modification of TOM, which allows for $M$ in (\ref{bvp1}a).  Extra
  arguments $M$ and {\tt lu,vsw} in {\tt opt}.  If {\tt opt.lu=0}, then
  $\backslash$ is used for solving linear systems instead of an
  LU--decomposition, which becomes too slow when $2nN\times m$ becomes
  too large. See the TOM documentation for all other arguments
  included in {\tt opt}, and note that the modifications in mtom can
  be identified by searching ``HU'' in {\tt mtom.m}. Of course {\tt
    mtom} (as any other function) can also be called directly, which
  for instance can be useful to postprocess the output of some
  continuation by changing parameters by hand.
\item[{\tt sol=bvphdw(ODE,BC,solinit,tomopt,opt).}] A simple Newton
  solver for CPs, which was mainly used for testing but is sometimes
  more robust than the sophisticated methods 
(error estimation and mesh refinement) of {\tt mtom}. 
\item[{\tt f=mrhs(t,u,\k,opt); J=fjac(t,u,opt);} and {\tt
    f=mrhse(t,u,\k,opt); J=fjace(t,u,opt)}.] The rhs \linebreak and
  its Jacobian to be called within {\tt mtom} resp. {\tt bvphdw} (see
  the respective wrapper files). These are wrappers which
  calculate $f$ and $J$ by calling the resp.~functions in the
  \pdep--struct {\tt p.opt.s1}, which were already set up and used to
  calculate the CSS/CPS.  Similar remarks apply to {\tt mrhse} and
  {\tt fjace} for the arclength continuation.
\item[{\tt bc=cbcf(ya,yb,opt);[ja,jb]=cbcjac(ya,yb,opt);
    bc=cbcfe(ya,yb);[ja,jb]=cbcjace(ya,yb).}] The boundary conditions
  (in time) for \reff{bvp1} resp. \reff{bv_cps} and the associated
  Jacobians again with wrappers to be used for {\tt mtom} and {\tt
    bvphdw} at once. The {\tt *e} (as in {\tt e}xtended) versions are
  for arclength continuation again.  
\item[{\tt [jval,jcav,jcavd]=jcaiT(s1,cp,rho)} and {\tt
    djca=disjcaT(s1,cp,rho).}] Computes the  value
  \huga{\label{jadef} {\tt jval}=J(u)=\int_0^T \er^{-\rho
      t}J_{ca}(v(t,\cdot),\k(t,\cdot))\dd t
} of the solution $u$ in
  {\tt cp} (with $J_c$ taken from {\tt s1.fuha.jcf}), and also returns 
$J_{ca}$ and $\er^{-\rho t}J_{ca}$ along the CP for easy plotting, 
cf.~Fig.~\ref{slfig2}(b). 
\item[{\tt [di,d2]=tadev(p).}] Deviations (in sup norm and 
euclidean norm) of endpoint of CP from target, see \reff{csss0} and 
\reff{cssTbc}. 
\eci 

  \brem{\rm There are also some plotting function in {\tt oclib}, which however should be seen as templates for plotting of canonical paths, including diagnostic plots to check the convergence behavior of
    the canonical path as $t\ra T$, cf.~(d),(f) in
    Fig.~\ref{slfig2}. 
The function {\tt slsolplot(p,view)} in the {\tt sloc} demo
    directory can serve as a template how to set up such plots. }\eex\erem

\teskip

\subsubsection{Skiba points}\label{skibasec}
In ODE OC applications, if there are several locally stable CSS, 
then often an important issue is to identify their domains of attractions. 
These are separated by so called threshold or Skiba--points (if $N=1$) 
or Skiba--manifolds (if $N>1$), see \cite{Skiba78} and \cite[Chapter 5]{grassetal2008}. 
Roughly speaking, these are initial states from which there are several 
optimal paths with the same value but leading to different CSS. Here we 
give an example for the SLOC model 
how to compute a patterned Skiba point between FSC and FSM. 

In Cell 3 of {\tt cpdemo1D.m} we attempt to find a path from $v_{\text{PS}}$ 
given by 
{\tt p1/pt68} to $(v,q)_{\text{FSC}}$ given by {\tt FSC/pt13}; this fails 
due to the fold in 
$\al$. However, for given $\al$ we can also try to find a path from 
the initial state 
$v_\al(0):=\al v_{\text{PS}}+(1-\al)v_{\text{FSC}}$ to the FSM, and compare to 
the path to the FSC. For this, we can use the problem structure {\tt p} computed in Cell 3. The initial states for the k'th $\al$ value of {\tt p.hist.alpha} are, due to {\tt p.opt.retsw=1}, stored in {\tt p.hist.u\{k\}(:,1)} i.e. as the starting point of the canonical path associated to the given $\al$ value.

\begin{figure}[ht]
\bce{\small	\begin{tabular}{ll}
(a) A Skiba point at $\al=0.454$&(b) Paths to FSC (blue) and FSM (red)\\
\ig[width=0.23\tew]{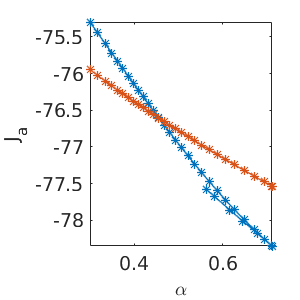}&\ig[width=0.22\tew]{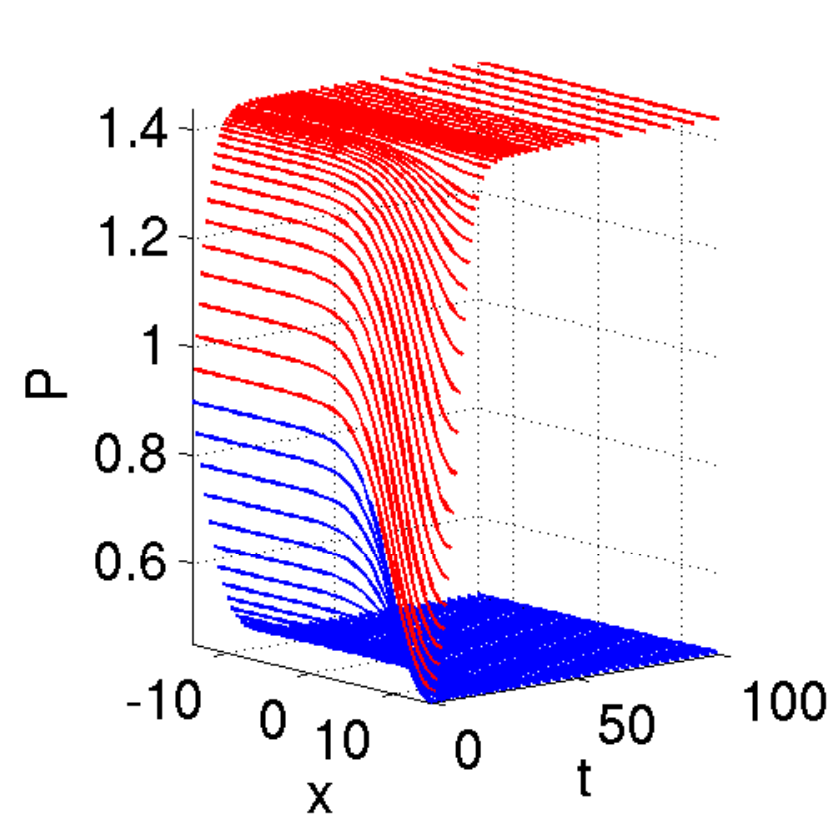}
\raisebox{0mm}{\ig[width=0.2\tew]{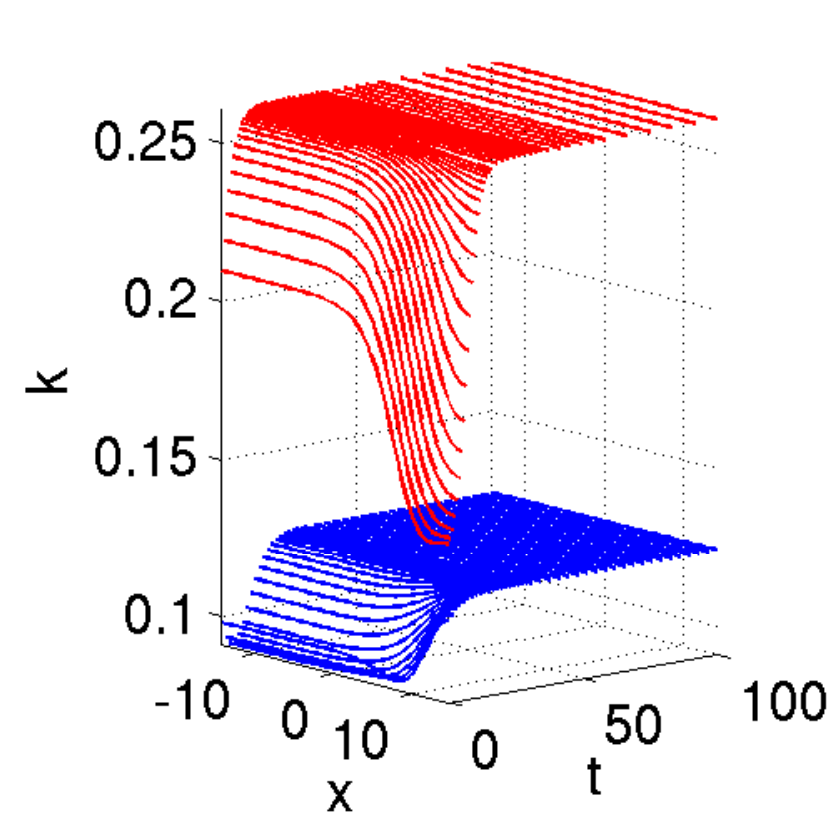}}
\end{tabular}
\caption{Example outputs from {\tt skibademo.m}. For the 
($\al$--dependent) initial states 
from the path from {\tt p1/pt68} to the FSC from  cpdemo1D (blue curve in (a)), we compute the CPs to the FSM, yielding the red curve in (a). 
Near $\al=0.43$ the two CPs have the same values, which makes these 
initial states a so--called Skiba candidate. (b) shows the two 
associated CPs of equal value. \label{slfig3}}}
\ece
\end{figure}

In {\tt skibademo.m} (Listing \ref{l66}, which is a rather elaborate application of the OC facilities of \pdep, and can be skipped on first reading) we find paths from these initial states from {\tt cpdemo1D.m} to the other flat steady state 
with the SPP, namely FSM, and compare the values with the values 
of the paths to FSC, stored in {\tt p.hist.vv}. 
See Fig.~\ref{slfig3} for illustration.


\hulst{caption={{\small {\tt \dname/skibademo.m}. Using the 
data  stored in {\tt p.hist.*}  in Cell 3 of  {\tt cpdemo1D.m} (for $b=0.65$), we compute canonical paths from the initial states 
to the FSM at $b=0.65$, and compare the objective values with those stored 
in {\tt p.hist.vv} (for the path to the FSC). If both are sufficiently close, then we have a good approximation of a Skiba point. Thus, in line 6 we start a loop over the $\al$ values generated in {\tt cpdemo1D.m}; in line 8 we put the associated initial states into {\tt p2.opt.s0}, i.e. overwrite the loaded point, and then (line 10) find the canonical path to the FSM. In line 14 we check if we found a Skiba point approximation, in which case we plot both paths.
Also the second cell then deals with plotting. }},
label=l66,language=matlab,stepnumber=5, linerange=1-21}
{\dhome/sloc/skibademo.m}

\brem{\rm 
The directory {\tt slocdemo} also contains the script files 
{\tt bdcmds2D.m} and {\tt cpdemo2D.m}, used to compute CSS and canonical 
paths for \reff{slcan2} over the 2D domain 
$\Om=(-L,L)\times (-\frac L 2, \frac L 2)$ (based on exactly the same 
init file {\tt slinit.m}), and some modified plotting functions 
{\tt plotsolf.m} and {\tt plotsolfu.m}, see, e.g., \cite[Fig.~4,5]{GU17} 
for some 2D results. 
}\eex 
\erem

\subsection{Optimal harvesting patterns in a vegetation model}\label{vegoc-sec}
\def\dname{vegoc}
Our second example, from \cite{U16}, 
concerns the optimal control of a reaction diffusion 
system used to model harvesting (or grazing by herbivores) in a system for biomass (vegetation) $v$ and soil water $w$, following \cite{BX10}. 
Denoting the harvesting (grazing) effort as the control by $E$, we consider 
\begin{subequations}\label{v1}
\hual{V(v_0,w_0)&=\max_{E(\cdot,\cdot)} J(v_0,w_0,E), \\
\pa_t v&=d_1\Delta v+[gwp^\eta-d(1+\del v)]v-H,\\
\pa_t w&=d_2\Delta w+R(\beta+\xi v)-(r_u v+r_w)w, 
}
with harvest $H=v^\al E^{1-\al}$, and current value 
objective function $J_{c}=J_c(v,E)=pH-cE$, which thus depends on the price $p$, 
the costs $c$ for harvesting/grazing, and $v$, $E$ in a classical 
Cobb--Douglas form with elasticity parameter $0<\al<1$. 
Furthermore, we have the boundary conditions and 
initial conditions 
\huga{\label{nbc}
\pa_\nu v=\pa_\nu w=0 \text{ on } \pa \Om,\quad (v,w)|_{t=0}=(v_0,w_0).  
}
\end{subequations}
Again 
we want to maximize the discounted profit 
\huga{\label{vegJ}
J=\int_0^\infty \er^{-\rho t}J_{ca}(v,E)\dd t. 
}
For the modeling, and the meaning and values of the parameters 
$(g,\eta,d,\delta,\beta,\xi, r_u,r_w, d_{1,2})$ we refer to 
\cite{BX10,U16}, and here only remark that the model aims at a 
realistic description of certain semi--arid systems, that, e.g., 
the discount rate $\rho=0.03$ is in the pertinent economic regime, 
and that, like in most studies of semi--arid systems, 
we take the rainfall $R$ as the main bifurcation parameter. 

Denoting the co-states by $(\lam,\mu)$ we obtain the canonical system 
\begin{subequations}\label{vegcs}
\hual{
\pa_t v&=d_1\Delta v+[gwp^\eta-d(1+\del v)]v-H,\\
\pa_t w&=d_2\Delta w+R(\beta+\xi v)-(r_u v+r_w)w,\\
\pa_t\lam&=\rho\lam-p\al v^{\al-1}E^{1-\al}-
\lam\bigl[g(\eta+1)wv^\eta-2d\del v-d-\al v^{\al-1}E^{1-\al}]
-\mu(R\xi-r_u)w-d_1\Delta \lam, \\
\pa_t\mu&=\rho\mu-\lam g v^{\eta+1}+\mu(r_u v+r_w)-d_2\Delta \mu, 
}
where $\ds E\left(\frac{c}{(p-\lam)(1-\al)}\right)^{-1/\al}v$. 
\end{subequations}

The system \reff{vegcs} has a similar structure as \reff{slcan2}, 
with the immediate difference that 
\reff{vegcs} has four components and many parameters, 
and thus looks somewhat complicated. However, it is still 
convenient 
to implement in \pdep, and leads to many patterned {\em optimal} 
steady states, see \cite{U16} for further discussion. 
Thus, besides documenting the implementation of \reff{vegcs} 
underlying the results in \cite{U16}, our aim here is 
to illustrate that also rather complicated systems can be implemented 
and studied in the \pdep\ OC setting in a simple way.  Writing  
\reff{vegcs} as 
$\pa_t u=-G(u)$, $u=(v,w,\lam,\mu)$, we basically 
need to set up the domain, $G$ and the BCs, and the objective 
function. Table \ref{vegoctab} lists and 
comments on the scripts and functions 
in {\tt ocdemos/vegoc}. The implementation of \reff{vegcs} follows 
the general \pdep\ settings with 
the OC related modifications already explained in \S\ref{slsec}, and 
thus we only give a few remarks in the Listing captions. 

\taskip
\begin{table}[ht]\caption{{\small Scripts (1D) and functions in {\tt ocdemos/vegoc}; 
oosetfemops and veginit (which also contains the parameter values) 
as usual; 
the bottom part contains 'helper' functions for plotting.}}\label{vegoctab}
{\small 
\begin{tabular}{p{0.2\tew}|p{0.73\tew}} 
script/function&purpose,remarks\\\hline
bdcmds, cpcmds&scripts to compute CSS and canonical paths\\
efu, vegjcf&functions to compute [E,H] from u, and the current value 
$J_c$ (by calling efu)\\
\hline 
vegsolplot, vegdiagn&functions to 
plot CPs, and compute and plot diagnostics for CPs\\
vegcm.asc, watcm.asc&colormaps for CP plots (and state plots in 2D)\\ 
\end{tabular}
}
\end{table}
\teskip

\hulst{label=l15,language=matlab,stepnumber=5, firstnumber=1}
{\dhome/vegoc/efu.m}

\hulst{caption={{\small {\tt \dname/efu.m} and {\tt vegsG}. 
{\tt efu} computes the harvesting effort (control) $E$, the 
harvest $h$ and the 
current value $J$, hence also called in {\tt \dname/vegjcf.m}. 
In {\tt vegsG} we first extract the parameters and 
the solution components, and compute $H$ (from {\tt efu}). 
Then we implement the 'nonlinearity' in a straightforward way, and 
compute the residual $G$ using the preassembled stiffness and 
mass matrix (see {\tt oosetfemops}. }},
label=l15b,language=matlab,stepnumber=5, firstnumber=1}
{\dhome/vegoc/vegsG.m}

\begin{figure}\small
(a)\hspace{78mm}(b)\hspace{45mm}(c)\\[-0mm]
\ig[width=78mm]{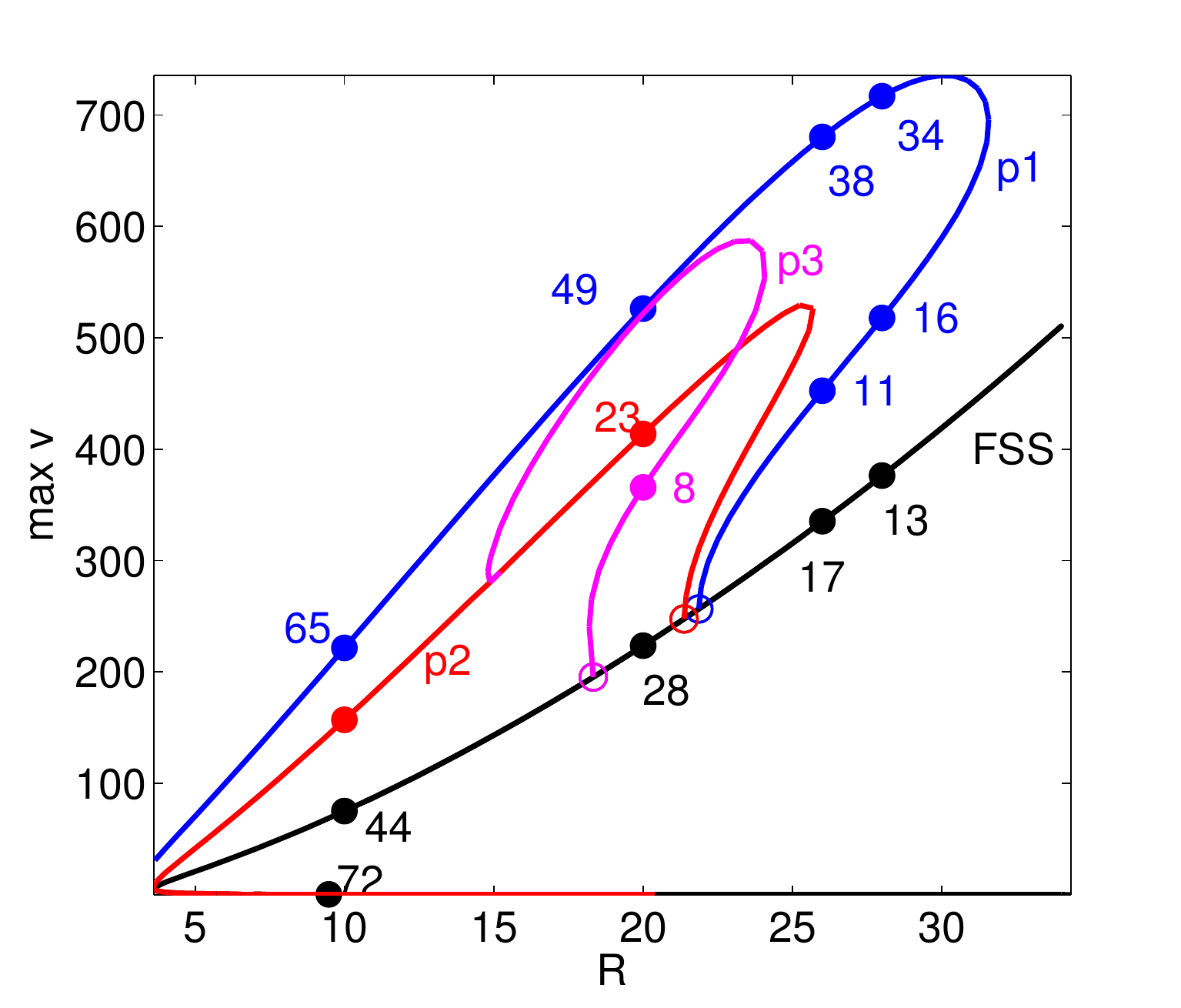}\quad\ig[width=52mm]{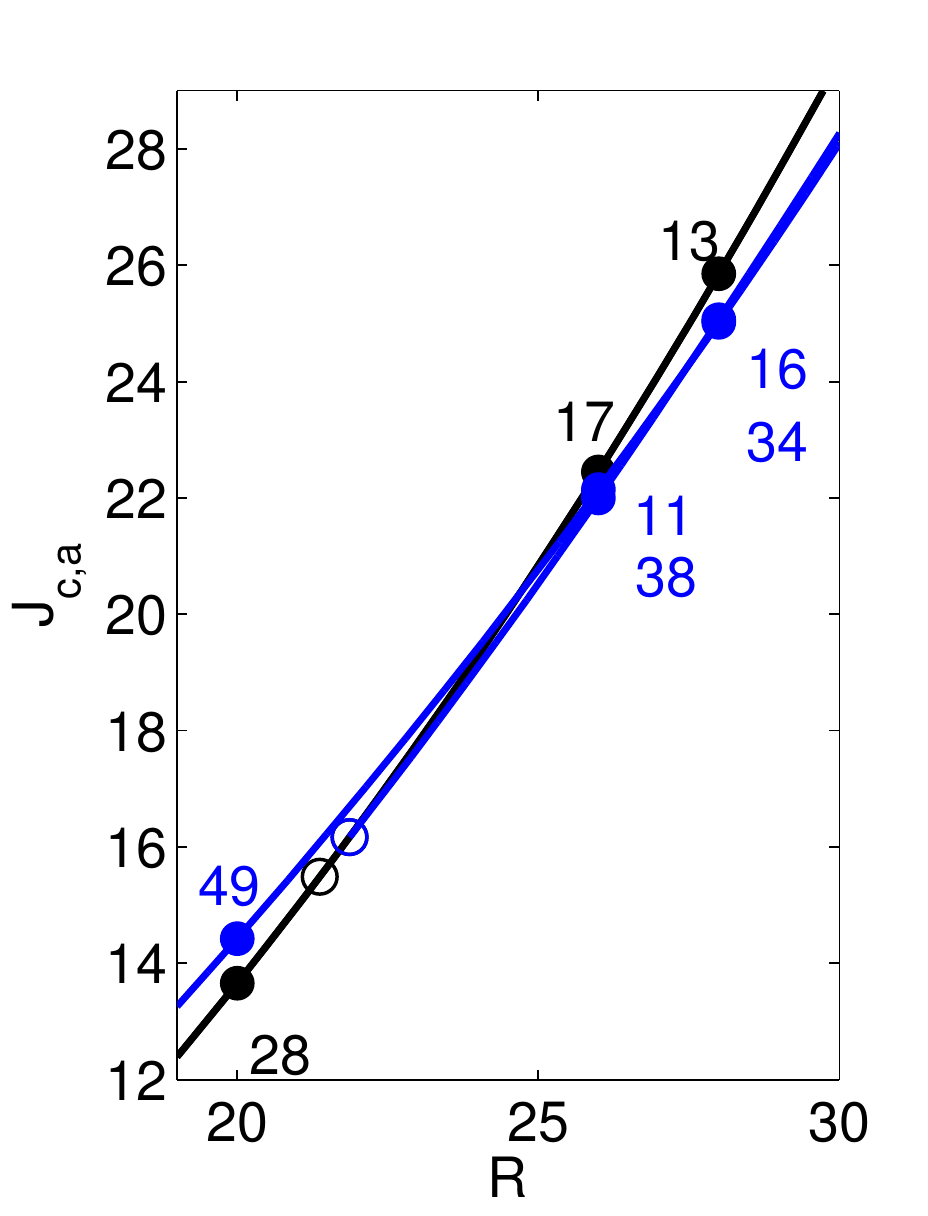} 
\raisebox{30mm}{\begin{tabular}{l}
\ig[width=39mm]{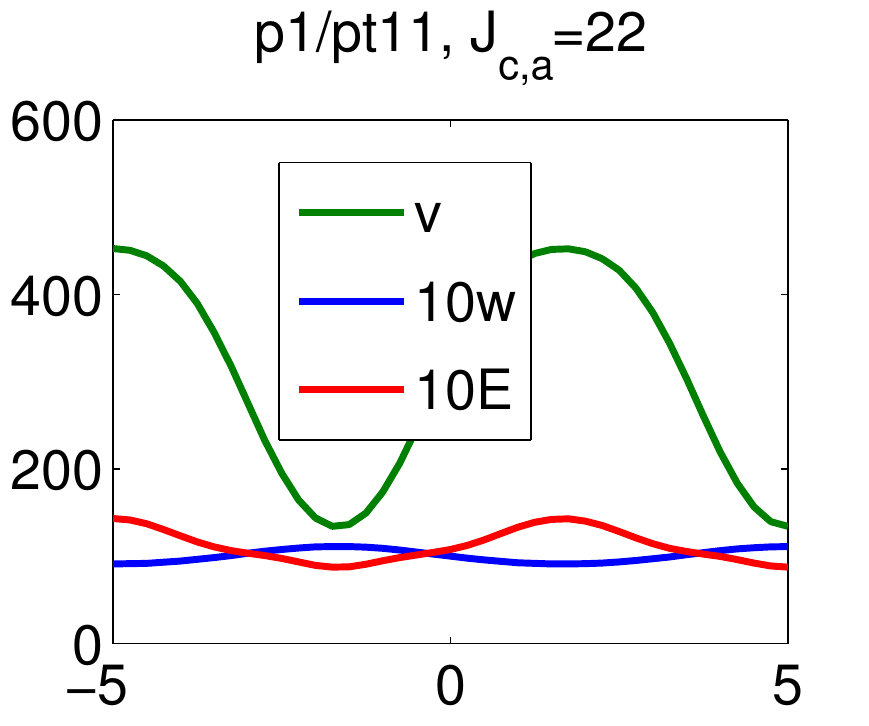}\\\ig[width=39mm]{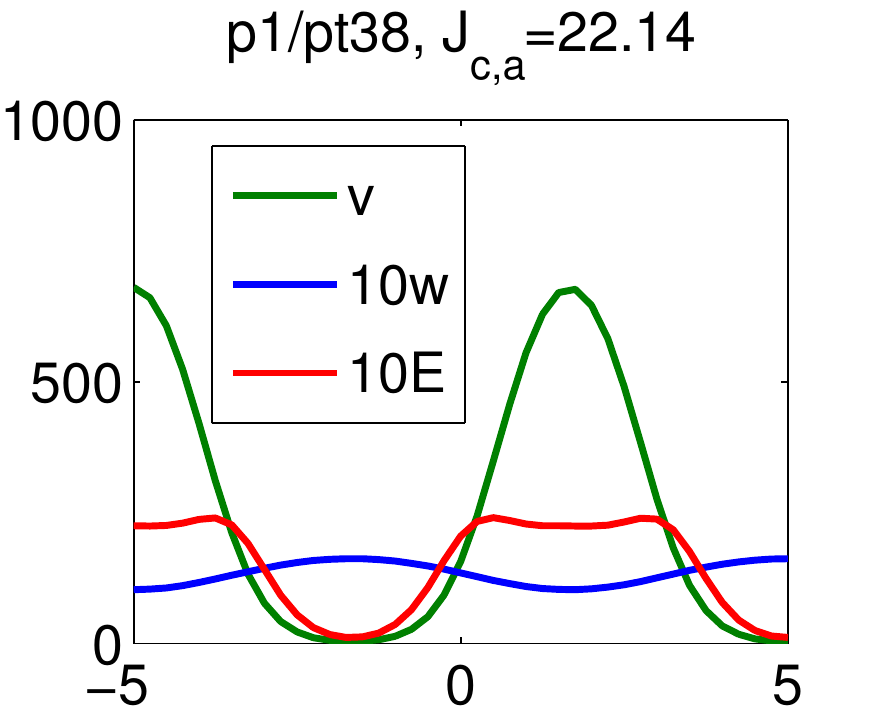}\end{tabular}}\\
\caption{{\small Example outputs of {\tt bdcmds.m}. (a),(b) bifurcation diagrams of CSS in 1D; (c) example solutions. }\label{vf1}}
\end{figure}

Figure \ref{vf1} shows a basic bifurcation diagram of CSS in 1D 
with $\Om=(-L,L)$, $L=5$, from the script file {\tt bdcmds.m}, which 
follows the same principles as the one for the SLOC demo. Again we start 
with a spatially flat (i.e., homogeneous) canonical steady FSS (black), 
on which we find a number of Turing-like bifurcations. 
The blue branch in (a) represents 
the primary bifurcation of PCSS (patterned canonical steady states), which for certain $R$ have the SPP, and, moreover, turn out to be POSS (patterned 
{\em optimal} steady states). See \cite{U16} for more details, 
including a comparison with the uncontrolled case of so called 
``private optimization'', and 2D results for 
$\Om=(-L,L)\times (-\sqrt{3}L/2,\sqrt{3}L/2)$ yielding various POSS, 
including hexagonal patterns. 

\begin{figure}[ht]\small
\begin{tabular}{ll}
(a)&(b)\\
\hs{-5mm}\ig[width=50mm]{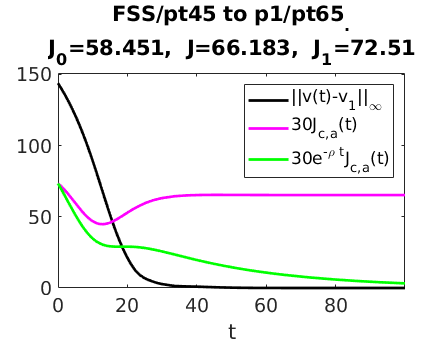}&\hs{-5mm}\ig[width=43mm]{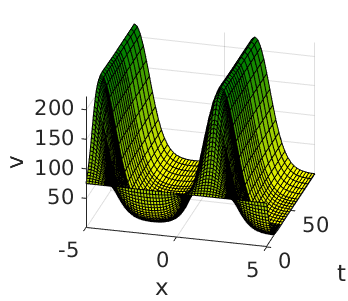}
\ig[width=43mm]{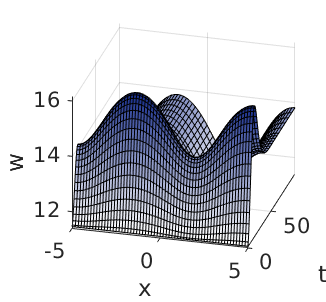}\ig[width=43mm]{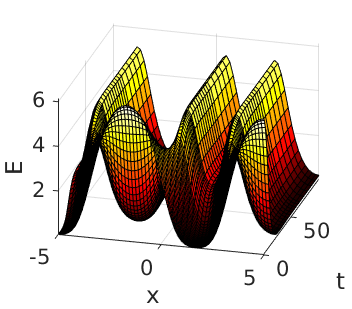}
\end{tabular}
\\[-2mm]
\caption{{\small Example output of {\tt vegoc/cpcmds.m}, namely the 
canonical path from the (states of the) 
lower FCSS(FSS/pt45) (cf.~Fig.~\ref{vf1}) to the PCSS 
(p1/pt65) at $R=10$. 
(a) shows the convergence behavior, the 
current value profit, and obtained objective value. (b) 
show $(v,w)$ and the harvesting strategy $E$. In particular, 
the values $J_0$ of the starting CSS, $J$ of the CP, and $J_1$ of 
the target CSS show that here controlling the system from the 
flat CSS to the PCSS yields a significantly higher value. 
See \cite{U16} for further comments and more details. \label{vf2}}}
\end{figure}

The script files {\tt cpcmds.m} for CPs, and {\tt skibacmds.m} for a Skiba 
point between the flat optimal steady state {\tt FSS/pt13} and the 
POSS {\tt p1/pt34}, 
again follow the same principles as in the the SLOC demo. See Figure \ref{vf2}
for an example output. We use 
customized colormaps for vegetation (green) and water (blue), 
which are provided as vegcm.asc, watcm.asc and whitecm.asc, respectively. 

\subsection{Optimal boundary catch as an example of boundary 
control}\label{bdoc-sec}
\def\dname{lvoc} Our third example, taken from \cite{GUU19}, considers
the optimization of the discounted fishing profit
$J=\int_0^\infty \er^{-\rho t} J_c(v(0,t),\k(t))\dd t$,
$\ds J_c(v,\k)=\sum_{j=1}^2 p_j h_j(v_j,\k_j)-c_j \k_j$.  Here
$v=(v_1,v_2)$ are the populations of two fish species ($v_1=$prey,
$v_2=$predator) in a (1D) lake or ocean $\Om=(0,l_x)$,
$\k=(\k_1,\k_2)$ are the fishing (harvesting) efforts (controls) of
$v_1$ and $v_2$, respectively, {\em at the shore}, and $p_{1,2}$ and
$c_{1,2}$ are the prices for the fishes and the costs for fishing,
respectively, and we again choose a Cobb--Douglas form for the
harvests $h_j(v_j,\k_j)=v_j^{\al_j} \k_j^{1-\al_j}$, with parameters
$\al_j\in(0,1)$. We assume that the fish populations evolve according to a
standard Lotka-Volterra model, namely \huga{\begin{split}
    \pa_t v_1&=d_1\Delta v_1+(1-\beta v_1-v_2)v_1,\\
    \pa_t v_2&=d_2\Delta v_2+(v_1-1)v_2,
	\end{split} \quad\text{or $\pa_t v=-G_1(v)=D\Delta v+f(v)$,} \label{lv1}
}
in $\Om$, with $D=\bpm d_1&0\\0&d_2\epm$
and growth function
$f(v)=\bpm (1-\beta v_1-v_2)v_1 \\  {(v_1-1)v_2} \epm$, with parameter $\beta>0$. The controls $\k_{1,2}$ (fishing efforts for species $v_{1,2}$, 
respectively) occur in the BCs for \reff{lv1}, i.e., we 
assume the Robin BCs 
\huga{
	d_j\pa_n v_j=-g_j:=-\ga_j h_j,\quad j=1,2, \text{ at $x=0$}, 
} 
and zero flux BCs $\pa_n v_j=0$ at $x=l_x$.  Thus, in contrast to the 
{\tt sloc} and {\tt vegoc} examples we no longer have a 
spatially distributed control, but a (two component) boundary control. 

Introducing the co-states 
$\lam_{1,2}:\Om\ra\R$, \PMAXP\ 
yields the evolution and the BCs of the co-states (combining
with~\reff{lv1}, to have it all together)
\begin{subequations}\label{csv}
\hual{
&\left.\barr{l}
\pa_t v=D\Delta v+f(v), \\
\pa_t \lam=\rho\lam-D\Delta\lam-(\pa_vf(v))^T\lam
\earr\right\}\quad\text{ in $\Om=(0,l_x)$}, \\
&\left.\barr{l}
D\pa_n v+g=0, \\
D\pa_n \lam+\pa_v g(v)\lam-\pa_v J_c=0,
\earr\right\}
\quad \text{ on the left boundary $x=0$}, \\
&\left.\barr{l}
D\pa_n v=0, \\
D\pa_n \lam=0,
\earr\right\}
\quad\text{ on the right boundary $x=l_x$},
}
and 
\huga{
\k_j=\left(\frac{(1-\al_j)^2(p_j-\ga_j\lam_j)}{c_j}\right)^{1/\al_j}v_j, 
\text{ evaluated at the left boundary $x=0$},
\quad j=1,2.\label{kjdef}
}
\end{subequations}
See \cite{GUU19} for details on the derivation of \reff{csv}. 

Thus, we again have a 4 component reaction diffusion system (\ref{csv}a-c) for the 
states $v$ and the costates $\lam$, but now the controls live on 
the boundary at $x=0$, leading to nonlinear flux boundary conditions. 
Also, from the modeling point of view, the pertinent questions for 
\reff{csv} are slightly different than for \reff{vegcs}, since for 
\reff{csv} we are not so much interested in  bifurcations and 
pattern formation (which do not occur for the parameters chosen below), 
but rather in the dependence of the (unique) CSS on the parameters, and 
mostly in the canonical paths leading to these CSS. 

In any case, we now focus on how to put the 
nonlinear BCs into \pdep. 
Table \ref{lvoctab} lists the scripts and functions in {\tt ocdemos/lvoc}, 
and some further comments are given in the listing captions below. 

\taskip 
{\small \begin{longtable}{p{0.17\tew}|p{0.8\tew}} 
\caption{{\small Scripts and functions in {\tt ocdemos/lvoc}; the first four follow closely the same \pdep\ principles as the respective functions in 
{\tt sloc} and {\tt vegoc}. {\tt lvsG} is a rather non--generic implementation of the rhs of \reff{csv}, and in particular strongly uses the 1D nature of the problem. The others are naturally also problem specific, and in particular the last three overload some \pdep\ standard functions. } }\label{lvoctab}\endfirsthead
\endhead\endfoot\endlastfoot
script/function&purpose,remarks\\\hline
bdcmds,cpcmds&scripts to compute CSS and CPs\\%
lvinit,oosetfemops&init routine, and setting of FEM matrices \\
lvsG&the rhs for (\ref{csv}a), also explicitly implementing 
the BC (\ref{csv}b), while (\ref{csv}c) are naturally fulfilled with K  
the Neumann Laplacian. \\
lvbra&branch-output, here substantially modified from stanocbra, i.e., writing specific data like profits/harvest/controls per fish-species on the 
		branch \\
hfu&returns harvest $h\in\R^2$, and derivatives (needed 
for {\tt lvsG}) $\pa_v h\in\R^{2\times 2}$, $\pa_k h\in\R^{2\times 2}$; 
straightforward implementation.\\%
lvjcf&returns $J_c$, as required in the standard \pdep\ setup for OC problems, but also the individual profits $J_{1,2}$ per fish-species.\\
jca&$J_c$ average, overloaded here since $J_c$ is on boundary, hence  
the default averaging by $1/|\Om|$ makes no sense\\ 
lvdiagn&diagnostics functions for canonical paths\\
\hline
stanpdeo1D${<}$pde &classdef, overloads the \pdep\ classdef stanpdeo1D in order to have $\Om=(0,l_x)$ instead of $\Om=(-l_x,l_x)$ (standard).\\
plotsol, psol3D&some more overloads of \pdep\ standard functions 
for plotting. 
	\end{longtable}
}
\teskip

\begin{figure}[ht]{\small
\bce\begin{tabular}{ll}
(a) $J,v$ and $\k$ at $x=0$; continuation in $c_1$ 
&(b) CSS at $c=(0.1,0.1)$\\
\ig[width=33mm, height=44mm]{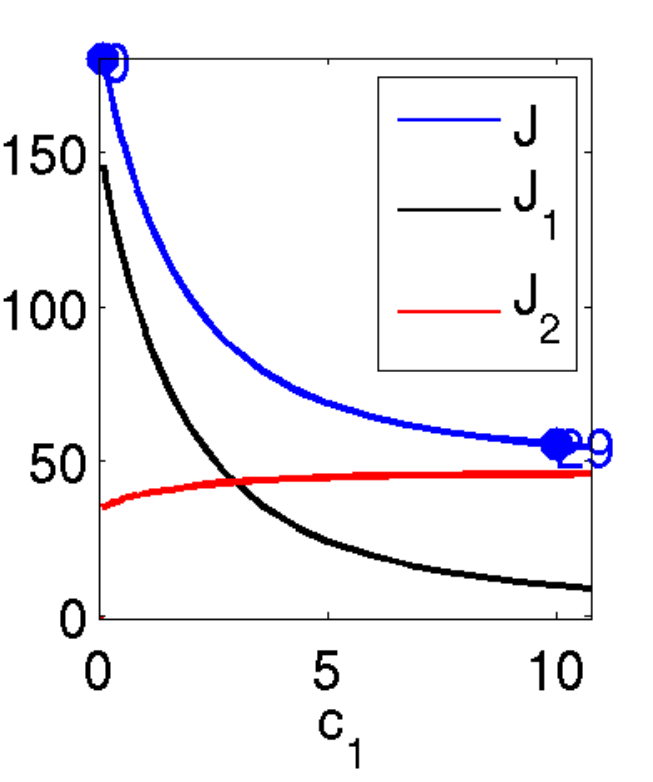}\ \ 
\ig[width=33mm, height=44mm]{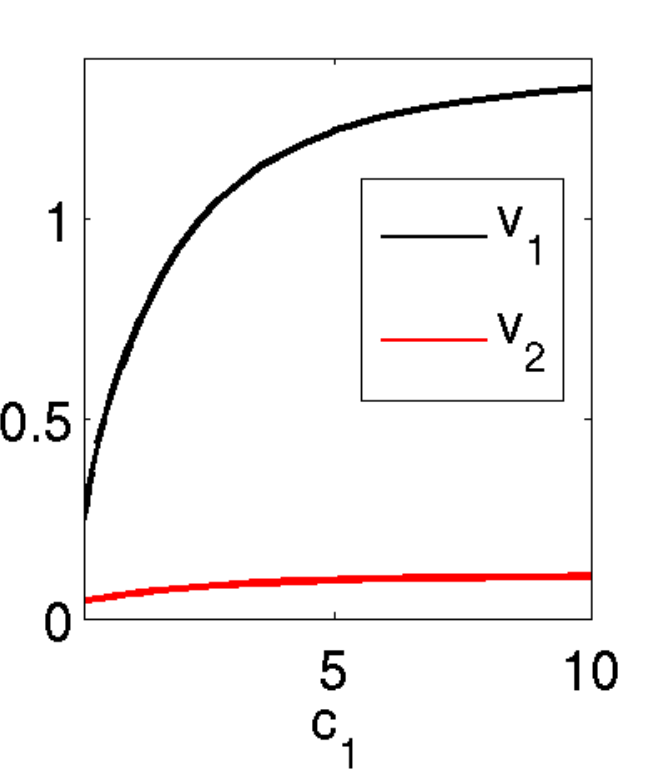}\ \ 
\ig[width=33mm, height=44mm]{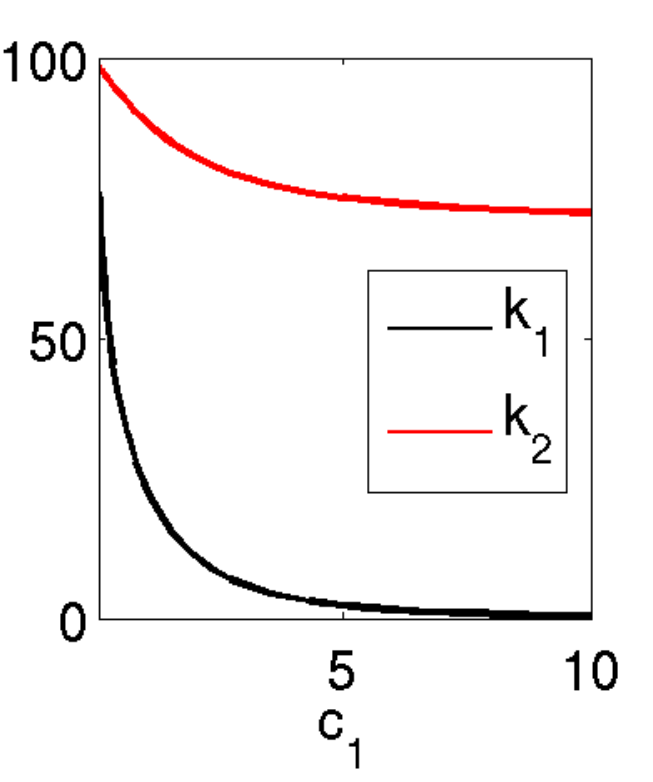}&
\hs{-4mm}\ig[width=75mm, height=53mm]{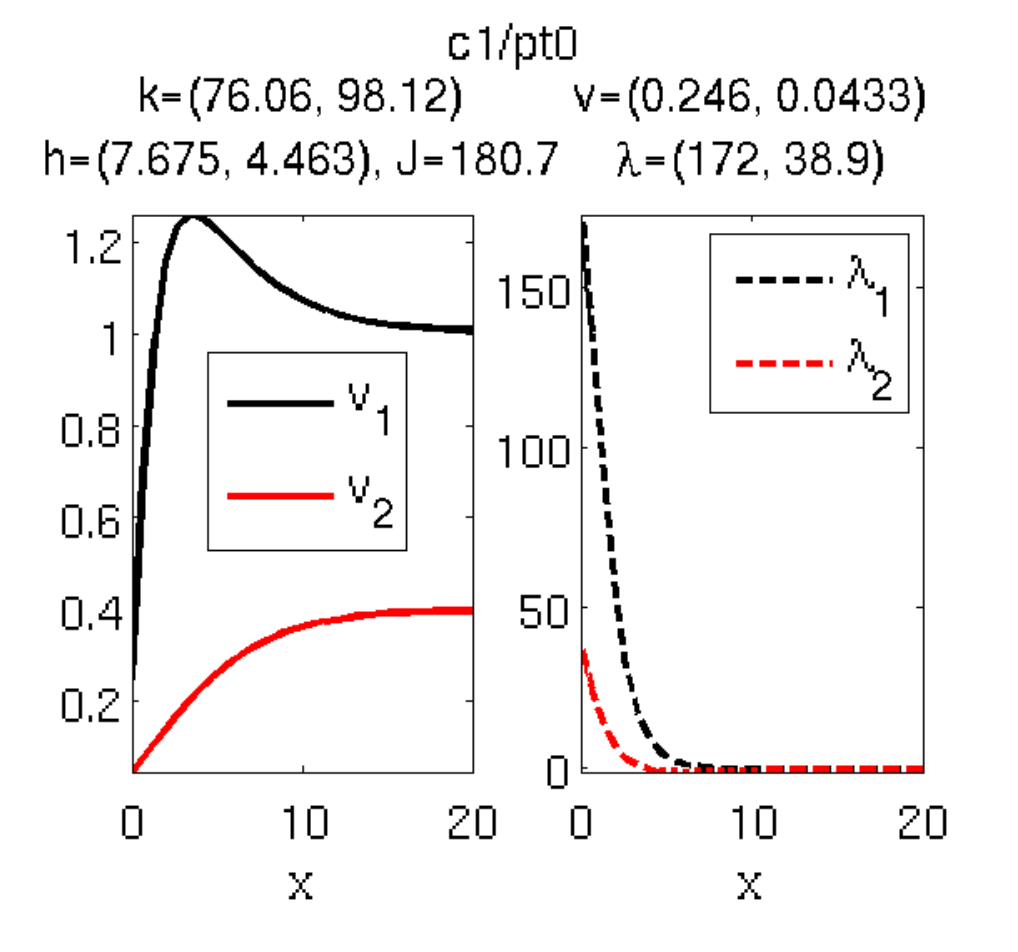}
\end{tabular}
\ece}
\vs{-6mm}
\caption{{\small (a) continuation diagrams in $c_1$ (costs for prey fishing); $J_j=p_jh_j-c_j\k_j$, $J=J_1+J_2$. (b) An example CSS plots.  \label{lvf1}}}
\end{figure}

\begin{figure}[ht]{\small
\bce\begin{tabular}{l}
\hs{-3mm}
\ig[width=48mm, height=48mm]{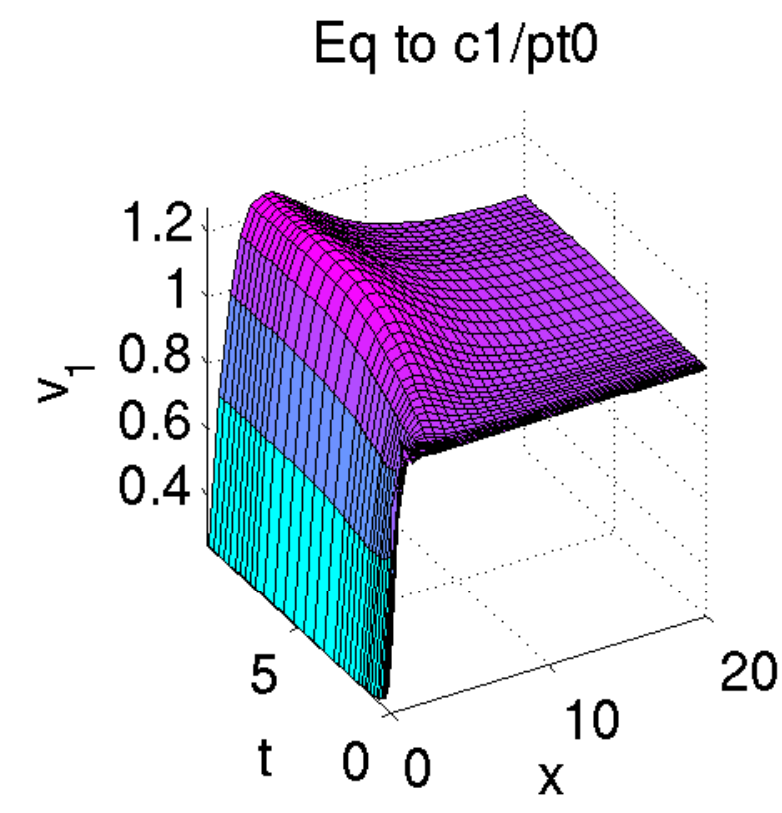}
\hs{-3mm}\ig[width=45mm, height=48mm]{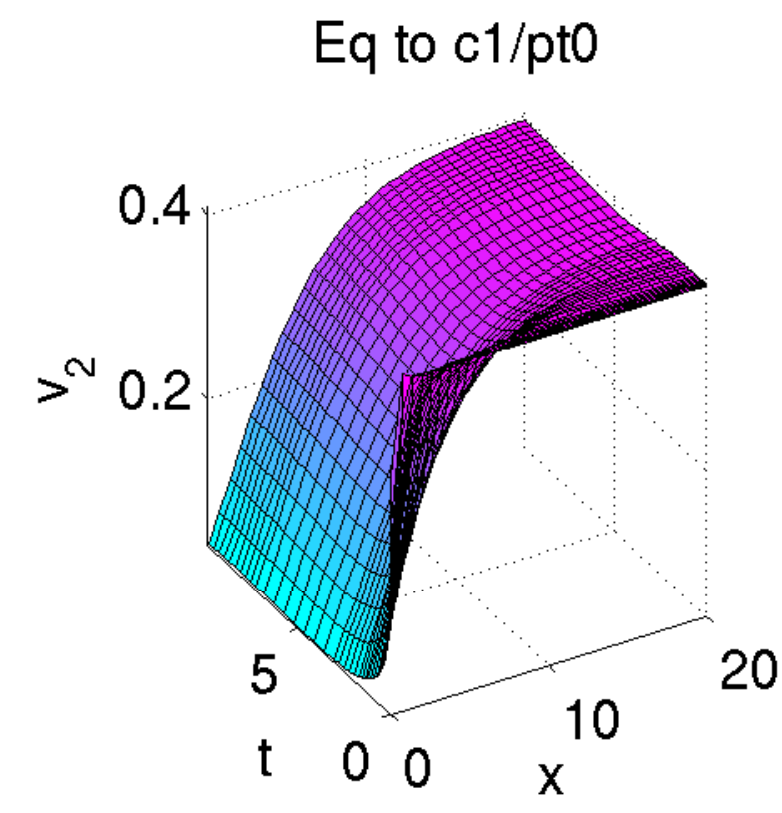}\hs{3mm}
\ig[width=28mm, height=44mm]{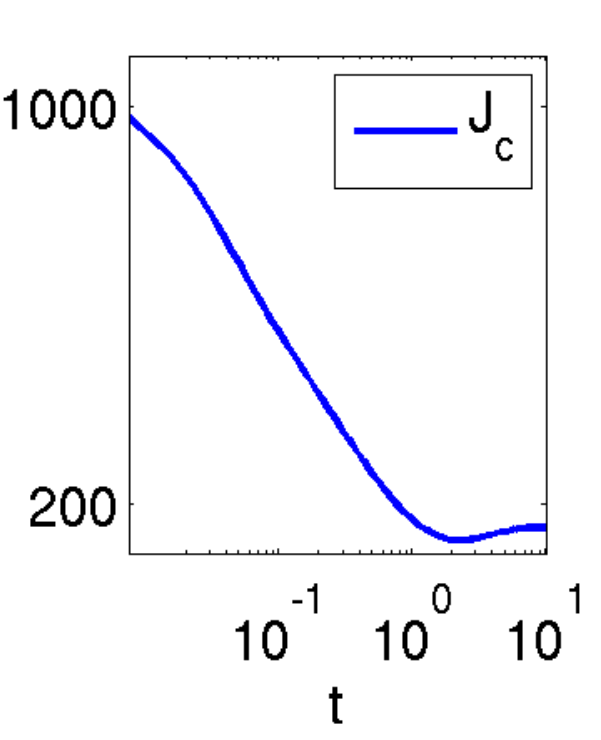}\ig[width=28mm, height=44mm]{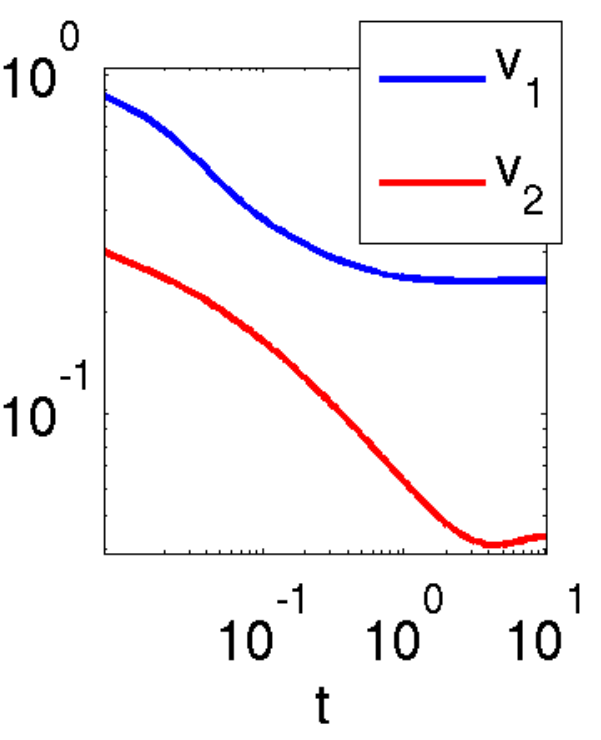}
\ig[width=28mm, height=44mm]{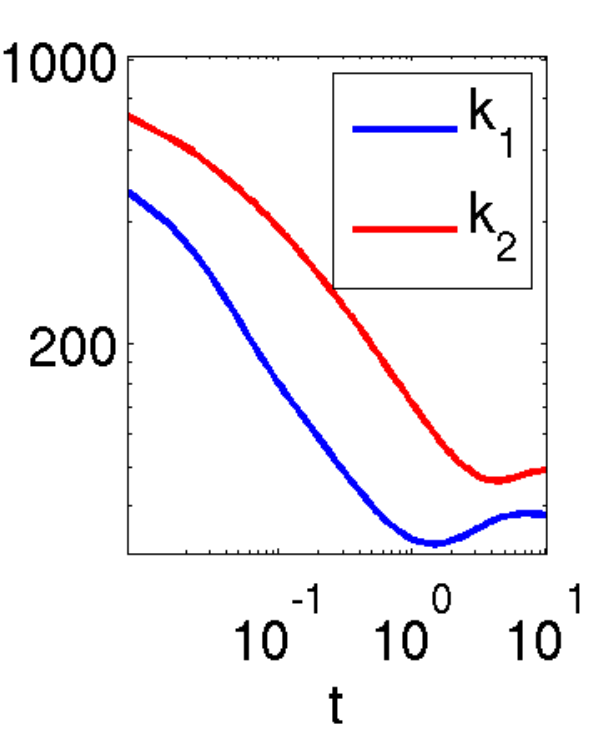}
\end{tabular}
\ece}
\vs{-6mm}
\caption{{\small A CP starting from the spatially homogeneous steady 
state $V^*$ of \reff{lv1} to the CSS at $(c_1,c_2)=(0.1,0.1)$. 
Note the logarithmic scales in the time-series of the values at the left boundary. \label{lvf2}}}
\end{figure}

Figure \ref{lvf1} gives some example plots from running {\tt bdcmds.m}. 
This script is rather lengthy, with the main point here to illustrate how to plot various data of 
interest by first putting it on the branch (here via {\tt lvbra.m}) and then 
choosing the pertinent component for plotting.  
The script {\tt cpcmds.m} for computing canonical paths 
follows the same outline as those for the {\tt sloc} and {\tt vegoc} examples. 
Figure \ref{lvf2} 
shows a canonical path  to the CSS at $c=(0.1,0.1)$, starting from the homogeneous fixed point $V^*=(1,1-\beta)$ of \reff{lv1}. As already said, we refer to \cite{GUU19} for discussion of the results. 

\hulst{caption={{\small {\tt \dname/lvsG.m}. As usual we first 
extract the relevant parameters and fields from $u$, and compute the 
'bulk' nonlinearity $f$ (i.e., the nonlinearity in $\Om$). Then, to compute 
the BCs we first extract values of the co-states on the boundary, and 
the associated harvests and their derivatives. These are needed to 
compute the BCs $g_j$, $j=1,\ldots,4$ (lines 13,14), which can then 
directly be added to the rhs on the boundary (lines 16,17). Here we strongly use the 1D setup, i.e., that the left boundary value of component $j$ is at 
u((j-1)*np+1), where np is the number of spatial discretization points. 
In line 19 we then assemble the system $K$. Compared to the setup in {\tt vegoc/vegsG} this has the advantage that the diffusion constants $d_{1,2}$ 
can be used like any other parameter at this point. The actual computation 
of r in line 20 then works as usual.  }},
label=33,language=matlab,stepnumber=5, firstnumber=1}
{\dhome/lvoc/lvsG.m} 

\section{Examples and implementation details for CPs to CPSs}
\label{cps-sec}
The basic idea for the computation of CPSs and of CPs to CPSs  
is similar to that for the computation of CSSs and CPs to CSSs. 
We first search for CPSs, usually via Hopf bifurcations from
CSSs, and then aim to 
compute CPs $u$ to such CPSs with the SPP, again using the main 
\pdep\ OC user interface {\tt isc}, which now implements 
Algorithm \ref{CPSalg}.  
To illustrate the setup we first discuss an ODE toy model.  
Subsequently we come to a PDE model for pollution mitigation. 

\subsection{An ODE toy problem}
We start with the ODE toy model
\begin{subequations}\label{ode}
\hual{& \dot{x_1} = \rho \left(-x_1 - \frac{\theta x_2}{\rho} + x_1 y_1 r^2 \right), \quad \dot{x_2} = \rho \left( -x_2 + \frac{\theta x_1}{\rho} + x_2 y_1 r^2 \right), \\
& \dot{y_1} =  \om y_2, \quad \dot{y_2} = \om \sin(2\pi y_1), 
}
\end{subequations}
with parameters $\rho,\om,\theta>0$ and $ r=\sqrt{x_1^2+x_2^2} $. Although 
\reff{ode} is not derived as a canonical system for an OC problem, 
we interpret $ x=(x_1,x_2) $ as states and $ y=(y_1,y_2) $ as costates 
and call a time periodic solution of \reff{ode} a CPS. 

\subsubsection{Preliminary analytical remarks}
\begin{figure}[ht]
\bce 
\begin{tabular}{lll}
{\small (a)}&{\small (b)}&{\small (c)}\\
\ig[width=0.27\tew]{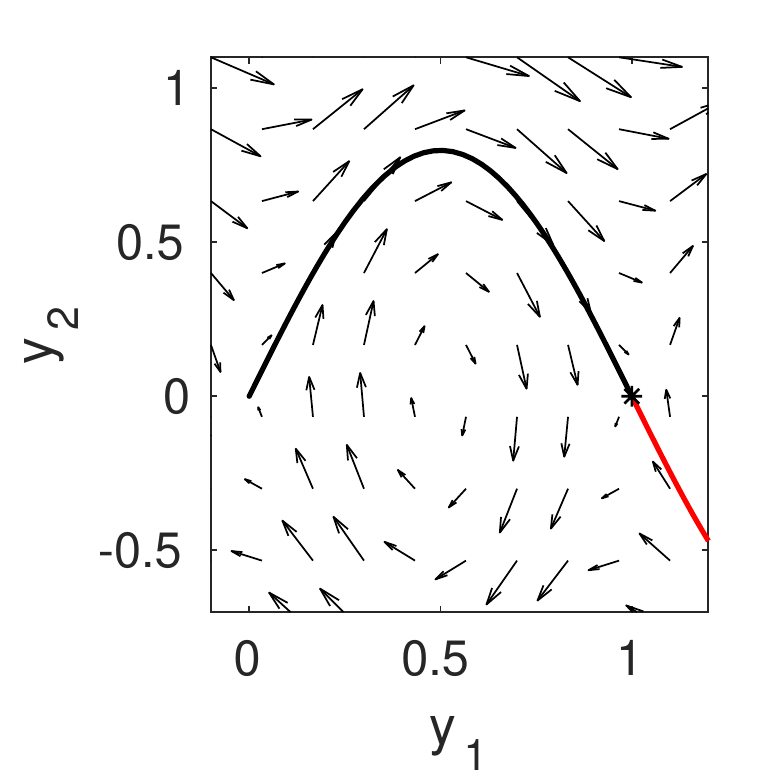}
&\raisebox{5mm}{\ig[width=0.14\tew]{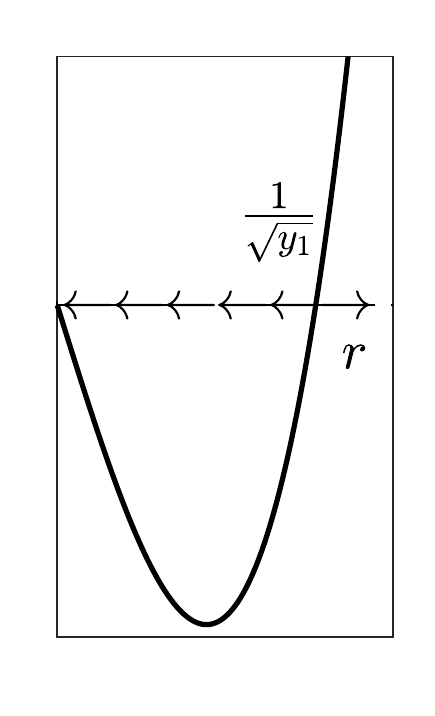}}
&\ig[width=0.3\tew]{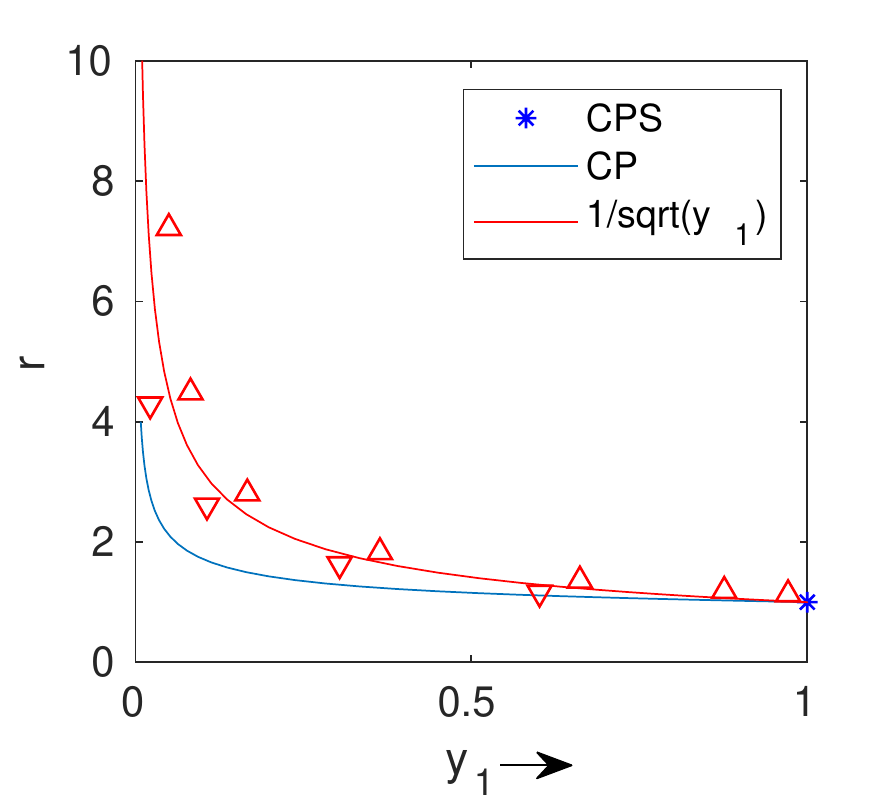}
\end{tabular}
\ece 

\vs{-5mm}
\caption{\small (a) Phase portraits 
of the ``control system'' for $y$ 
with the heteroclinic orbits from $(0,0)$ to $(1,0)$ and from 
$(2,0)$ to $(1,0)$ (red orbit, partial plot). 
(b) The behavior of the $r$ system. 
(c) Sketch of the expected CPs in the $y_1$--$r$--plane 
(using actual numerics for 
parameter values $(\om,\rho,\theta)=(1,1,1)$ and initial state 
$(x_1,x_2)(0)=(4,0)$. The red line shows the nullcline $r=1/\sqrt{y_1}$, 
above (below) which he have $\dot r>0$ ($\dot r<0$). As $r(0)=\|x(0)\|>1$ 
we must start below the red--line by choice of $0<y_1(0)<1$. Together, 
the costates $(y_1,y_2)(0)$ must be on the heteroclinic to $y=(1,0)$ 
such that $y_1(t)$ increases in $t$ in such a way that $r(t)\to 1$ 
as $t\to\infty$. \label{ps}}
\end{figure}

Concerning our points of interest, the model \reff{ode} can almost completely be treated analytically and thus 
can be used to test our numerical methods. For fixed $y_1>0$, 
the nonlinear system (\ref{ode}a) has the unstable periodic orbit 
$r=1/\sqrt{y_1}$ of period $2\pi/\theta$, and (\ref{ode}a) is 
coupled to (or driven by) by the nonlinear pendulum (\ref{ode}b). 
In detail, by polar coordinates in $ \left( x_1,x_2 \right) $, \reff{ode} transforms to
\begin{subequations}\label{toypolar}
\hual{&\dot{r} =\rho( -r+y_1 r^3), \quad 
	\dot{\varphi} = \theta,  \\
	& \dot{y_1} =  \omega y_2, \quad 
	 \dot{y_2} = \omega \sin(2\pi y_1), 
}
\end{subequations}
with $ \varphi = \arg(x_1,x_2) $ and $ r=\sqrt{x_1^2+x_2^2} $, 
with the phase portraits of the $r$ ODE (for fixed $y_1$) 
and the $y$ system sketched in Fig.\ref{ps}(a,b).  Thus, to 
find a CPS we look for CSS of the reduced system 
\begin{align*}
& \dot{r} =\rho( -r+y_1 r^3), \\
&\dot{y_1} =  \omega y_2, \quad 
 \dot{y_2} = \omega \sin(2\pi y_1). 
\end{align*}
The costates are independent of the states and the $\dot{y}$ system has 
the first integral 
$E(y_1,y_2)=\frac{1}{2} y_2^2+\frac{\omega^2}{2 \pi} \cos(2\pi y_1)$,  
i.e.~solutions of the system lie on contour lines of $ E $, 
see Fig.\ref{ps}(a). We choose 
$ y_1 \in \N $, $ y_2=0 $ and $r=\frac{1}{\sqrt{y_1}}$, which 
yields a $2\pi$ periodic solution in the full system. 
We fix $(\rho,\theta)=(1,1)$ and one of these CPS, namely 
\huga{
\uh(t)= \left( r(t),\varphi(t),y_1(t),y_2(t) \right) = 
\left( 1,t,1,0 \right)
}  
as our CPS of interest and aim to compute canonical paths to $\uh$. 

Now given an initial state $(x_1,x_2)(0)$ with, e.g., $\|x\|>1$ and 
aiming at a CP to $\uh$, i.e., 
the CPS associated with $(r,y_1)=(1,1)$, we 
have the situation sketched in Fig.\ref{ps}(c). The only possible 
co-state choice to end in the CPS lies on the
heteroclinic connection from $ (y_1,y_2)=(0,0) $ 
to $ (1,0) $. Thus, $ y_1 $ which is the only
costate which influences the states, can take values in
$\left[0,2\right]$. The state dynamics are sketched by the 
red triangles for $\dot r$  
with $\dot r<0$ ($\dot r>0$) for $r<1/\sqrt{y_1}$ ($r>\sqrt{y_1}$). 
Since $r(0)>1$ we need to choose $y_1(0)>0$ sufficiently small 
such that $r(0)<1/\sqrt{y_1(0)}$ to have $\dot r<0$, and at 
the same time we need $(y_1(0),y_2(0))$ on the black $y$--heteroclinic 
to  $(1,0)$. The argument for 
$ \frac{1}{\sqrt{2}}< \|x_0\| < 1 $ works similarly. 

We can also explicitly compute the Floquet multipliers of $\uh$ and 
the associated projections. In polar coordinates, the variational 
equation \reff{mon123} is autonomous, namely 
\huga{\dot v=J_f(\uh)v, \quad J_f(\uh)=\bpm 
2\rho&0&0&0\\
0&0&0&0\\
0&0&0&\om\\
0&0&2\pi\om&0
\epm, \label{toyvar} 
} 
with eigenvalues $\mu=(0,-\sqrt{2\pi}\om,2\rho,\sqrt{2\pi}\om)$. 
Since \reff{toyvar} is autonomous we obtain the multipliers by 
exponentiation of $\mu$, namely 
\huga{\label{mult}
\ga=(\ga_1,\ga_2,\ga_3,\ga_4)=
(1,\exp(-\sqrt{2\pi}2\pi\om/\theta),\exp(4\pi\rho/\theta), 
\exp(-\sqrt{2\pi}2\pi\om/\theta)).
}
Clearly, additional to the trivial multiplier we 
have one stable multiplier $\ga_2$ and two unstable 
multipliers $\ga_{3,4}$. Similarly, we can also compute 
the projection $P$ onto the center unstable eigenspace in 
Cartesian coordinates (which depends on the target point $\uh_0$) 
analytically. Moreover,  
\huga{\label{ga2}
\text{$\ga_2\nearrow 1$ as $\om\to 0$ or $\theta\to \infty$,}
}
and conversely $\ga_2\searrow 0$ as $\om\to \infty$ or $\theta\to 0$, 
and this (and the analytical projection, 
implemented in a testing function {\tt anaproj}) can 
be used to tune and test the convergence behavior 
of the CPs in the numerics.

\subsubsection{\pdep\ implementation and results}
\begin{figure}[h!]\centering
	\ig[width=4cm, height=4cm]{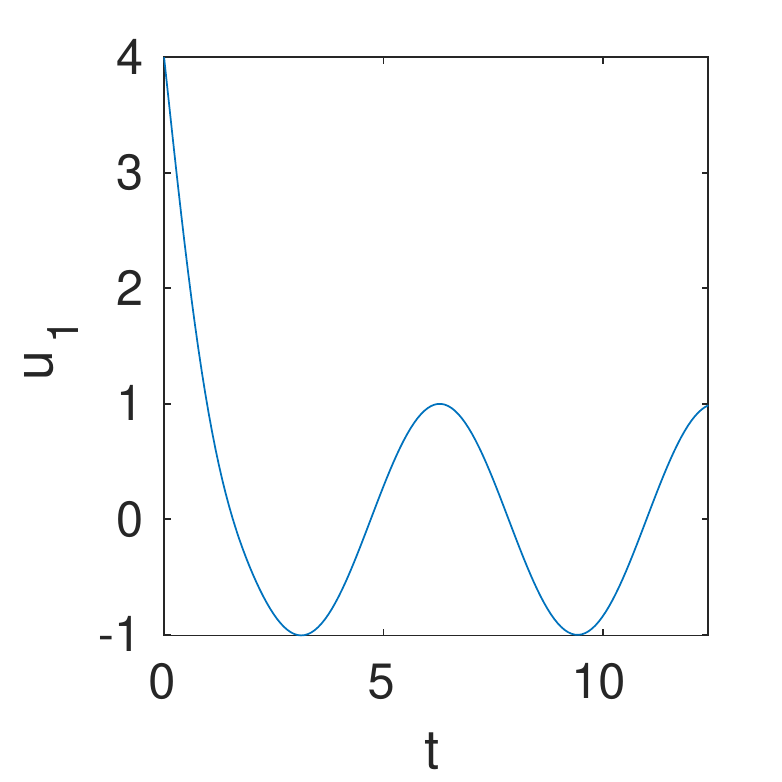}
	\ig[width=4cm, height=4cm]{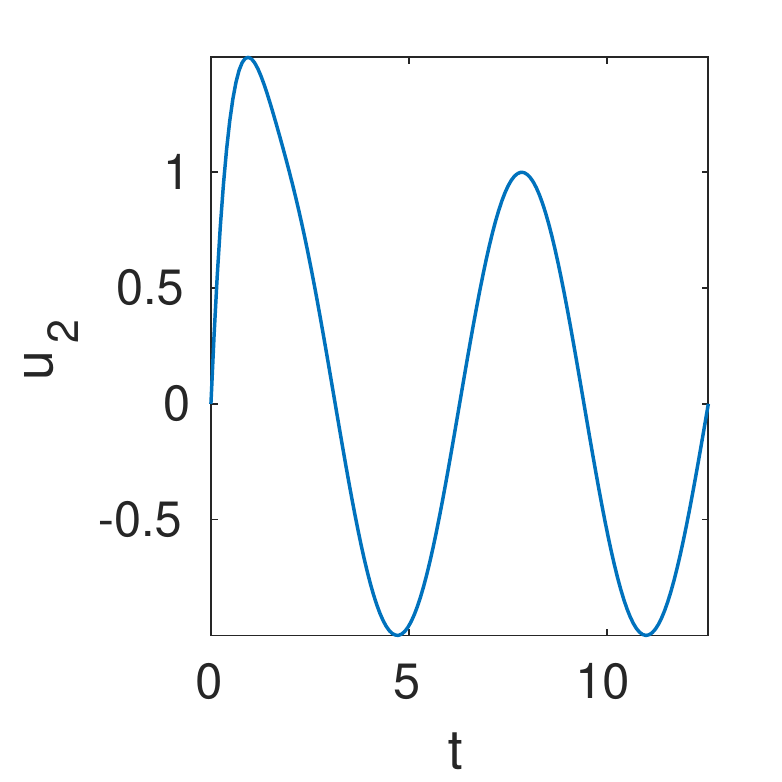}
	\ig[width=4cm, height=4cm]{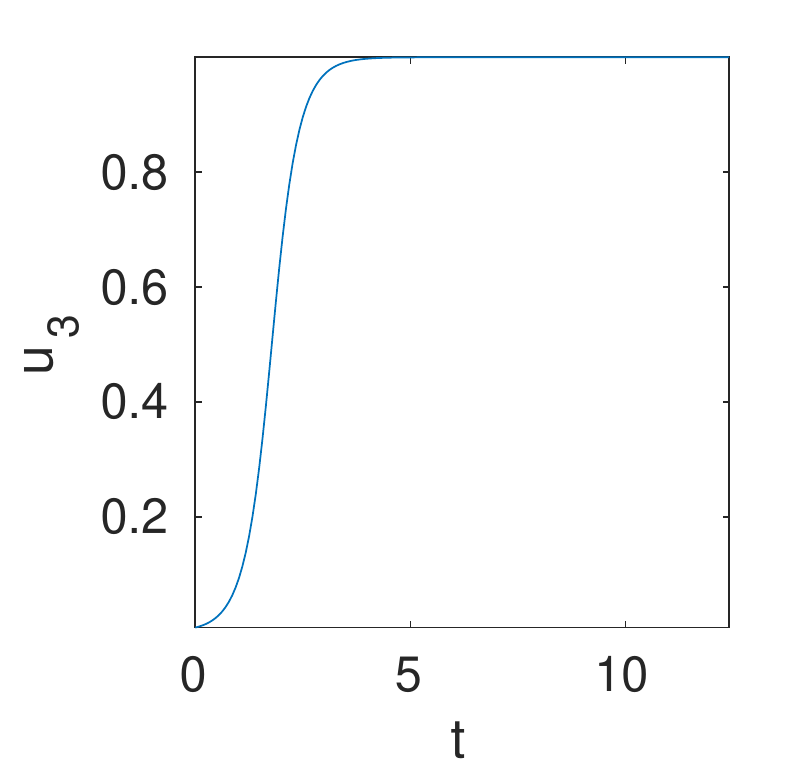}
	\ig[width=4cm, height=4cm]{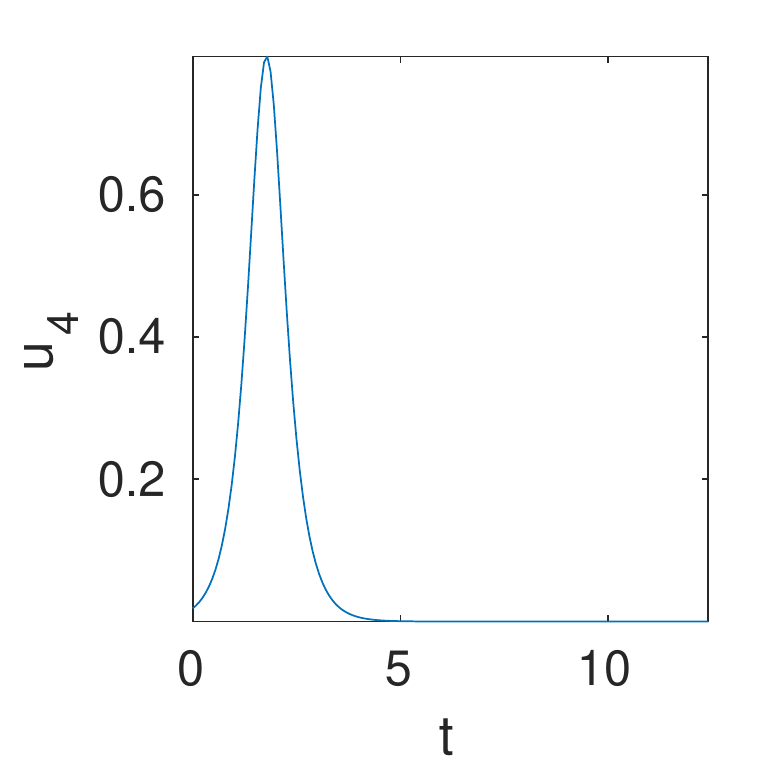}
	\caption{\small Canonical path from $(x_1,x_2)=(4,0)$ to $\uh_0=(1,0,1,0)$ on CPS $\uh$, $(\rho,\om,\theta)=(1,1,1)$. \label{toy_b1}}
\end{figure}

Table \ref{toytab} lists the main files for the implementation of
\reff{ode}. Even though we do not need the bifurcation methods of
\pdep, we implement the ODE \reff{ode} as a \pdep\ problem via the
convenience function {\tt toyinit}, because the OC
routines reuse these basic \pdep\ data structures. To show the setup, 
in {\tt cmds\_basic} (Listing \ref{toy1}), 
 we compute some CPs with easy parameter settings, namely 
$(\rho,\om,\theta)=(1,1,1)$, which yields 
$\ga_2\approx 10^{-7}$ for the ``leading'' Floquet multiplier, 
and hence fast convergence to the CPS. See Fig.~\ref{toy_b1} 
for some basic results, which were also used to generate 
the blue curve in Fig.~\ref{ps}(c). 
For convenience we outsourced the main setup in {\tt  ocinit\_sp}, 
see Listing \ref{toy1b}, which also recalls the meaning of the 
most important parameters. Additionally, there are some helper functions 
for plotting. 

\taskip
\begin{table}[ht]
{\small\begin{tabular}{p{0.15\tew}|p{0.8\tew}}
cmds\_basic &computes canonical paths with easy parameter setting.\\%
cmds\_advanced &canonical paths with advanced parameter setting, i.e. Floquet-multipliers near $1$. \\%
toyinit& init routine; set the limit cycle as a {\tt pde2path} struct with standard parameters. \\%
ocinit\_sp&local extension of {\tt ocinit}, resetting a number of parameters \\
sG, sGjac&rhs side of \reff{ode} resp. the Jacobian\\
anaproj& computes the monodromy matrix analytically 
\end{tabular}}
	\caption{{\small Main scripts and functions in {\tt ocdemos/toy}. }\label{toytab}}
\end{table}
\teskip

\hulst{language=matlab,stepnumber=5, linerange=1-9}
{\dhome/toy/cmds_basic.m}
\hulst{caption={{\small Selection from {\tt ochopftriv/cmds\_basic.m}. 
Cell 1: Initialization. Cell 2: A first canonical path, see Figure \ref{toy_b1} for solution plots. Cell 3: A path which needs to adapt the truncation time {\tt T}. A full circle around the origin needs time $2\pi$ independent of the other states and thus a start in $(4,4,0,0)$ to the fixed end point $(1,0,1,0)$ has a truncation time of $\frac{7}{4}\pi+2\pi\N$. In Cell 4 we shift the target point on the CPS and 
see similar behavior. 
The omitted rest of the script is plotting.}},
label=toy1,language=matlab,stepnumber=5, firstnumber=17, linerange=18-26}
{\dhome/toy/cmds_basic.m}

\hulst{caption={{\small {\tt toy/ocinit\_sp.m}, initialization 
for the computation of canonical paths adapted to this problem. 
First runs {\tt ocinit}, which sets most values to defaults, 
but sets start and end-point manually because no structure 
via bifurcation analysis has been constructed before.  }},
label=toy1b,language=matlab,stepnumber=5, firstnumber=1}
{\dhome/toy/ocinit_sp.m}

\begin{figure}[ht]
\bce
\begin{tabular}{lll}
(a)&(b)&(c)\\
\ig[width=33mm]{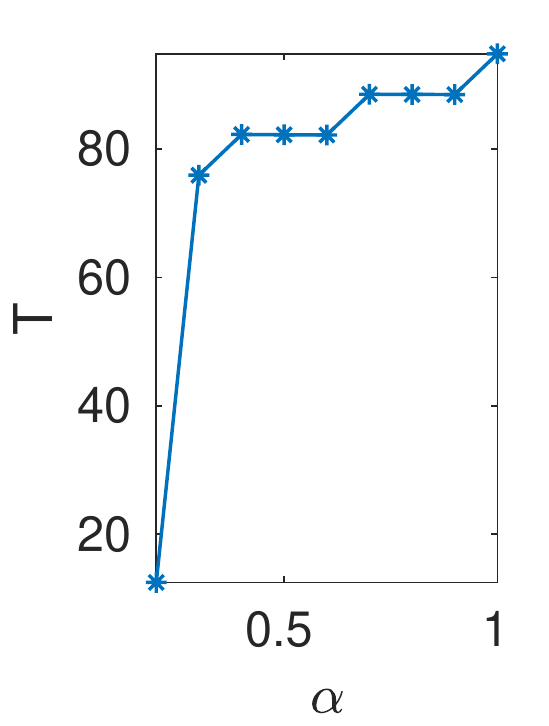}&
\ig[width=43mm]{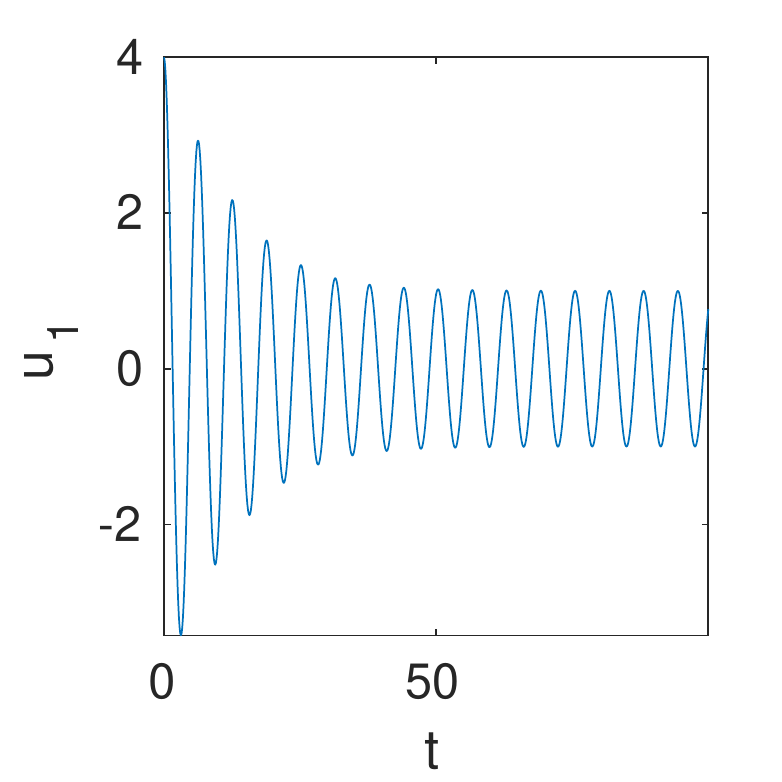}
\ig[width=43mm]{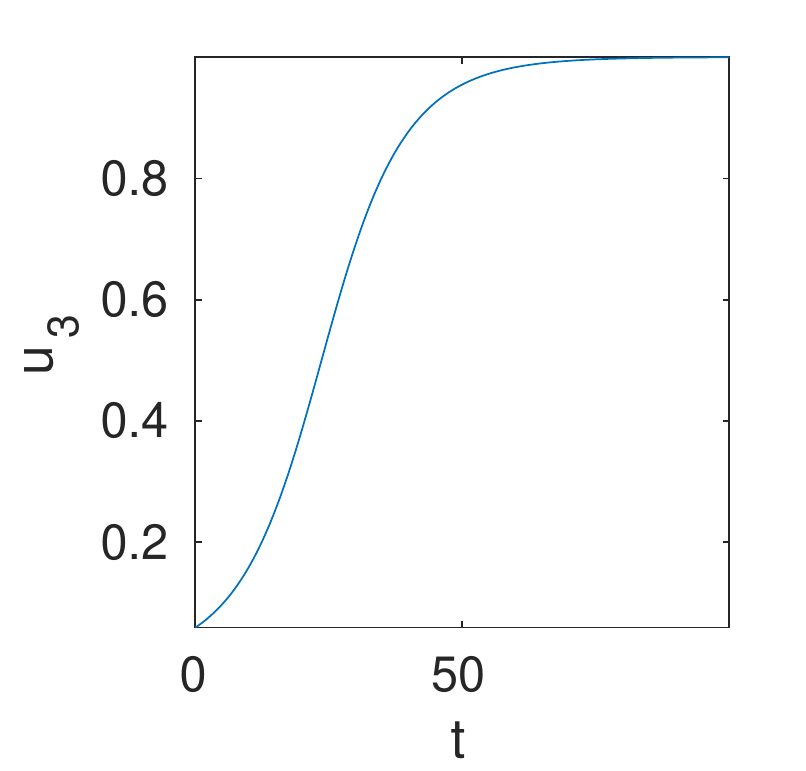}
&\ig[width=43mm]{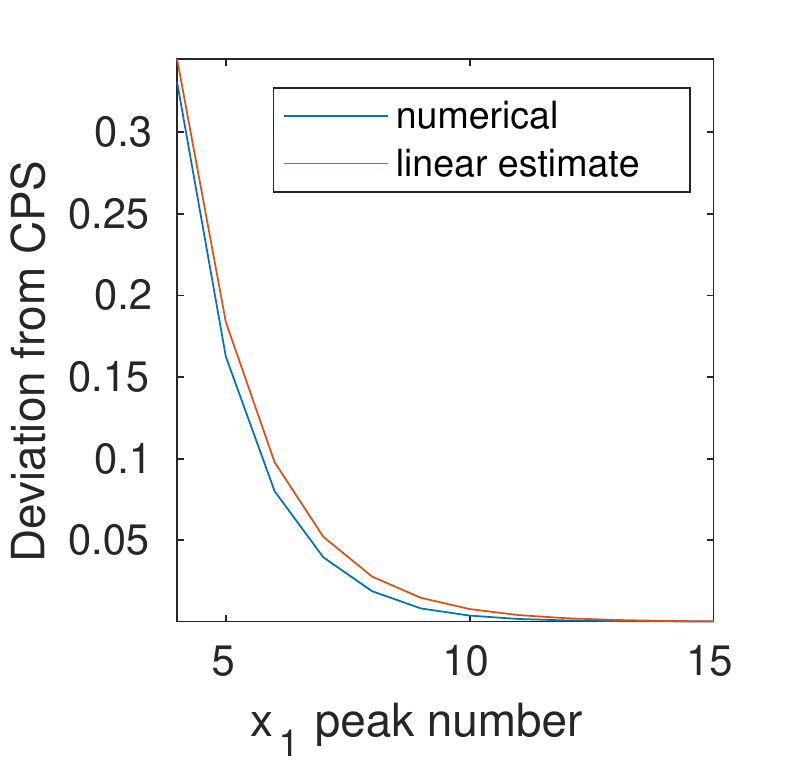}
\end{tabular}
\ece 

\vs{-4mm}
\caption{\small CP from states $(4,0)$ to $\uh_0$ for 
$\om=0.04$, yielding $\ga_2\approx0.52$ and hence 
slow convergence. 
(a) Adaptation of the truncation time $T$ during the continuation 
(initialized with $T=2T_p$). At $\al=0.3, 0.4, 0.7$ and $\al=1$ 
additional periods are added, 
while during the Newton loops $T$ only changes 
slightly. 
(b) $x_1=u_1$ and $y_1=u_3$ from the CP. 
(b) Deviation of the maxima of $u_1$ from $1$ (the maximum of $\uh_1$) 
on the CPS (blue line), and deviation predicted by the 
leading stable multiplier (red line). \label{toy_lf}}
\end{figure}

In {\tt cmds\_advanced} we 
essentially decrease $\om$ (to $0.04$) which makes the problem 
more expensive due to slow convergence to $\uh$, cf.~\reff{ga2}. 
For $\om=0.04$ we obtain $\ga_2\approx 0.52$ for 
the leading  stable multiplier. Intuitively, 
a change of $ \omega $ correspond to a rescaling of time in the 
costates by $ \frac{1}{\omega} $, i.e., small $ \omega \ll 1 $ reduces 
the speed of the costates. Then we expect that a canonical path 
spirals around the CPS several times while approaching it, and hence a long 
truncation time will be necessary for its computation. 
Figure \ref{toy_lf}(a) 
depicts typical results. We also compare the analytical 
and numerical Floquet multipliers, and generally find good 
agreement only for reasonably fine $t$--discretizations. 
In Fig.~\ref{toy_lf}(b) we compare the deviation of the 
$x_1$--maxima of $u$ from $\uh_{0,1}$ 
with the asymptotic analytical prediction, showing rather good agreement, see 
Cell 2 of {\tt cmds\_advanced}. 

\subsection{Optimal pollution mitigation}\label{hoocsec}
As an example for an OC problem with Hopf bifurcations we consider 
\begin{subequations} \label{oc2}
\hual{
    &V(v_0(\cdot))\stackrel!=\max_{\k(\cdot,\cdot)}J(v_0(\cdot),\k(\cdot,\cdot)), \qquad 
    J(v_0(\cdot),\k(\cdot,\cdot)):=\int_0^\infty\er^{-\rho t}
J_{ca}(v(t),\k(t)) \dd t,
}
with discount rate $\rho>0$, and 
where $\ds J_{ca}(v(\cdot,t),\k(\cdot,t))=\frac 1{|\Om|}
\int_\Om J_c(v(x,t),\k(x,t))\dd x$ as in \S\ref{slsec} and \S\ref{vegoc-sec} 
is the spatially averaged current value function, with here 
$J_c(v,\k)=pv_1-\beta v_2-C(\k)$ the local current value,
$C(\k)=\k+\frac 1 {2\ga} \k^2$. 
The state evolution is
\hual{
&\pa_t v_1=-\k+d_1\Delta v_1, \quad \pa_t v_2=v_1-\al(v_2)+d_2\Delta v_2, 
\label{woc} 
}
\end{subequations}
with Neumann BCs $\pa_{\bf n} v=0$ on $\pa\Om$, where 
$v_1=v_1(t,x)$ models the emissions of some firms, and 
$v_2=v_2(t,x)$ is the pollution stock, while 
the control $\k=\k(t,x)$ models the firms' abatement policies. 
In $J_c$, 
$pv_1$ and $\beta v_2$ are the firms' value 
of emissions and costs of pollution, 
and $C(\k)$ are the costs for abatement, and 
$\al(v_2)=v_2(1-v_2)$ in (\ref{oc2}b) is the recovery function of the  environment. 
Again, the $\max$ in (\ref{oc2}a) runs over all admissible controls 
$\k$, meaning that $\k\in L^\infty((0,\infty)\times\Om,\R)$, 
and we do not consider active control or state constraints. 
The associated ODE OC problem (no $x$--dependence of $v,\k$) 
was set up and analyzed in \cite{TH96,wirl00}; in suitable parameter regimes 
it shows Hopf bifurcations of periodic orbits for the associated 
canonical (ODE) system. See also, e.g., \cite{DF91, HMN92, wirl96,
KGF02, grassetal2008} 
for results about the occurrence of Hopf bifurcations 
and optimal periodic solutions in ODE OC problems. 

Setting $D=\bpm d_1&0\\0&d_2\epm$, 
$g_1(v,\k)=\bpm  -\k\\v-\al(w)\epm$, and introducing  
the co--states (Lagrange multipliers) 
$$
\lam:\Om\times(0,\infty)\ra \R^2,
$$ 
and the (local current value) Hamiltonian
$\CH=\CH(v,\lam,\k)=J_c(v,\k)+\spr{\lam,D\Delta v+g_1(v,\k)}$, 
by Pontryagins Maximum Principle we obtain 
\begin{subequations} \label{csho}
\hual{
\pa_t v&=\pa_\lam\CH=D\Delta v+g_1(v,\k), \quad v|_{t=0}=v_0, \\
\pa_t \lam&=\rho\lam-\pa_v\CH=\rho\lam+g_2(v,\lam)-D\Delta\lam, 
}
\end{subequations}
where $\pa_{\bf n} \lam=0$ on $\pa\Om$, and 
\huga{\label{kform}
\k=\k(\lam_1)=-(1+\lam_1)/\ga.
}
Finally we set 
$u(t,\cdot):=(v(t,\cdot),\lam(t,\cdot)): \Om\ra\R^{4}$,  
and write \reff{csho} as 
 \hual{\label{cs2ho}
&\pa_t u=-G(u):=\CD\Delta u+f(u),
}
where $\ds 
\CD=$diag$(d_1,d_2,-d_1,-d_2)$, 
$\ds f(u)=\biggl(-\k, v_1-\al(v_2), 
\rho\lam_1-p-\lam_2, (\rho+\al'(v_2))\lam_2+\beta\biggr)^T$.

For all parameter values, \reff{cs2ho} has the spatially homogeneous 
CSS 
$$u^*=(z_*(1-z_*),z_*,-1,-(p+\rho)), \quad \text{where}\quad 
z_*=\frac 1 2\left(1+\rho-\frac \beta{p+\rho}\right). 
$$
 We use similar parameter ranges as in \cite{wirl00}, namely 
\huga{
(p,\beta,\ga)=(1,0.2,300), \text{ and } \rho\in[0.5,0.65] \text{ as 
a continuation parameter}, 
}
consider \reff{cs2ho} over 
$\Om=(-\pi/2,\pi/2)$, 
and set the diffusion constants to $d_1=0.001, d_2=0.2$.
\brem\label{prem1}{\rm 
a) The motivation for the choice of $d_{1,2}$ is to have the first 
(for increasing $\rho$) 
Hopf bifurcation to a spatially patterned branch, and the second to  
a spatially uniform Hopf branch, because the former is 
more interesting from the PDE point of view. 
We use that the Hopf bifurcations 
for the model \reff{cs2ho} 
can be analyzed by a simple modification of 
\cite[Appendix A]{wirl00}. We find that 
for branches with spatial wave number $l\in\N$ the 
necessary condition for Hopf bifurcation, $K>0$ from \cite[(A.5)]{wirl00}, 
becomes $K=-(\al'+d_2l^2)(\rho+\al'+d_2l^2)-d_1l^2(\rho+d_1l^2)>0$. 
Since $\al'=\al'(z_*)<0$, a convenient way to first fulfill $K>0$ 
for $l=1$ is to choose $0<d_1\ll d_2<1$, such that for $l=0,1$ the factor 
$\rho+\al'+d_2l^2$ is the crucial one.

b) Even though we do not specify the units, 
$\rho\in[0.5,0.65]$ may be considered quite large, in the 
following sense. Typical periods of the CPS will be between 20 and 40, 
and, moreover, CPs starting not close to these CPSs will need times scales 
$T\ge 100$ (and larger) for convergence to the CPSs, but $\rho>0.5$ 
means that the large time ($T\ge 100$) behavior of a CP hardly plays 
a role for the value of the CP, as the discounted current value 
drops to $\er^{-\rho t}J_c(t)<\er^{-50}$. Thus, our example 
turns out to be somewhat academic, but nevertheless it will show 
the robustness of our approach. 

c) In the literature, most of the (ODE) OC examples with canonical 
periodic states 
show  these at rather large discount rates, see, e.g., \cite{HMN92, KGF02}. 
An exception is for instance the resource management model in 
\cite{BPS01}, where (ODE) CPSs are found at discount rates near $\rho=0.1$. 
We have also implemented this example, including a PDE setting, 
but its main drawback, already hinted at in \cite{BPS01}, is that 
already in the ODE setting it is 
extremely rich in CPSs, which undergo several period doubling and 
fold bifurcations, and the smaller $\rho$ is again offset by rather long 
periods (between 20 and 60). 
In summary, we focus on the pollution example because it gives a clear 
and robust bifurcation picture.  
}
\eex\erem 

\input{pollfig1}

The implementation of \reff{cs2ho} works as usual, and, moreover, 
the computation of the bifurcation diagram of CSS and CPS, 
and of the Floquet multipliers, is already explained 
in \cite{hotheo, hotut}, and the novelty here is the computation of CPs 
to the CPSs in a full PDE setting. Table \ref{polltab1} gives a few 
comments on the used files.
In Figure \ref{ocf1} (essentially already contained in \cite{hotheo}) 
we give some basic results for \reff{cs2ho} 
with a coarse spatial discretization of $\Om$ by only $n_p=21$ 
points (and thus $n_u=84$).  
(a) shows the full spectrum of the linearization 
of \reff{cs2ho} around $u^*$ at $\rho=0.5$. 
(b) shows a basic 
bifurcation diagram. At $\rho=\rho_1\approx 0.53$ there bifurcates a 
Hopf branch {\tt h1} with spatial wave number $l=1$, and at 
$\rho=\rho_2\approx 0.58$ a 
spatially homogeneous ($l=0$) Hopf branch {\tt h2} bifurcates 
subcritically with a fold at $\rho=\rho_f\approx 0.56$. 
(c) shows the pertinent time series on h2/pt17. As should be expected, 
$J_c$ is large when the pollution stock is low and emissions are high, 
and the pollution stock follows the emissions with some delay. 
In (b) we plot $J$ over 
$\rho$. For the CSS $u^*$ this is again simply 
$J(u^*)=\frac 1 \rho J_{c,a}(u^*)$, but for the 
periodic orbits we take into account the 
phase, which is free for \reff{cs2ho}. If $u_H$ is a $T_p$ periodic solution 
of \reff{cs2ho}, then, for $\phi\in [0,T_p)$, we consider 
$$
J(u_H;\phi):=\int_0^\infty \er^{-\rho t}J_{c,a}(u_H(t+\phi))\dd t
=\frac 1 {1-\er^{-\rho T_p}}\int_0^{T_p}
\er^{-\rho t}J_{c,a}(u_H(t+\phi))\dd t, 
$$
which in general may depend on the phase, 
and for {\tt h2} in (c) we plot $J(u_H;\phi)$ for $\phi=0$ (full red line) 
and $\phi=T/2$ (dashed red line).  For the 
spatially periodic branch {\tt h1}, $J_{c,a}(t)$ averages out in $x$ and 
hence $J(u_H;\phi)$ only weakly depends on $\phi$.  
Thus, we first conclude that for $\rho\in(\rho_1,\rho_f)$ the 
spatially patterned periodic orbits from {\tt h1} give the 
highest $J$, while for $\rho\ge \rho_f$ this is obtained from 
{\tt h2} with the correct phase. 
The example plots (d) at {\tt h1/pt8} illustrate the 
spatio-temporal dependence of $\k$, $v$, and $J_c$ on the 
patterned CPS. 

\taskip\begin{table}
{\small\begin{tabular}{p{0.18\tew}|p{0.8\tew}}
\hline
bdcmds,cpcmds&bifurcation diagram of CSS and CPS, and 
computation of some CPs\\
cpplot, polldiagn&plot of CPs, computation and plotting of diagnostics for CPs.  \\
\end{tabular}}
\caption{{\small Main scripts and functions in {\tt ocdemos/pollution}. Additionally, we locally modify some standard \pdep\ functions for convenience 
(plotting).}\label{polltab1}}
\end{table}\teskip

It remains to 
\bci
\item compute the defects $d(u^*)$ of the CSS and 
$d(u_H)$ of periodic orbits on the bifurcating branches, 
\item compute CPs to saddle point CSSs and CPSs. 
\eci 
For $d(u^*)$ we find 
that it starts with $0$ at $\rho=0.5$, and, as expected, 
increases by 2 at each Hopf point. Below we shall focus on CPs to the 
CPSs {\tt h1/pt8} and {\tt h2/pt17}, 
and in Fig.~\ref{ocf1}(e--g) 
we illustrate typical multiplier spectra, 
computed with {\tt pqzschur}, which yields 
$|\ga_1-1|<10^{-8}$ for all computations, i.e., a 
very accurate trivial multiplier, and hence we trust it. 
The large multipliers are {\em very} large, i.e., $10^{40}$ and larger, 
even for the coarse space discretization, as should be expected 
from the spectrum in Fig.~\ref{ocf1}(a). 

On {\tt h1} we find $d(u_H)=0$ up to {\tt pt9}, 
see (e) for the $n_u/2$ smallest multipliers at {\tt pt8}, and (f) 
for $|\ga_j|$ for the large ones, which are mostly real. 
For larger $\rho$ the {\tt h1} branch looses stability by a (second) 
multiplier going through 1, and in fact at 
{\tt h1/pt8} we have $\ga_2\approx 0.948$, which suggests a slow 
convergence of CPs to the CPS. 
On {\tt h2} we start with $d(u_H)=3$, but $d(u_H)=0$ after the fold 
until $\rho=\rho_1\approx 0.6$, after which $d(u_H)$ increases 
again by multipliers going though 1. At {\tt pt17} we have 
$\ga_2\approx 0.905$, again suggesting slow convergence. 

\begin{figure}[ht]
\bce{\small
\begin{tabular}{llll}
(a)&(b)&(c)&(d)\\
\hs{-2mm}\ig[width=0.24\tew]{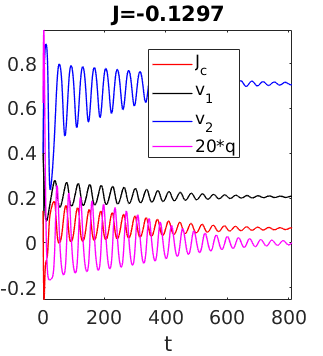}
\hs{-2mm}\ig[width=0.13\tew]{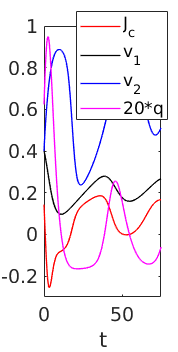}&
\hs{-2mm}\ig[width=0.125\tew]{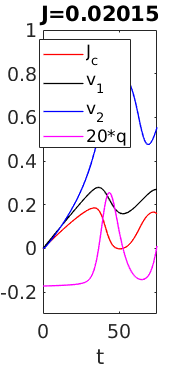}&
\hs{-2mm}\ig[width=0.24\tew]{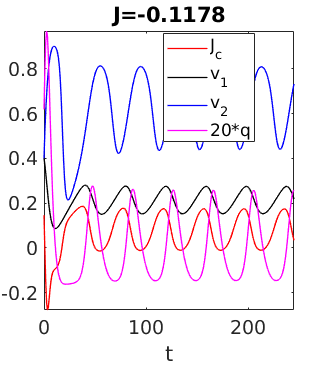}&
\hs{-2mm}\ig[width=0.238\tew]{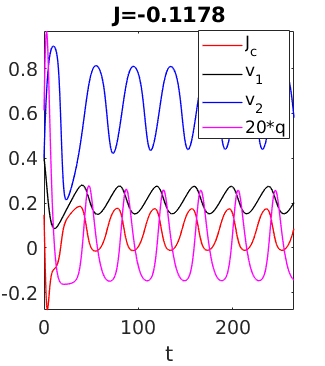}
\end{tabular}\\
\begin{tabular}{l}
\hs{-4mm}\mbox{(e)}\\
\hs{-2mm}\ig[width=0.34\tew]{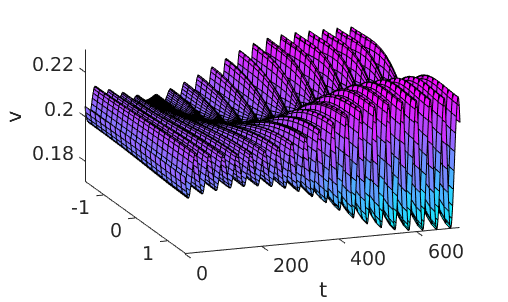}
\hs{-2mm}\ig[width=0.34\tew]{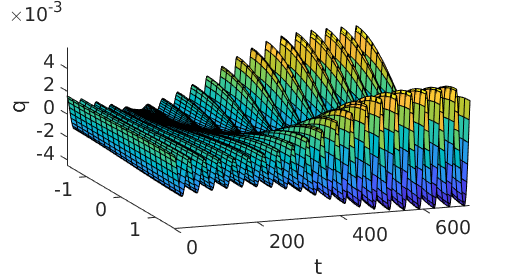}
\hs{-2mm}\ig[width=0.32\tew]{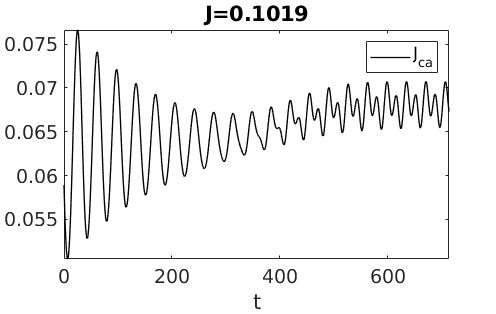}
\end{tabular}
}
\ece 

\vs{-5mm}
   \caption{{\small CPs to the CSS and to different CPS from different 
initial states. (a,b) CPs to the CSS at $\rho=0.55$, 
  ICs $v\equiv(0.4,0.4)$ in (a), with zoom of the initial phase on the right, and ICs $v\equiv(0,0)$ in (b) showing the different initial behavior, 
and hence a different value.  (c) 
CP to the homogeneous CPS at $\rho=0.57$, 
ICs $v\equiv(0.4,0.4)$. (d) Same as (c), but with 
different target $\uh_0$ on the CPS (phase shift by half a period). In (b--d) we start with a short 
initial $T=2T_p$ and set $\epsi=10^{-2}$, leading to repeated 
extension of the CPs by the CPS $\uh$ via \reff{addT}. 
(e) CP to the patterned CPS 
  {\tt h1/pt8} at $\rho=0.56$, 
starting close ($\al=0.975$) to the states 
$v\equiv (0.205,0.72)$ which is near the CSS at $\rho=0.56$. 
The patterning instability of the CSS then only manifests after a 
rather long transient. See text for further discussion. 
  \label{ocf1b}}}
\end{figure}

Nevertheless, the computation of CPs (from various initial states) 
to the CPSs works quite robustly, and in Fig.~\ref{ocf1b} 
we present some sample results from {\tt cpcmds}, and from 
{\tt pollODE/cpcmdsode}, which treats the same problem 
as an ODE. The idea is that the behavior of the spatially homogeneous 
CPs can be studied much faster in the ODE setting due 
to much less DoF. 
Also, for instance the ODE multipliers $\ga$ at {\tt h2/pt17} are 
$\ga\approx (0.303,\    1,\   1.012*10^{10},\    3.789*10^{10})$, 
i.e., $\ga_2=0.303$. 
The associated ODE paths then also exist in the PDE 
as spatially homogeneous paths, and show the same 
convergence behavior
as long as the instability which yields the existence of the 
patterned CPS {\tt h1} plays no role. 

In Fig.~\ref{ocf1b}(a)--(d) 
we show CPs to the CSS at $\rho=0.55$ (starting from two 
different initial states), and to 
the homogeneous CPS {\tt h2/pt17} at $\rho=0.57$. 
The convergence to the CSS is very slow as we are close to 
the Hopf bifurcation and hence the slowest decay rate is $\mu=0.0059$. 
For $(v_0,w_0)=(0.4,0.4)$ (significant initial 
emissions and pollution) we obtain a negative value $J\approx -0.1297$, 
as initially the emission are strongly reduced and 
the initial abatement investments are high. 
For $(v_0,w_0)=(0,0)$ (no initial emissions and pollution) we have 
$J\approx 0.0202$ due to increasing initial emission and 
negative initial abatement investments. 
To compute CPs to the CPS, we start with a rather small $T=2T_p$ 
and for $\|u(1)-\uh_0\|>\eps_\infty=10^{-2}$ use the extension of the CP 
by copies of the CPS during the continuation in $\al$. 
This leads to the extension of $T$ to about $10 T_p$, and to 
$\|u(1)-\uh_0\|_\infty\approx 10^{-4}$ for the 
final deviation from the target $\uh_0$, and we also illustrate 
how to a posteriori decrease this deviation to $10^{-6}$. 
We also run a few further tests, for instance computing CPs from the same 
ICs to shifted CPSs, e.g., {\tt h2} shifted by half a period. 
As expected, shifting the base 
point $\uh_0$ on the CPS just expands or shortens the truncation time 
$T$ by half a period, 
see Fig.~\ref{ocf1b}(d).

In  Fig.~\ref{ocf1b}(e) we give one exemplary CP to the inhomogeneous CPS at $\rho=0.56$, 
starting with ICs $0.975*(0.21,0.71)+0.025*\uh_0$. The ICs are thus 
quite close to the CSS, which is stable in the ODE, but (very weakly) 
unstable in the PDE, as we are beyond the primary Hopf bifurcation. 
Consequently, the associated CP transiently decays towards the CSS, 
before the inhomogeneous instability 
manifests and the CP converges to the inhomogeneous CPS. 
In summary, these examples show that our algorithms allow a 
robust control towards CSSs and CPSs with the SPP.

\section{Summary and outlook}\label{dsec}
We explained how to study OC problems of class \reff{oc1} 
in \pdep. The class \reff{oc1}  is quite general, 
and with the \pdep\ machinery we have a powerful tool to first study the 
bifurcations of CSS/CPS. For the computation of canonical paths to CSSs 
and CPSs, our Algorithms \ref{CSSalg} and \ref{CPSalg} implement for 
the class \reff{oc1} variants of the 
connecting orbits methods explained for ODE problems in 
\cite[Chapter 7]{grassetal2008}. For the CPS case, because of the 
very small and very large mutipliers present due diffusion and anti--diffusion, an important technical issue is the use of {\tt pqzschur} to compute 
the projection onto the center--unstable eigenspace. 
Similarly, the idea to start with a rather small truncation time $T$ 
and then using \reff{addT}, i.e., adding copies of the CPS 
to the CP to ensure convergence, seems crucial to have a fast 
and robust algorithm. 

There also is a number of issues we do not address (yet), 
for instance inequality constraints that frequently 
occur in OC problems.  In our examples we can simply check the 
natural constraints (such as $v,\k\ge 0$ in the {\tt sloc} example) 
a posteriori and find them to be always fulfilled, 
i.e., {\em inactive}. 
If such constraints become {\em active} 
the problem becomes much more complicated.

\renewcommand{\refname}{References}
\renewcommand{\arraystretch}{1.05}\renewcommand{\baselinestretch}{1}
\small
\input{octut.bbl}

\end{document}

%% file: hudef.tex
\def\bexe{\begin{exercise}}\def\eexe{\eex\end{exercise}}
\def\bsol{\begin{solution}}\def\esol{\eex\end{solution}}
\def\bexa{\begin{example}}\def\eexa{\end{example}}
\def\brem{\begin{remark}}\def\erem{\end{remark}}
\def\bthm{\begin{theorem}}\def\ethm{\end{theorem}}
\def\blem{\begin{lemma}}\def\elem{\end{lemma}}
\def\bcor{\begin{corollary}}\def\ecor{\end{corollary}}
\def\bdefi{\begin{definition}}\def\edefi{\end{definition}}
\newcommand{\IDEA}{\textbf{Idea of the Proof.} }

\newcommand{\abs}[1]{\lvert#1\rvert}

\def\bmip{\begin{minipage}{\textwidth}}\def\emip{\end{minipage}}
\def\huga#1{\begin{gather} #1 \end{gather}}

\def\hual#1{\begin{align} #1 \end{align}}

\newcommand{\R}{{\mathbb R}}
\newcommand{\C}{{\mathbb C}}\newcommand{\N}{{\mathbb N}}

\def\CK{{\cal K}}
\def\CA{{\cal A}}\def\CD{{\cal D}}  
\def\CG{{\cal G}}\def\CH{{\cal H}}\def\CL{{\cal L}}

\def\ga{\gamma}\def\om{\omega}\def\uh{\hat{u}}\def\vh{\hat{v}}
\def\noi{\noindent}\def\ds{\displaystyle}
\def\pa{{\partial}}\def\lam{\lambda}
\newcommand{\bi}{\begin{itemize}}\newcommand{\ei}{\end{itemize}}
\newcommand{\ben}{\begin{enumerate}}\newcommand{\een}{\end{enumerate}}
\newcommand{\bce}{\begin{center}}\newcommand{\ece}{\end{center}}
\newcommand{\bci}{\begin{compactitem}}\newcommand{\eci}{\end{compactitem}}
\newcommand{\bcen}{\begin{compactenum}}\newcommand{\ecen}{\end{compactenum}}
\newcommand{\reff}[1]{(\ref{#1})}
\newcommand{\ov}[1]{{\overline {#1}}}

\newcommand{\spr}[1]{\left\langle #1 \right\rangle}
\newcommand{\hs}[1]{{\hspace{#1}}}\newcommand{\vs}[1]{{\vspace{#1}}}

\def\eps{\varepsilon}\def\epsi{\epsilon}
\def\ra{\rightarrow}

\newcommand{\barr}{\begin{array}}\newcommand{\earr}{\end{array}}
\newcommand{\bpm}{\begin{pmatrix}}\newcommand{\epm}{\end{pmatrix}}
\newcommand{\bsm}{\left(\begin{smallmatrix}}
\newcommand{\esm}{\end{smallmatrix}\right)}
\newcommand{\ba}{\begin{array}}\newcommand{\ea}{\end{array}}
\def\dd{\, {\rm d}}

\def\er{{\rm e}}
\def\re{{\rm Re}}
\def\om{\omega}\def\Om{\Omega}
\def\ddt{\frac{\rm d}{{\rm d}t}}
\def\del{\delta}

\def\eex{\hfill\mbox{$\rfloor$}}

\def\sig{\sigma}
\def\al{\alpha}

\def\Lam{\Lambda}

\def\bd{\begin{displaymath}} \def\ed{\end{displaymath}}
\def\ba{\begin{array}} \def\ea{\end{array}}
  
\def\eps{\varepsilon}


%% file: p2pdef.tex
\def\pdep{{\tt pde2path}}
\def\PMAXP{Pontryagin's maximum principle}
\def\poc{{\tt p2poc}}
 
\def\oop{{\tt OOPDE}}

\def\mlab{{\tt Matlab}}

\def\dhome{/home/hannes/path/pde2path/ocdemos/}

%% file: lst.tex
\definecolor{codegreen}{rgb}{0,0.6,0}
\definecolor{codegray}{rgb}{0.5,0.5,0.5}
\definecolor{codepurple}{rgb}{0.58,0,0.82}
\definecolor{backcolour}{rgb}{0.95,0.95,0.92}
 
\lstdefinestyle{mystyle}{
    backgroundcolor=\color{backcolour},   
    commentstyle=\color{codegreen},
    keywordstyle=\color{black},
    numberstyle=\small\color{codegray},
    stringstyle=\color{codepurple},
    basicstyle=\footnotesize\ttfamily,
    breakatwhitespace=false,         
    breaklines=true,                 
    captionpos=b,                    
    keepspaces=true,                 
    numbers=left,                    
    numbersep=5pt,                  
    showspaces=false,                
    showstringspaces=false,
    showtabs=false,                  
    tabsize=2, 
  xleftmargin=4mm,
}
 

%% file: slfig2.tex
\begin{figure}[ht]
\bce{\small
\begin{tabular}{lll}
{\small (a) $P$ and $\k$ on the CP from p3/pt19 to FSC}
&{\small (b) diagnostics for (a)}&{\small (c) Using $T$ adaptation}\\
\hs{-2mm}\ig[width=0.25\tew]{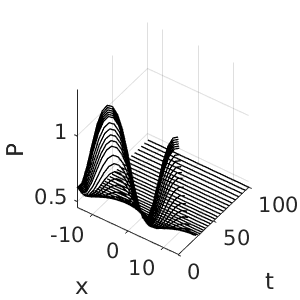}\hs{-1mm}\ig[width=0.25\tew]{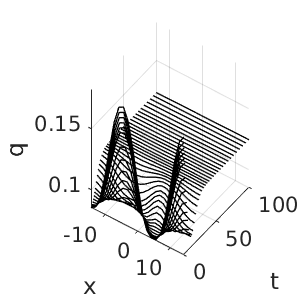}
&\raisebox{0mm}{\ig[width=0.23\tew]{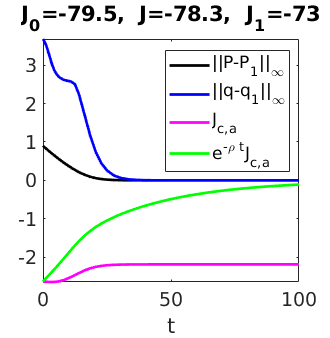}}
&{\ig[width=0.23\tew]{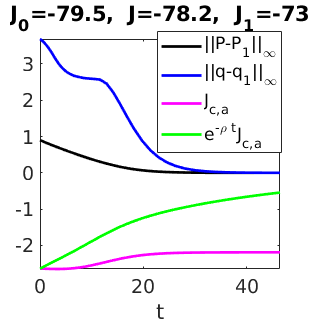}}
\\
\end{tabular}
{\small (d) Fold in $\al$ for continuation for CP from {\tt p1/pt68} to FSS, and ``upper'' canonical path at $\al=0.6$}\\
\raisebox{0mm}{\ig[width=0.26\tew]{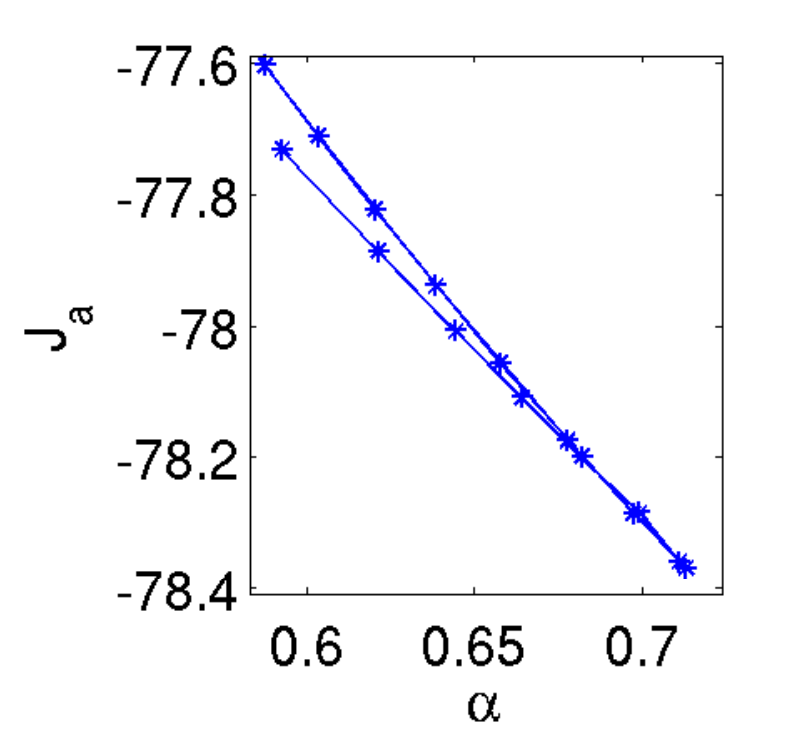}}
\ig[width=0.27\tew]{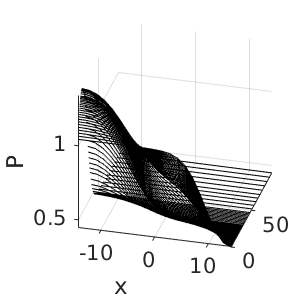}\ig[width=0.27\tew]{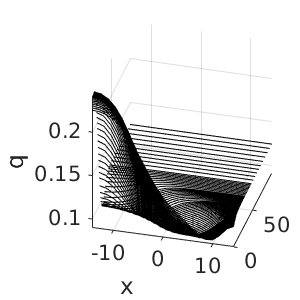}
}\ece 
\caption{Example outputs from {\tt cpdemo1D.m}. (a) shows $P,q$ 
on the CP 
from the patterned CSS p3/pt19 (see Fig.~\ref{slfig1}) to the 'Clean Flat 
State' FSC, and (b) shows typical associated 'diagnostics', namely 
the convergence behavior to $\uh$ (in $\|\cdot\|_\infty$ norm), the 
current value, and the 
discounted current value along the CP. On top we plot the value $J_0$ of the 
'starting' CSS (from which we take the states), the value $J$ of 
the path, and the value $J_1$ of the target CSS. (c) shows essentially 
the same as (b), but starting with a rather small truncation time $T$, 
and then adapted during the continuation in the initial states. 
(d) shows a case   
where a CP from some states (here taken from p1/pt68) to a target CSS 
(here again the FSC) does not seem to exist, or at least cannot 
be computed by 
continuation in $\al$, due to a fold. \label{slfig2}}
\end{figure}

%% file: octab.tex
\taskip 
\begin{table}[ht]
 {\small\begin{tabular}{p{0.15\tew}|p{0.075\tew}|p{0.7\tew}} 
struct/varname & default & description/comment \\ \hline
\textbf{{\tt oc/}}  & &struct with controls for {\tt isc} external to 
the BVP solvers {\tt mtom} and {\tt bvphdw} \\
nti && Initial number of mesh points. Highly problem dependent, 
so should be set by user. {\tt mtom} has automatic mesh refinement, so rather try a small {\tt nti}, while for {\tt bvphdw} a somewhat 
larger {\tt nti} should be used. \\
T &  & first guess for truncation time - if empty it is set by {\tt isc} \\
nTp& 2&for setting $T{=}{\tt nTp}{*}T_p$ as a guess for $T$ for CPs to a CPS, 
$T_p{=}$period of CPS\\ 
freeT &0&if 1, then truncation time is set free for CPs to CSS, and 
\reff{cssTbc} is included in the BVP, with $\eps$ in {\tt oc.tadev2}\\
tv & [] & initial t-mesh if not empty, otherwise generated by 
{\tt isc}\\ 
retsw & 0 & return-switch for isc: 0: only final soln, 1: 
solutions for all $\al$ \\ 
msw & 0 & predictor in natural continuation. 0: trivial, 1: secant \\
rhoi & 1&\pdep\ index of the discount rate $\rho$ \\
tadevs& inf & target error in sup--norm, i.e., $\eps_\infty$ in \reff{csss0}. Can be set initially, but we recommend a first 
step without $T$ adaptation (small $\al$), i.e., also {\tt freeT=0}.\\
tadev2&*&target error in euclidean norm, i.e., $\eps$ in \reff{cssTbc}, 
initialized by \reff{epsini} once \reff{csss0} is violated. Can be reset 
later to decrease the target deviation.\\
mtom & 1 & switch between {\tt mtom} and {\tt bvphdw} 
(in--house bvp solver) \\ 
sig & 0.1 & stepsize for arclength continuation \\
sigmin/sigmax & 1e-2,10 & minimal and maximal stepsizes 
for arclength continuation \\
fn&&file--names, for instance for the initial and target states of the CP\\ 
u0&&initial states; can be provided by file (see examples), or be set after 
{\tt ocinit}\\
s1&&classical \pdep\ data struct which will typically contain 
the data for the target $\uh$ (CSS in {\tt s1.u}, or CPS in {\tt s1.hopf}))\\
\hline 
\textbf{{\tt tomopt/}} & &struct with controls for 
the BVP solvers {\tt mtom} and {\tt bvphdw}\\
lu & 0 & if 0, then use \textbackslash\ (usually faster) instead of $LU$ decomposition 
in {\tt mtom}\\
* & * & Standard mtom-parameters can be given here, e.g. Itlimax, Itnlmax and Nmax for maximal number of linear and nonlinear iterations and maximal number of mesh points. See {\tt mtom} documentation. \\
tol, maxIt & 1e-8,10 & tolerance (in $\|\cdot\|_{\infty}$ norm) 
and max nr of iterations in {\tt bvphdw}.        
\end{tabular}}
\caption{{\small Data in the struct {\tt p} for 
computing CPs. 
Most parameters are set to standard values via {\tt ocinit}, but some are 
problem specific and have to be set explicitly. Others can and
some usually have to be overwritten for the specific
 problem. Some options are not listed here or in Table \ref{tab_isc}, which are for internal or expert use only. See {\tt ocinit} and {\tt isc} for  
further comments and details. \label{tab_user_ocinit}}}
\end{table}
\teskip 

%% file: pollfig1.tex
\begin{figure}[ht]
\bce{\small
\begin{tabular}{p{42mm}p{48mm}p{45mm}p{45mm}}
\hs{0mm}\mbox{(a) spectrum of} $\pa_u G(u^*)$, 
$\rho=0.5$ &(b) bifurcation diagram
&\mbox{\hs{-4mm}(c) time series on h2/pt17 (spat.~homogen.~branch)}\\
\hs{0mm}\ig[width=0.22\tew, height=35mm]{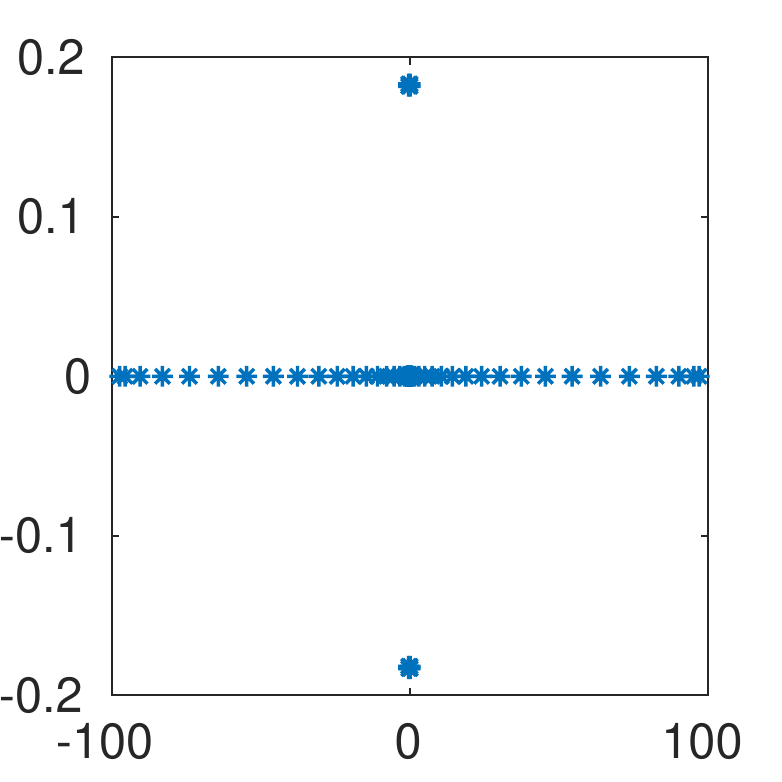}&
\hs{-4mm}\ig[width=0.22\tew]{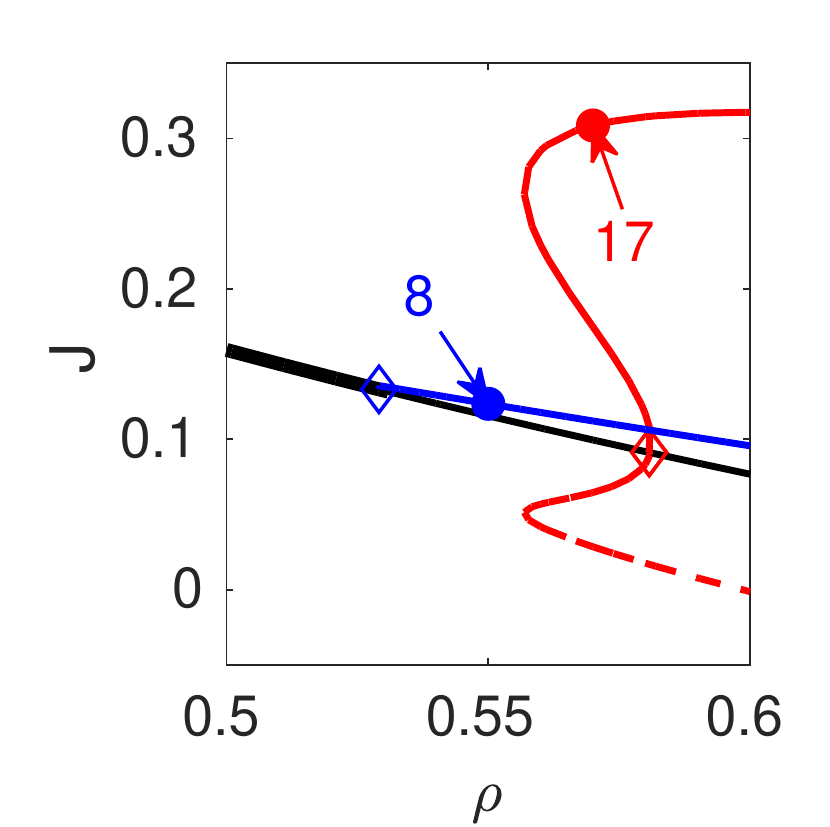}&
\hs{-8mm}\ig[width=0.22\tew]{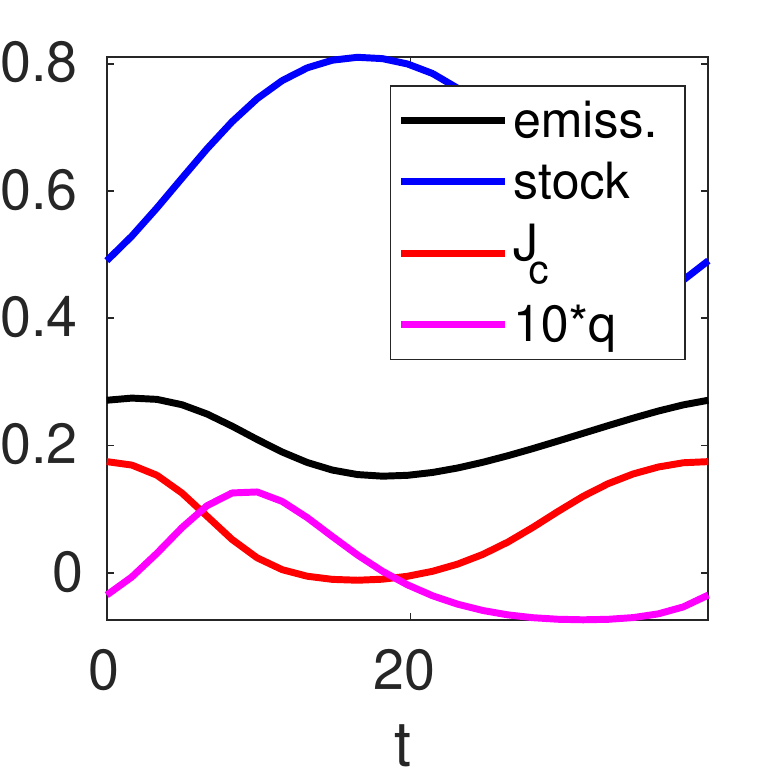}&
\hs{-8mm}\ig[width=0.22\tew]{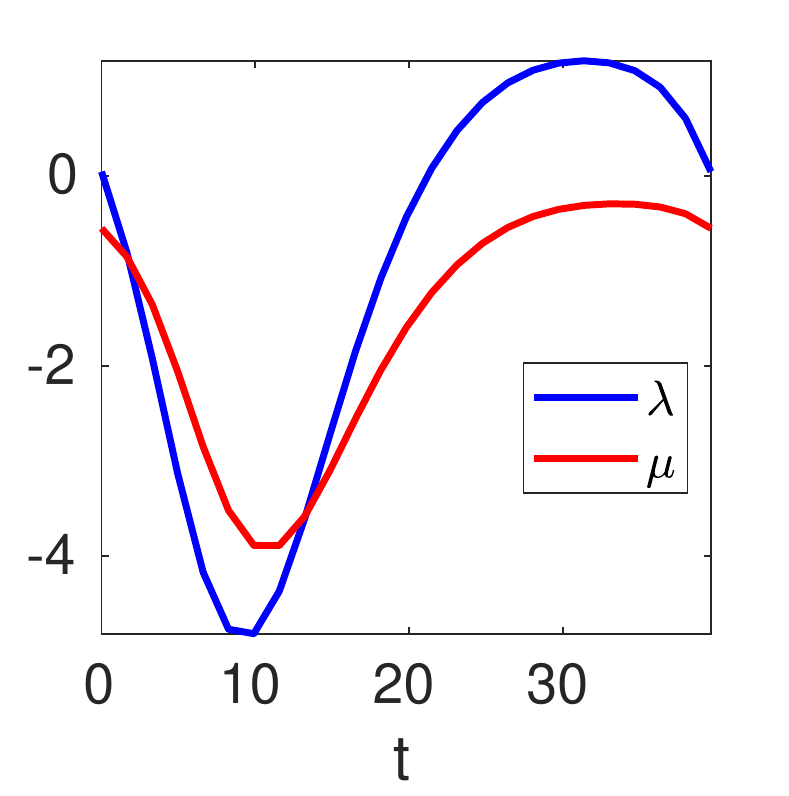}\\
\mbox{(d) sample plots at h1/pt8} \\
\ig[width=0.22\tew]{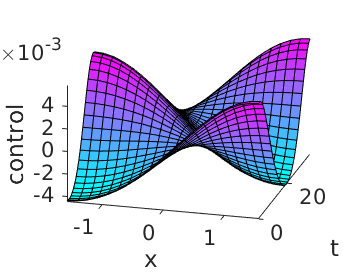}&
\hs{-4mm}\ig[width=0.22\tew]{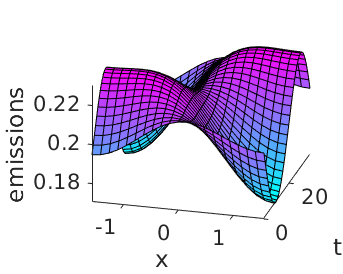}&
\hs{-10mm}\ig[width=0.22\tew]{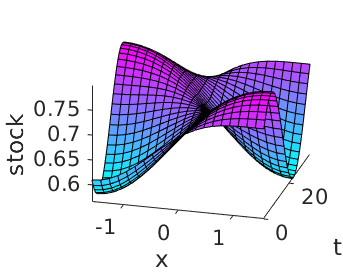}&
\hs{-10mm}\ig[width=0.22\tew]{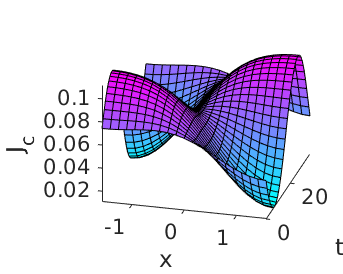}\\
\mbox{(e) the $\frac {n_u}2$ smallest $\ga_j$} at h1/pt8&
\mbox{(f) $|\ga_j|$ for the $\frac {n_u}2$} largest $\ga_j$ at h1/pt8
&\mbox{(g) the $\frac {n_u}2$ smallest $\ga_j$ at h1/pt10 and at h2/pt17}\\
\ig[height=35mm]{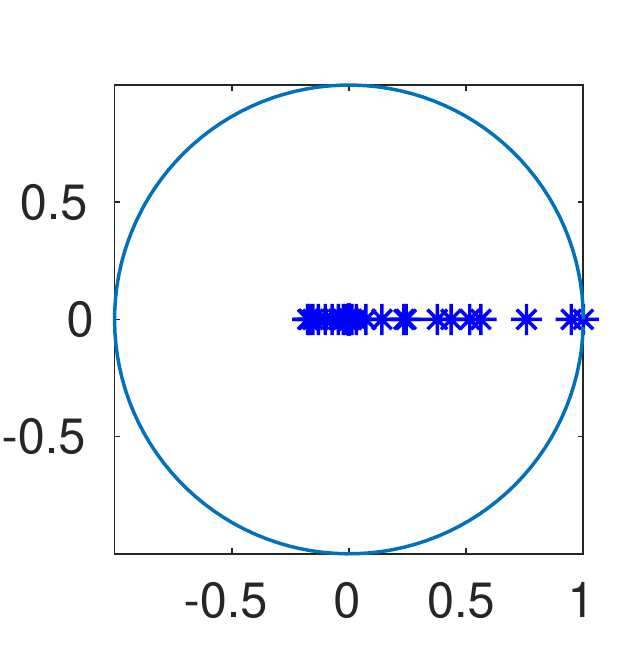}&
\ig[height=33mm]{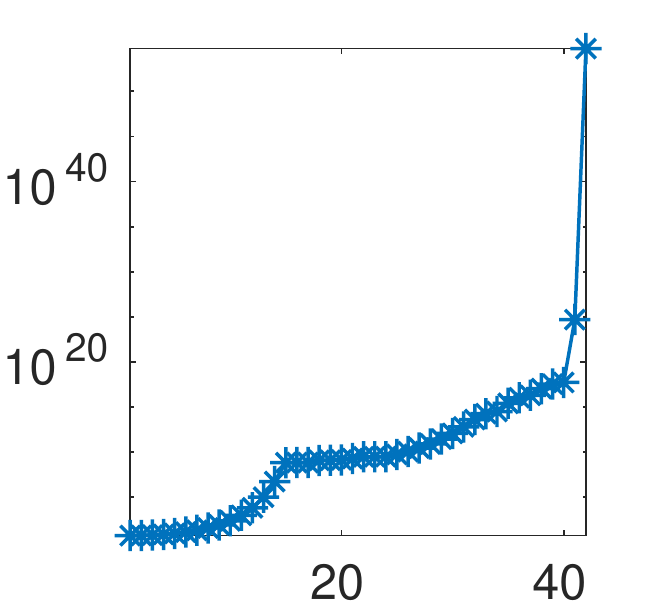}&
\mbox{\ig[height=36mm]{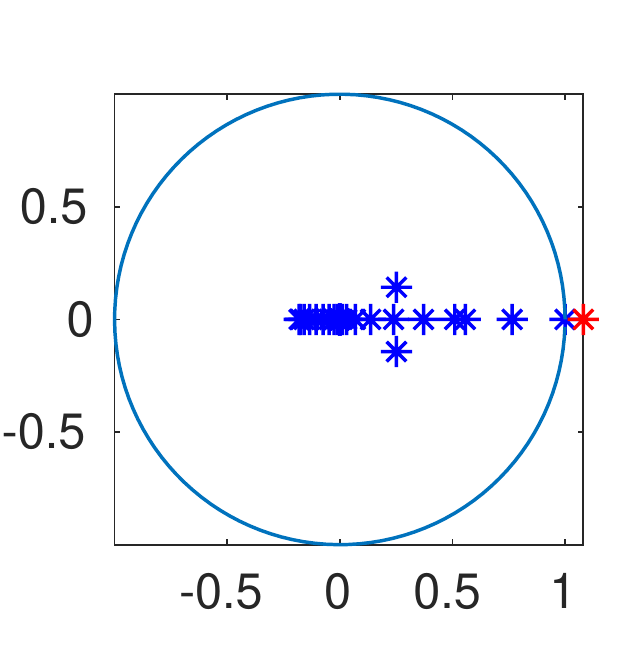}
\raisebox{3mm}{\ig[height=33mm]{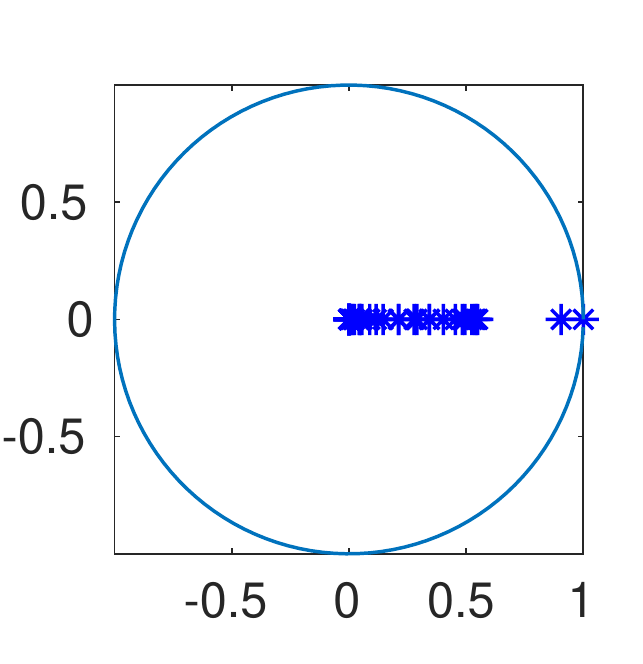}}}
\end{tabular}
}
\ece 

\vs{-5mm}
   \caption{{\small (a) full spectrum of the linearization of \reff{cs2ho} 
around $u^*$ at $\rho=0.5$ on a coarse mesh with $n_p=21$. 
(b) Bifurcation 
diagram, value $J$ over $\rho$. Black: $u^*$; blue: {\tt h1}, 
  red: {\tt h2}, $J(u^H;0)$ (full line) and $J(u^H;T/2)$ (dashed line).  
(c) Time series of a spatially homogeneous solution, including current value  
$J_c$, control $\k$, and co--states $\lam_{1,2}$. 
(d) Example plots  of $u_H$ at {\tt h1/pt8}. (e)-(g) Floquet multipliers 
at selected CPS, with $\rho=0.56$ at {\tt h1/pt10}. The largest $\ga$ in (f) 
is $\ga_{84}\approx 10^{79}$. 
  \label{ocf1}}}
\end{figure}


%% file: octut.bbl
\newcommand{\etalchar}[1]{$^{#1}$}

%% file: octut.bbl
\begin{thebibliography}{dWDR{\etalchar{+}}18}

\bibitem[AAC11]{AAC11}
S.~Ani{\c{t}}a, V.~Arn{\u{a}}utu, and V.~Capasso.
\newblock {\em An introduction to optimal control problems in life sciences and
  economics}.
\newblock Birkh\"auser/Springer, New York, 2011.

\bibitem[Bey90]{beyn90}
W.-J. Beyn.
\newblock The numerical computation of connecting orbits in dynamical systems.
\newblock {\em IMA J. Numer. Anal.}, 10(3):379--405, 1990.

\bibitem[BPS01]{BPS01}
W.J. Beyn, Th. Pampel, and W.~Semmler.
\newblock Dynamic optimization and {S}kiba sets in economic examples.
\newblock {\em Optimal Control Applications and Methods}, 22(5--6):251--280,
  2001.

\bibitem[BX08]{BX08}
W.A. Brock and A.~Xepapadeas.
\newblock Diffusion-induced instability and pattern formation in infinite
  horizon recursive optimal control.
\newblock {\em Journal of Economic Dynamics and Control}, 32(9):2745--2787,
  2008.

\bibitem[BX10]{BX10}
W.~Brock and A.~Xepapadeas.
\newblock Pattern formation, spatial externalities and regulation in coupled
  economic--ecological systems.
\newblock {\em Journal of Environmental Economics and Management},
  59(2):149--164, 2010.

\bibitem[DCF{\etalchar{+}}97]{auto}
E.~Doedel, A.~R. Champneys, Th.~F. Fairgrieve, Y.~A. Kuznetsov, Bj. Sandstede,
  and X.~Wang.
\newblock {AUTO}: Continuation and bifurcation software for ordinary
  differential equations (with {HomCont}).
\newblock \url{http://indy.cs.concordia.ca/auto/}, 1997.

\bibitem[DF91]{DF91}
E.~Dockner and G.~Feichtinger.
\newblock On the optimality of limit cycles in dynamic economic systems.
\newblock {\em Journal of Economics}, 53:31--50, 1991.

\bibitem[DKvVK08]{DKK18a}
E.~J. Doedel, B.~W. Kooi, G.~A.~K. van Voorn, and Yu.~A. Kuznetsov.
\newblock Continuation of connecting orbits in 3{D}-{ODE}s. {I}.
  {P}oint-to-cycle connections.
\newblock {\em Internat. J. Bifur. Chaos Appl. Sci. Engrg.}, 18(7):1889--1903,
  2008.

\bibitem[DKvVK09]{DKK18b}
E.~J. Doedel, B.~W. Kooi, G.~A.~K. van Voorn, and Yu.~A. Kuznetsov.
\newblock Continuation of connecting orbits in 3{D}-{ODE}s. {II}.
  {C}ycle-to-cycle connections.
\newblock {\em Internat. J. Bifur. Chaos Appl. Sci. Engrg.}, 19(1):159--169,
  2009.

\bibitem[dWDR{\etalchar{+}}18]{qsrc}
H.~de~Witt, T.~Dohnal, J.D.M. Rademacher, H.~Uecker, and D.~Wetzel.
\newblock {pde2path - Quickstart guide and reference card}, 2018.

\bibitem[GCF{\etalchar{+}}08]{grassetal2008}
D.~Grass, J.P. Caulkins, G.~Feichtinger, G.~Tragler, and D.A. Behrens.
\newblock {\em Optimal Control of Nonlinear Processes: With Applications in
  Drugs, Corruption, and Terror}.
\newblock Springer, 2008.

\bibitem[Gra15]{grass2014}
D.~Grass.
\newblock From {0D} to {1D} spatial models using {OCMat}.
\newblock Technical report, ORCOS, 2015.

\bibitem[GU17]{GU17}
D.~Grass and H.~Uecker.
\newblock Optimal management and spatial patterns in a distributed shallow lake
  model.
\newblock {\em Electr.~J.~Differential Equations}, 2017(1):1--21, 2017.

\bibitem[GUU19]{GUU19}
D.~Grass, H.~Uecker, and T.~Upmann.
\newblock Optimal fishery with coastal catch.
\newblock {\em Natural Resource Modelling}, (e12235), 2019.

\bibitem[HMN92]{HMN92}
R.~F. Hartl, A.~Mehlmann, and A.~Novak.
\newblock Cycles of fear: periodic bloodsucking rates for vampires.
\newblock {\em J. Optim. Theory Appl.}, 75(3):559--568, 1992.

\bibitem[KGF{\etalchar{+}}02]{KGF02}
P.~Kort, A.~Greiner, G.~Feichtinger, J~Haunschmied, A.~Novak, and R.~Hartl.
\newblock Environmental effects of tourism industry investments: an
  inter‐temporal trade‐off.
\newblock {\em Optim. Control -- Appl. and Methods}, 23(1):1--19, 2002.

\bibitem[Kre01]{kressner01}
D.~Kressner.
\newblock An efficient and reliable implementation of the periodic qz
  algorithm.
\newblock In {\em IFAC Workshop on Periodic Control Systems}. 2001.

\bibitem[KW10]{KW10}
T.~Kiseleva and F.O.O. Wagener.
\newblock Bifurcations of optimal vector fields in the shallow lake system.
\newblock {\em Journal of Economic Dynamics and Control}, 34(5):825--843, 2010.

\bibitem[MS02]{mazS2002}
F.~Mazzia and I.~Sgura.
\newblock Numerical approximation of nonlinear {BVPs} by means of {BVMs}.
\newblock {\em Applied Numerical Mathematics}, 42(1--3):337--352, 2002.
\newblock Numerical Solution of Differential and Differential-Algebraic
  Equations, 4-9 September 2000, Halle, Germany.

\bibitem[MST09]{MT2009}
F.~Mazzia, A.~Sestini, and D.~Trigiante.
\newblock The continuous extension of the {B}-spline linear multistep methods
  for {BVPs} on non-uniform meshes.
\newblock {\em Applied Numerical Mathematics}, 59(3--4):723--738, 2009.

\bibitem[MT04]{MT04}
F.~Mazzia and D.~Trigiante.
\newblock A hybrid mesh selection strategy based on conditioning for boundary
  value {ODE} problems.
\newblock {\em Numerical Algorithms}, 36(2):169--187, 2004.

\bibitem[Pam01]{Pam01}
Th. Pampel.
\newblock Numerical approximation of connecting orbits with asymptotic rate.
\newblock {\em Numerische Mathematik}, 90(2):309--348, 2001.

\bibitem[RU18]{actut}
J.D.M. Rademacher and H.~Uecker.
\newblock {The OOPDE setting of pde2path -- a tutorial via some Allen-Cahn
  models}, 2018.

\bibitem[RZ99a]{RZ99}
J.~P. Raymond and H.~Zidani.
\newblock Hamiltonian {P}ontryagin's principles for control problems governed
  by semilinear parabolic equations.
\newblock {\em Appl. Math. Optim.}, 39(2):143--177, 1999.

\bibitem[RZ99b]{RZ99b}
J.~P. Raymond and H.~Zidani.
\newblock Pontryagin's principle for time-optimal problems.
\newblock {\em J. Optim. Theory Appl.}, 101(2):375--402, 1999.

\bibitem[Ski78]{Skiba78}
A.~K. Skiba.
\newblock Optimal growth with a convex-concave production function.
\newblock {\em Econometrica}, 46(3):527--539, 1978.

\bibitem[Tau15]{T15}
N.~Tauchnitz.
\newblock The {P}ontryagin maximum principle for nonlinear optimal control
  problems with infinite horizon.
\newblock {\em J. Optim. Theory Appl.}, 167(1):27--48, 2015.

\bibitem[Tr{\"o}10]{Tr10}
Fredi Tr{\"o}ltzsch.
\newblock {\em Optimal control of partial differential equations}, volume 112
  of {\em Graduate Studies in Mathematics}.
\newblock American Mathematical Society, Providence, RI, 2010.

\bibitem[TW96]{TH96}
O.~Tahvonen and C.~Withagen.
\newblock Optimality of irreversible pollution accumulation.
\newblock {\em Journal of Environmental Economics and Management},
  20:1775--1795, 1996.

\bibitem[Uec16]{U16}
H.~Uecker.
\newblock Optimal harvesting and spatial patterns in a semi arid vegetation
  system.
\newblock {\em Natural Resource Modelling}, 29(2):229--258, 2016.

\bibitem[Uec19a]{hotheo}
H.~Uecker.
\newblock {Hopf bifurcation and time periodic orbits with pde2path --
  algorithms and applications}.
\newblock {\em Comm.~in Comp.~Phys}, 25(3):812--852, 2019.

\bibitem[Uec19b]{hotutb}
H.~Uecker.
\newblock User guide on {H}opf bifurcation and time periodic orbits with
  pde2path, 2019.
\newblock Available at \cite{p2phome}.

\bibitem[Uec19c]{hotut}
H.~Uecker.
\newblock {User guide on {H}opf bifurcation and time periodic orbits with
  pde2path}, 2019.

\bibitem[Uec19d]{p2phome}
H.~Uecker.
\newblock {\url{www.staff.uni-oldenburg.de/hannes.uecker/pde2path}}, 2019.

\bibitem[UWR14]{p2pure}
H.~Uecker, D.~Wetzel, and J.D.M. Rademacher.
\newblock {pde2path -- a Matlab package for continuation and bifurcation in 2D
  elliptic systems}.
\newblock {\em {NMTMA}}, 7:58--106, 2014.

\bibitem[Wir96]{wirl96}
Fr. Wirl.
\newblock Pathways to {H}opf bifurcation in dynamic, continuous time
  optimization problems.
\newblock {\em Journal of Optimization Theory and Applications}, 91:299--320,
  1996.

\bibitem[Wir00]{wirl00}
Fr. Wirl.
\newblock Optimal accumulation of pollution: Existence of limit cycles for the
  social optimum and the competitive equilibrium.
\newblock {\em Journal of Economic Dynamics and Control}, 24(2):297--306, 2000.

\end{thebibliography}
